\documentclass[12pt,fleqn]{article}
\usepackage{amssymb,latexsym,amsmath}
\usepackage{graphicx}

    \advance\textheight by \topskip

    \def\ds{\displaystyle}
    \def\Re{{\rm Re \,}}
    \def\Im{{\rm Im \,}}

    \newtheorem{theorem}{Theorem}[section]
    \newtheorem{lemma}[theorem]{Lemma}
    \newtheorem{corollary}[theorem]{Corollary}
    \newtheorem{proposition}[theorem]{Proposition}
    
    \newtheorem{Definition}[theorem]{Definition}
    \newenvironment{definition}{\begin{Definition}\rm}{\end{Definition}}
    \newtheorem{Remark}[theorem]{Remark}
    \newenvironment{remark}{\begin{Remark}\rm}{\end{Remark}}
    \newtheorem{Example}[theorem]{Example}

    \newenvironment{proof}%
    {\rm \trivlist \item[\hskip \labelsep{\bf Proof. }]}%
    {\hspace*{\fill}$\Box$\endtrivlist}
    \newenvironment{varproof}%
    {\rm \trivlist \item[\hskip \labelsep{\bf Proof}]}%
    {\hspace*{\fill}$\Box$\endtrivlist}

    \newcommand\PVint{\mathop{\setbox0\hbox{$\displaystyle\intop$}%
        \hskip0.2\wd0%
        \vcenter{\hrule width0.6\wd0height0.5pt depth0.5pt}%
        \hskip-0.8\wd0%
        }\mskip-\thinmuskip\intop\nolimits}

\begin{document}

    \begin{center} \large\bf
        The Riemann--Hilbert approach to strong asymptotics for
        orthogonal polynomials on $[-1,1]$
    \end{center}

    \

    \begin{center} \large
        A.B.J. Kuijlaars\footnote{Supported by FWO research projects G.0176.02
        and G.0455.04, and by INTAS project 00-272} \\
        \normalsize \em
        Department of Mathematics, Katholieke Universiteit Leuven,
        Celestijnenlaan 200 B, 3001 Leuven, Belgium \\
        \rm arno@wis.kuleuven.ac.be \\[3ex]
        \large
        K.T-R McLaughlin\footnote{Supported by NSF grants \#DMS-9970328
        and \#DMS-0200749}\\
        \normalsize \em
        Department of Mathematics, University of North Carolina, Chapel Hill, NC 27599, U.S.A. \\
        \rm and \\
        \em
        Department of Mathematics, University of Arizona, Tucson, AZ 85721, U.S.A.\\
        \rm mcl@amath.unc.edu \\[3ex]
        \large
        W. Van Assche\footnote{Supported by FWO research projects
        G.0184.01 and G.0455.04, and by INTAS project 00-272}\\
        \normalsize \em
        Department of Mathematics, Katholieke Universiteit Leuven,
        Celestijnenlaan 200 B, 3001 Leuven, Belgium \\
        \rm walter@wis.kuleuven.ac.be \\[3ex]
        \rm and \\[3ex]
        \large
        M. Vanlessen\footnote{Research Assistant of the Fund for Scientific Research -- Flanders (Belgium)}\\
        \normalsize \em
        Department of Mathematics, Katholieke Universiteit Leuven,
        Celestijnenlaan 200 B, 3001 Leuven, Belgium \\
        \rm maarten.vanlessen@wis.kuleuven.ac.be
    \end{center}

    \ \\[1ex]

\begin{abstract}
\noindent We consider polynomials that are orthogonal on $[-1,1]$
with respect to a modified Jacobi weight $(1-x)^\alpha (1+x)^\beta
h(x)$, with $\alpha,\beta>-1$ and $h$ real analytic and stricly
positive on $[-1,1]$. We obtain full asymptotic expansions for the
monic and orthonormal polynomials outside the interval $[-1,1]$,
for the recurrence coefficients and for the leading coefficients
of the orthonormal polynomials. We also deduce asymptotic behavior
for the Hankel determinants and for the monic orthogonal
polynomials on the interval $[-1,1]$. For the asymptotic analysis
we use the steepest descent technique for Riemann--Hilbert
problems developed by Deift and Zhou, and applied to orthogonal
polynomials on the real line by Deift, Kriecherbauer, McLaughlin,
Venakides, and Zhou. In the steepest descent method we will use
the Szeg\H{o} function associated with the weight and for the
local analysis around the endpoints $\pm 1$ we use Bessel
functions of appropriate order, whereas Deift et al.\ use Airy
functions.
\end{abstract}

    \section{Introduction and statement of results}
\setcounter{equation}{0}

\subsection{Introduction}
Strong asymptotics for the classical orthogonal polynomials
(Hermite, Laguerre and Jacobi polynomials) are known from the book of Szeg\H{o}
\cite[Section 8.21]{Szego}. Szeg\H{o} himself developed the asymptotic theory
of orthogonal polynomials on the unit circle, and on the unit interval
for weights that satisfy the Szeg\H{o} condition, see \cite{Szego} or \cite{Freud}.
During the last thirty years many asymptotic results were found
for polynomials that are orthogonal with respect to various classes of weights,
see the surveys \cite{Lubinsky1,Lubinsky2,Lubinsky3,Nevai},
and the monographs \cite{LevinLubinsky} and \cite{VanAssche}.
These developments were related to and partly motivated by questions
from approximation theory, in particular Pad\'e approximation
and weighted polynomial approximation.

Recent developments in the theory of random matrices \cite{Mehta}
provided a renewed interest in orthogonal polynomials from a
different perspective. For example, when the size of the matrices
tends to infinity, local statistical properties of the eigenvalues
are related to the strong asymptotics of the associated orthogonal
polynomials. Important new techniques for obtaining asymptotics
for orthogonal polynomials in all regions of the complex plane
were subsequently developed by Bleher and Its \cite{BleherIts} and
by Deift et al.\
\cite{Deift,DKMVZ1,DKMVZ2,KriecherbauerMcLaughlin}. The asymptotic
analysis in the latter papers is based on the characterization of
the orthogonal polynomials by means of a Riemann--Hilbert problem
for $2\times 2$ matrix valued functions due to Fokas, Its, and
Kitaev \cite{FokasItsKitaev}, together with an application of the
steepest descent method of Deift and Zhou, introduced in
\cite{DeiftZhou} and further developed in
\cite{BDJ,DVZ,DeiftZhou2} and the papers cited before. See
\cite{Deift,DKMVZ3,Kuijlaars} for an introduction, and
\cite{BKMM,BleherIts2,DeiftZhou3,DeiftZhou4,ErcolaniMcLaughlin,KMM,KM,KM2,KVW,KuijlaarsVanlessen,KV2,
Miller,StrahovFyodorov,Vanlessen1,Vanlessen2} for the latest
developments. A different method to obtain strong asymptotics for
orthogonal polynomials was recently developed by Wong and
co-authors, see e.g.\ \cite{LiWong,QuiWong,WangWong}.

The orthogonal polynomials  considered in \cite{DKMVZ1,DKMVZ2} are
orthogonal with respect to an exponential weight $e^{-Q(x)}$ on $\mathbb R$ or
with respect to varying weights $e^{-nV(x)}$ on $\mathbb R$.
Inspired by these papers,  we consider here
polynomials that are orthogonal on a finite interval $[-1,1]$ with
respect to a modified Jacobi weight
    \begin{equation} \label{Definitiew}
        w(x) = (1-x)^{\alpha} (1+x)^{\beta} h(x), \qquad x \in [-1,1],
    \end{equation}
where $\alpha, \beta > -1$ and $h(x)$ is real analytic and
strictly positive on $[-1,1]$. The main difference between the
weights (\ref{Definitiew}) and the weights $e^{-Q(x)}$ and $e^{-nV(x)}$
where $Q$ is a polynomial and $V$ is real analytic on $\mathbb R$,
lies in the behavior near the fixed endpoints $\pm 1$. In the language of random matrix
theory, this is the difference between the hard edge and the soft edge
of the spectrum \cite{NagaoSlevin,TracyWidom}. While the analysis near the soft edge of the
spectrum typically involves Airy functions, near the hard edge we
expect Bessel functions.

We use $\pi_n(x) = \pi_n(x;w)$ to denote the monic polynomial of degree $n$
orthogonal with respect to the weight $w$ on $[-1,1]$.  Thus
    \[
        \int_{-1}^1 \pi_n(x;w) x^k w(x) \ dx = 0, \qquad
        \mbox{for $k = 0,1, \ldots, n-1$,}
    \]
and $p_n(x) = p_n(x;w)$ to denote the corresponding orthonormal polynomials.
Thus
    \[
        p_n(x) = \gamma_n \pi_n(x),
    \]
where $\gamma_n > 0$ is the leading coefficient of $p_n$.

The weights (\ref{Definitiew}) belong to the class of weights that
satisfy the Szeg\H{o} condition
    \begin{equation} \label{SzegoConditie}
        \int_{-1}^1 \frac{\log w(x)}{\sqrt{1-x^2}} dx > - \infty.
    \end{equation}
Szeg\H{o}'s theory, see \cite{Freud,Szego}, provides the
strong asymptotics of the polynomials $p_n$ and $\pi_n$ and the
leading coefficients $\gamma_n$ as $n \to \infty$. To state the
results, we need the function
    \begin{equation}  \label{DefinitieVarphi}
        \varphi(z) = z + (z^2-1)^{1/2}, \qquad  z\in \mathbb C \setminus [-1,1],
    \end{equation}
which is the conformal map from $\mathbb C \setminus [-1,1]$ onto
the exterior of the unit circle. In (\ref{DefinitieVarphi}) we take
that branch of $(z^2-1)^{1/2}$ which is analytic in $\mathbb C \setminus [-1,1]$
and behaves like $z$ as $z \to \infty$. Similar conventions will be used
later on. We also need the so-called
Szeg\H{o} function associated with a weight $w$ satisfying (\ref{SzegoConditie}),
which for our purposes we define as the function $D(z) = D(z;w)$ given by
    \begin{equation} \label{Szegofunctie}
        D(z)=\exp\left( \frac{(z^{2}-1)^{1/2}}{2\pi}\int_{-1}^{1}
        \frac{\log w(x)}{\sqrt{1-x^{2}}} \frac{dx}{z-x}\right),\qquad\mbox{for
        $z\in\mathbb{C}\setminus[-1,1]$.}
    \end{equation}
The function $D(z)$ is a non-zero analytic function on $\mathbb C \setminus [-1,1]$
such that
\[ D_+(x) D_-(x) = w(x), \qquad \mbox{ for a.e.\ } x \in (-1,1), \]
where $D_+(x)$ and $D_-(x)$ denote the limiting values of $D(z)$ as $z$ approaches $x$ from above and
below, respectively. The limit at infinity
\[ D_{\infty} = \lim_{z \to \infty} D(z) =
    \exp \left( \frac{1}{2\pi} \int_{-1}^1 \frac{\log w(x)}{\sqrt{1-x^2}} dx \right) \]
exists and is a positive real number.

If the weight $w$ on $[-1,1]$ satisfies the Szeg\H{o} condition,
then the following hold, see e.g.\ \cite{Szego}, \cite{Freud}:
    \begin{equation} \label{Szegogamman}
        \lim_{n \to \infty} \frac{\gamma_n}{2^n} = \frac{1}{\sqrt{\pi} D_{\infty}}
        = \frac{1}{\sqrt{\pi}} \exp\left( -\frac{1}{2\pi} \int_{-1}^1 \frac{\log w(x)}{\sqrt{1-x^2}} dx \right),
    \end{equation}
    \begin{equation} \label{Szegopn1}
        \lim_{n \to \infty} \frac{p_n(z)}{\varphi(z)^n} =
        \frac{\varphi(z)^{1/2}}{\sqrt{2\pi} (z^2-1)^{1/4} D(z)},
    \end{equation}
and
    \begin{equation} \label{Szegophin}
        \lim_{n\to \infty} \frac{2^n \pi_n(z)}{\varphi(z)^n} = \frac{D_{\infty}}{D(z)}
        \frac{\varphi(z)^{1/2}}{\sqrt{2} (z^2-1)^{1/4}}.
    \end{equation}
The limits (\ref{Szegopn1}) and (\ref{Szegophin}) hold uniformly for
$z$ in compact subsets of $\mathbb C \setminus [-1,1]$.

\begin{remark}
    Since the Szeg\H{o} function of $\sqrt{1-x^2}$ is equal to
    $(z^2-1)^{1/4}\varphi(z)^{-1/2}$, we may also write (\ref{Szegopn1}) as
    \begin{equation}\label{Szegopn2}
        \lim_{n \to \infty} \frac{p_n(z)}{\varphi(z)^n} =
        \frac{1}{\sqrt{2\pi} D(z; \sqrt{1-x^2}w(x))},
    \end{equation}
    which is the form usually found in the literature, see e.g.\
    \cite[Lemma 1.8.]{VanAssche}.
    Because of (\ref{Szegopn2}),  $D(z; \sqrt{1-x^2}w(x))$ is sometimes
    called the Szeg\H{o} function associated with $w$, instead of $D(z; w)$.
    We found it more convenient to use (\ref{Szegopn1}). This form was also
    used, for example, in  \cite{BSW1,BSW2}.
\end{remark}

\begin{remark}
    There are also results for the orthogonal polynomials
    on the interval $[-1,1]$ under conditions that are somewhat
    stronger than the Szeg\H{o} condition, see e.g.\ \cite{Freud,LevinLubinsky,Szego}.
\end{remark}

In this paper, we give more precise asymptotic results than (\ref{Szegogamman}),
(\ref{Szegopn1}), (\ref{Szegophin}) for the special weights
(\ref{Definitiew}). We are able to obtain full asymptotic
expansions for $\gamma_n$, $\pi_n$, and $p_n$, as well as for the
coefficients $a_n$ and $b_n$ in the three-term recurrence relation
    \begin{equation} \label{RecRelphin}
        \pi_{n+1}(z) = (z-b_n) \pi_n(z) - a_n^2 \pi_{n-1}(z),
    \end{equation}
satisfied by the monic orthogonal polynomials.

From our analysis we are also able to derive strong asymptotics
for the orthogonal polynomials in the open interval $(-1,1)$, as
well as near the endpoints $\pm 1$. In a follow-up paper
\cite{KuijlaarsVanlessen} we discuss its consequences for the
level spacings of eigenvalues in the unitary random matrix
ensembles associated with the weights (\ref{Definitiew}).

\medskip

Besides the mapping function $\varphi$ and the Szeg\H{o} function
$D$, the statement of the results involves certain coefficients that
are determined by the extra factor $h$ in the weight.
Since $h$ is real analytic and strictly positive on $[-1,1]$,
there is a neighborhood $U$ of
$[-1,1]$ such that $h$ has an analytic extension to $U$ (also
denoted by $h$) with positive real part. Then $\log h$ is defined
and analytic on $U$, where we take the branch of the logarithm
which is real on $[-1,1]$. For definiteness, we take $U$ of the
form
    \begin{equation} \label{DefinitieU}
        U := \{ z \in \mathbb C \mid d(z,[-1,1]) < r \}
    \end{equation}
for some $0 < r < 1$, where $d(z,[-1,1])$ denotes the distance from $z$ to
$[-1,1]$.

\begin{definition}
    Let $\gamma$ be a closed contour in $U \setminus [-1,1]$,
    encircling the interval $[-1,1]$ once in the positive direction.
    We define two sequences $(c_n)$ and $(d_n)$ of coefficients by
    \begin{equation} \label{Definitiecn}
        c_n = \frac{1}{2\pi i} \int_{\gamma} \frac{\log h(\zeta)}{(\zeta^2-1)^{1/2}}
        \frac{d\zeta}{(\zeta-1)^{n+1}},
    \end{equation}
    and
    \begin{equation} \label{Definitiedn}
        d_n = \frac{1}{2\pi i} \int_{\gamma} \frac{\log h(\zeta)}{(\zeta^2-1)^{1/2}}
        \frac{d\zeta}{(\zeta+1)^{n+1}},
    \end{equation}
for $n = 0, 1, 2, \ldots$.
\end{definition}
Note that $c_n$ and $d_n$ do not depend on the precise form of the contour $\gamma$.

\subsection{Asymptotics for monic polynomials}

Our first result concerns the monic orthogonal polynomials $\pi_n$.
It is similar to Theorem 8.21.9 of \cite{Szego} where a full asymptotic expansion
for the usual Jacobi polynomials $P_n^{(\alpha,\beta)}(z)$
(that is, $h \equiv 1$) was found.

\begin{theorem} \label{theorempin}
    For $z \in \mathbb C \setminus [-1,1]$, we have that
    $2^n \pi_n(z) \varphi(z)^{-n}$ has an asymptotic expansion in
    powers of $1/n$ of the form
    \[
        \frac{2^n \pi_n(z)}{\varphi(z)^{n}} \sim
        \frac{D_{\infty}}{D(z)} \frac{\varphi(z)^{1/2}}{\sqrt{2}(z^2-1)^{1/4}}
        \left[ 1 + \sum_{k=1}^{\infty} \frac{\Pi_k(z)}{n^k}
        \right], \qquad \mbox{as $n\to\infty$.}
    \]
    This expansion is valid uniformly on compact subsets
    of $\overline{\mathbb C} \setminus [-1,1]$. The functions $\Pi_k(z)$
    are analytic on $\mathbb C \setminus [-1,1]$, and are explicitly computable.
    The first two are
    \begin{eqnarray} \label{DefinitiePi1}
        \Pi_1(z) & = &  -\, \frac{4\alpha^2-1}{8(\varphi(z)-1)}+ \frac{4\beta^2-1}{8(\varphi(z)+1)}\ , \\[2ex]
        \nonumber
        \Pi_2(z) & = & \frac{(4 \alpha^2-1)(\alpha+\beta+c_0)}{16(\varphi(z)-1)}
            -\frac{(4\beta^2-1)(\alpha + \beta + d_0)}{16(\varphi(z)+1)} -
            \frac{(4\alpha^2-1)(4\beta^2-1)}{128(z^2-1)} \\[1ex]
            \label{DefinitiePi2}
            & & +\, \frac{2\alpha^2+2\beta^2-5}{64}
            \left[ \frac{4\alpha^2-1}{(\varphi(z)-1)^2} +
            \frac{4\beta^2-1}{(\varphi(z)+1)^2} \right].
    \end{eqnarray}
\end{theorem}

\begin{remark}
    The functions $\Pi_k(z)$ are all rational function in $\varphi(z)$,
    that vanish at $\infty$. With increasing $k$, the functions $\Pi_k(z)$
    are increasingly more complicated.

    If $4\alpha^2 - 1 = 4 \beta^2 - 1 = 0$, then it turns out that
    all functions $\Pi_k(z)$ are identically zero. In that case one can in fact prove that
    \[
        \frac{2^n \pi_n(z)}{\varphi(z)^{n}} =
        \frac{D_{\infty}}{D(z)} \frac{\varphi(z)^{1/2}}{\sqrt{2}(z^2-1)^{1/4}}
        \left[ 1 + O\left(e^{-cn}\right) \right], \qquad \mbox{as $n\to\infty$}
    \]
    for some $c > 0$, cf.\ \cite{Kuijlaars}.
\end{remark}

\subsection{Asymptotics for leading coefficients}
The leading coefficients $\gamma_n$ of the orthonormal polynomials $p_n$
also have an asymptotic expansion in powers of $1/n$:
\begin{theorem} \label{theoremgamman}
    We have the following asymptotic expansion for $2^{-n} \gamma_n$
    \begin{equation} \label{AsympOntGamma}
        \frac{\gamma_n}{2^n} \sim  \frac{1}{\sqrt{\pi} D_{\infty}}
        \left[ 1 + \sum_{k=1}^{\infty}
        \frac{\Gamma_k}{n^k}\right], \qquad \mbox{as $n\to\infty$,}
    \end{equation}
    where the coefficients $\Gamma_k$ are explicitly computable.
    The first two are
    \begin{eqnarray} \label{DefinitieGamma1}
        \Gamma_1 &=& -\, \frac{4\alpha^2-1}{16} - \frac{4\beta^2-1}{16}\ , \\[2ex] \nonumber
        \Gamma_2 &=&  \frac{4\alpha^2-1}{32}\ (\alpha+\beta+c_0) +
        \frac{4\beta^2-1}{32}\ (\alpha+\beta+d_0) \\[1ex]
        \label{DefinitieGamma2}
        & & +\, \left(2\alpha^2+2\beta^2+7\right)\left[\frac{4\alpha^2-1}{256}+\frac{4\beta^2-1}{256}\right].
    \end{eqnarray}
\end{theorem}

The coefficients $\Gamma_k$ are polynomial expressions in
$\alpha$, $\beta$, and the coefficients $c_n$ and $d_n$, defined
in (\ref{Definitiecn}) and (\ref{Definitiedn}). It turns out that
$\Gamma_{2k}$ and $\Gamma_{2k+1}$ depend only on
\[
    \alpha, \beta, c_0, d_0, \ldots ,c_{k -1}, d_{k-1}.
\]

\subsection{Asymptotics for orthonormal polynomials}

Since $p_n = \gamma_n \pi_n$, we obtain from Theorem \ref{theoremgamman} and
Theorem \ref{theorempin}:
\begin{corollary} \label{corollarypn}
    For $z \in \mathbb C \setminus [-1,1]$, we have that
    $p_n(z) \varphi(z)^{-n}$ has an asymptotic expansion in
    powers of $1/n$ of the form
    \[ \frac{p_n(z)}{\varphi(z)^n} \sim \frac{\varphi(z)^{1/2}}{\sqrt{2\pi} (z^2-1)^{1/4} D(z)}
    \left[ 1 + \sum_{k=1}^{\infty} \frac{P_k(z)}{n^k} \right], \qquad \mbox{as $n\to\infty$.} \]
    The expansion is valid uniformly on compact subsets of $\overline{\mathbb C} \setminus [-1,1]$.
    The functions $P_k(z)$ are analytic on $\mathbb C \setminus [-1,1]$,
    and are explicitly computable. The first two are
    \begin{eqnarray}
        \nonumber
        P_1(z) & = & -\, \frac{4 \alpha^2-1}{16} \frac{\varphi(z)+1}{\varphi(z)-1}
            -\frac{4 \beta^2-1}{16} \frac{\varphi(z)-1}{\varphi(z)+1}, \\[3ex]
        \nonumber
        P_2(z) & = &
        \frac{(4\alpha^2-1)(\alpha+\beta+c_0)}{32}\frac{\varphi(z)+1}{\varphi(z)-1}
        +\frac{(4\beta^2-1)(\alpha+\beta+d_0)}{32}\frac{\varphi(z)-1}{\varphi(z)+1} \\[2ex]
        \nonumber
        & &
        + \, \frac{(4\alpha^2-1)^2}{128(\varphi(z)-1)}+\frac{(4\beta^2-1)^2}{128(\varphi(z)+1)}
        - \frac{(4\alpha^2-1)(4\beta^2-1)}{64}\frac{\varphi^2(z)+1}{(\varphi^2(z)-1)^2} \\[2ex]
        \nonumber
        & & + \, (2\alpha^2+2\beta^2-5)\left[\frac{4\alpha^2-1}{64\left(\varphi(z)-1\right)^2}+
        \frac{4\beta^2-1}{64\left(\varphi(z)+1\right)^2}\right] \\[2ex]
        \nonumber
        & & +\, (2\alpha^2+2\beta^2+7)\left[\frac{4\alpha^2-1}{256}
        +\frac{4\beta^2-1}{256}\right].
    \end{eqnarray}
\end{corollary}

\subsection{Asymptotics for Hankel determinants}

The Hankel determinant $D_n$ associated with $w$ is
\begin{equation} \label{DefinitieHankel}
    D_n = \det \left( \int_{-1}^1 x^{j+k} w(x) dx \right)_{j,k=0,\ldots,n},
\end{equation}
see \cite[Chapter II]{Szego}.
Since $\gamma_n = \sqrt{\frac{D_{n-1}}{D_n}}$ \cite[Section 2.2]{Szego}, we have
that
\begin{equation} \label{Dningammaj}
    \frac{D_n}{D_0} = \prod_{j=1}^n \gamma_j^{-2}.
\end{equation}
From Theorem \ref{theoremgamman} we then obtain the following asymptotic
formula for $D_n$ as $n\to \infty$.
\begin{corollary}
There is a constant $C > 0$ such that
\begin{equation} \label{AsympOntDn}
    D_n = C \left( \frac{1}{2} \right)^{n^2} \left( \frac{\pi D_{\infty}^2}{2} \right)^n
    n^{\frac{4\alpha^2 -1}{8} + \frac{4\beta^2-1}{8}}
    \left(1 + O\left(\frac{1}{n}\right)\right)
\end{equation}
as $n \to \infty$.
\end{corollary}
\begin{proof}
Re-write (\ref{Dningammaj}) as
\begin{equation} \label{Dnproduct}
    \frac{D_n}{D_0} = \prod_{j=1}^{n} \left(  \frac{2^{j} }{ \sqrt{\pi} D_{\infty} } \right)^{-2}
\times
    \prod_{j=1}^{n}  \left( 1 + \frac{ \Gamma_{1} }{j}  \right)^{-2}
\times
    \prod_{j=1}^{n}   \left(  \frac{\sqrt{\pi} D_{\infty} }{2^j}
 \frac{\gamma_{j}}{1 + \frac{\Gamma_{1}}{j} }  \right)^{-2}
\end{equation}
The first two factors are explicit. The first one is equal to
\begin{equation} \label{factor1}
\left( \frac{1}{2} \right)^{n^2} \left( \frac{\pi D_{\infty}^2}{2} \right)^n
\end{equation}
and the second one is equal to
\begin{equation} \label{factor2}
    C_1 n^{-2\Gamma_1}  \left(1 + O\left(\frac{1}{n}\right)\right).
\end{equation}
The third factor is the partial product of a convergent infinite product,
since by (\ref{AsympOntGamma}) we have
\[
\frac{\sqrt{\pi} D_{\infty} }{2^j} \frac{\gamma_{j}}{1 + \frac{\Gamma_{1}}{j} }
= \frac{1 + \frac{\Gamma_1}{j} + O(\frac{1}{j^2})}{1 + \frac{\Gamma_1}{j}}
=  1 + O \left(
\frac{1}{j^{2}} \right) \qquad \mbox{ as } j \to \infty. \]
If $C_2$ is the value of the infinite product, then the third factor
in (\ref{Dnproduct}) is
\begin{equation} \label{factor3}
    C_2  \prod_{j=n+1}^{\infty}   \left(  \frac{\sqrt{\pi} D_{\infty} }{2^j}
 \frac{\gamma_{j}}{1 + \frac{\Gamma_{1}}{j} }  \right)^{2}
 = C_2 \left( 1 + O\left(\frac{1}{n}\right)\right).
\end{equation}
Inserting (\ref{factor1}), (\ref{factor2}), and (\ref{factor3}) into
(\ref{Dnproduct}), and using (\ref{DefinitieGamma1}),
we obtain (\ref{AsympOntDn}) with constant $C = C_1 C_2 D(0)$.
\end{proof}

\begin{remark}
The $O(1/n)$ term in (\ref{factor2}) and (\ref{factor3}) can be expanded into
a complete asymptotic expansion, with explicitly computable constants
depending only on the coefficients $\Gamma_j$ in the expansion of $\gamma_n$.
On the other hand, we cannot evaluate explicitly the constant $C$ with our methods.
For the determination of $C$ in certain cases, see \cite{Johansson}
and the references cited therein. This is known as Szeg\H{o}'s strong limit
for Hankel determinants.
\end{remark}
\subsection{Asymptotics for recurrence coefficients}

The monic orthogonal polynomials satisfy a three term recurrence
relation of the form
\[
    \pi_{n+1}(z) = (z-b_n) \pi_n(z) - a_n^2 \pi_{n-1}(z).
\]
In our next result we give an asymptotic expansion of these
recurrence coefficients.
\begin{theorem} \label{theoremanbn}
    We have the following asymptotic expansion for $a_n$ and $b_n$
    \begin{equation} \label{asympan}
    a_n \sim \frac{1}{2} + \sum_{k=2}^{\infty} \frac{A_k}{n^k}\ , \qquad \mbox{as } n\to\infty,
    \end{equation}
    \begin{equation} \label{asympbn}
    b_n \sim \sum_{k=2}^{\infty} \frac{B_k}{n^k}\ , \qquad \mbox{as } n\to\infty.
    \end{equation}
    The coefficients $A_k$ and $B_k$ are explicitly computable. The
    first few are
    \begin{eqnarray}
        A_2 & = & -\, \frac{4\alpha^2-1}{32}-\frac{4\beta^2-1}{32}, \\[2ex]
        A_3 & = & \frac{4\alpha^2-1}{32}\ (\alpha+\beta+c_0)+\frac{4\beta^2-1}{32}\ (\alpha+\beta+d_0), \\[2ex]
        \nonumber
        A_4 & = &
        -\,
        (3\alpha^2+3\beta^2+6\alpha\beta+1)\left[\frac{4\alpha^2-1}{128}+\frac{4\beta^2-1}{128}\right]
        -\frac{(4\alpha^2-1)(4\beta^2-1)}{256}\\[1ex]
        & & -\, \frac{4\alpha^2-1}{128}\ 3c_0(2\alpha+2\beta+c_0)
                - \frac{4\beta^2-1}{128}\ 3d_0(2\alpha+2\beta+d_0),
    \end{eqnarray}
    and
    \begin{eqnarray}
        B_2 & = & \frac{\beta^2-\alpha^2}{4}, \\[2ex]
        B_3 & = & -\, \frac{\beta^2-\alpha^2}{4}\ (1+\alpha+\beta)
        +c_0\ \frac{4\alpha^2-1}{16}
            -d_0\ \frac{4\beta^2-1}{16}, \\[2ex] \nonumber
        B_4 & = & -\frac{4\alpha^2-1}{64}\ 3c_0(2\alpha+2\beta+2+c_0) + \frac{4\beta^2-1}{64}\
            3d_0(2\alpha+2\beta+2+d_0) \\[1ex]
            & & +\, \frac{\beta^2-\alpha^2}{16}\left(6\alpha+6\beta+3\alpha^2+3\beta^2+6\alpha\beta+4\right).
    \end{eqnarray}
\end{theorem}
\begin{remark}
Note that $A_1 = B_1 = 0$. The coefficients $A_2$ and $B_2$ only depend on $\alpha$
and $\beta$. This means that up to order $1/n^2$ the recurrence coefficients agree
with the recurrence coefficients of the pure Jacobi weight
$(1-x)^{\alpha} (1+x)^{\beta}$ as $n \to \infty$. The effect of the extra factor
$h$ in $w$ starts playing a role in the $1/n^3$-term by means of the coefficients
$c_0$ and $d_0$.

It turns out that for the case $4 \alpha^2-1 = 4 \beta^2 -1 = 0$,
all coefficients $A_k$ and $B_k$ vanish. In that case one can
prove that $a_n = \frac{1}{2} + O(e^{-cn})$ and $b_n = O(e^{-cn})$
for some $c >0$, cf.\ \cite{Kuijlaars}.
\end{remark}

\subsection{Asymptotics for monic polynomials on $[-1,1]$}

Our final results concern strong asymptotics for $\pi_n$ on the
interval of orthogonality $[-1,1]$. The asymptotic behavior of
$\pi_n$ in $(-1,1)$ is given by the following theorem.

\begin{theorem} \label{theorem: monic orthogonal polynomial bulk: MJW}
    For $x\in (-1,1)$
    \begin{eqnarray}
        \nonumber
        \lefteqn{
        \pi_{n}(x)
        =
        \frac{\sqrt{2}D_{\infty}}{2^n\sqrt{w(x)}(1-x^{2})^{1/4}}}\\[2ex]
        \nonumber
        && \times\, \left[\Bigl(1+O(1/n)\Bigr)
        \cos\left((n+1/2)\arccos
        x+\psi(x)-\frac{\pi}{4}\right)
        \right. \\[2ex]
        \label{theorem: monic orthogonal polynomial bulk: MJW: equation}
        & &  \qquad\qquad\qquad \left. + O(1/n) \cos\left((n-1/2)\arccos
        x+\psi(x)-\frac{\pi}{4}\right)\right],
    \end{eqnarray}
    as $n\to\infty$. The error terms hold uniformly for
    $x$ in compact subsets of $(-1,1)$ and have a full asymptotic expansion in terms of $1/n$, which can be
    calculated explicitly. In {\rm (\ref{theorem: monic orthogonal polynomial bulk: MJW: equation})} the
    phase function $\psi$ is given by
    \begin{equation}\label{definition: psi: MJW}
        \psi(x) = \frac{1}{2}\Bigl(\alpha(\arccos x-\pi)+\beta \arccos
        x\Bigr) + \frac{\sqrt{1-x^2}}{2\pi}\PVint_{-1}^1\frac{\log
        h(t)}{\sqrt{1-t^2}} \frac{dt}{t-x}.
    \end{equation}
    The singular integral being understood in the sense of the principal value.
\end{theorem}

The asymptotic behavior of $\pi_n$ near the endpoint 1 is given by
the following theorem.

\begin{theorem} \label{theorem: monic orthogonal polynomial endpoint: MJW}
There exists $\delta>0$ so that for $x\in (1-\delta,1)$
    \begin{eqnarray}
        \nonumber
        \lefteqn{\pi_n(x) = \frac{\sqrt{\pi}D_\infty}{2^n \sqrt{w(x)}}\,
        \frac{(n\arccos x)^{1/2}}{\left(1-x^2\right)^{1/4}}} \\[2ex]
        \nonumber
        & &
        \times\,
            \Bigl[\Bigl(1+O(1/n)\Bigr)\Bigl(\cos \zeta_1(x)J_\alpha(n\arccos x)+
                \sin \zeta_1(x)J'_\alpha(n \arccos x)\Bigr)\Bigr. \\[2ex]
        \label{theorem: monic orthogonal polynomial endpoint: MJW: equation}
        & &
        \qquad\qquad\left. +\,
            O(1/n) \Bigl(\cos \zeta_2(x)J_\alpha(n\arccos x)+
            \sin \zeta_2(x)J'_\alpha(n\arccos x)\Bigr)\right],
    \end{eqnarray}
    as $n\to\infty$. The error terms hold uniformly for
    $x\in(1-\delta, 1)$ and have a full asymptotic expansion in powers of $1/n$,
    which can be calculated explicitly. In
    {\rm (\ref{theorem: monic orthogonal polynomial endpoint: MJW: equation})}, $J_\alpha$
    is the usual Bessel function of order
    $\alpha$, and
    \begin{equation} \label{definition: zeta12: MJW}
        \zeta_{1,2}(x)\, =\, \pm \frac{1}{2}\arccos x + \psi(x)+\frac{\alpha\pi}{2},
    \end{equation}
    where the $+$ holds for $\zeta_1$, the $-$ for $\zeta_2$, and where $\psi$ is
    given by \rm{(\ref{definition: psi: MJW})}.
\end{theorem}

Near the endpoint $-1$, we can obtain similar asymptotic behavior
in terms of Bessel functions of order $\beta$.

\begin{remark}
    If $\alpha=\pm 1/2$, it turns out, after an easy calculation using
    the facts
    \[
        J_{1/2}(z)=\sqrt{\frac{2}{\pi z}}\sin z, \qquad
        J_{-1/2}(z)=\sqrt{\frac{2}{\pi z}}\cos z,
    \]
    that the asymptotic
    formula (\ref{theorem: monic orthogonal polynomial endpoint: MJW:
    equation}) is the same as formula (\ref{theorem: monic orthogonal polynomial bulk: MJW:
    equation}), so that (\ref{theorem: monic orthogonal polynomial bulk: MJW:
    equation}) holds uniformly in $[-1,1]$ for these cases. This can be explained from the fact that for these
    special choices of $\alpha$ we do not need to construct a parametrix
    near the endpoint 1, see \cite{Kuijlaars}.
\end{remark}

\begin{remark}
    The strong asymptotics (\ref{theorem: monic orthogonal polynomial endpoint: MJW:
    equation}) allow us to determine the asymptotics for
    the largest zeros of $\pi_n$. Let $1>x_{n,1}>x_{n,2}>\cdots$ denote the zeros of $\pi_n$
    numbered in
    decreasing order, and let $0<j_{\alpha,1}<j_{\alpha,2}<\cdots$
    denote the positive zeros of the Bessel function $J_\alpha$ numbered in increasing order.
    Then, from (\ref{theorem: monic orthogonal polynomial endpoint: MJW:
    equation}) and Hurwitz' theorem we can
    deduce that, for every $k=1,2,\ldots$
    \begin{equation}
        x_{n,k}=1-\frac{j_{\alpha,k}^2}{2n^2}+O(1/n^3),\qquad\mbox{as $n\to\infty$.}
    \end{equation}
    This agrees with the well-known asymptotics for the largest zeros of
    Jacobi polynomials, see e.g.\ \cite[Section 22.16]{AbramowitzStegun} or
    \cite[Section 8.1]{Szego}.
\end{remark}

\begin{remark}
    The strong asymptotic formula (\ref{theorem: monic orthogonal polynomial endpoint: MJW:
    equation}) is also related to the behavior of spacings of eigenvalues
    near $1$ (the hard edge) of random matrices from the associated unitary
    ensemble, see \cite{KuijlaarsVanlessen}, where it was rigorously shown that
    the eigenvalue correlations can be expressed in terms of the Bessel kernel
    \[ \mathbb J_{\alpha}(x,y)
        = \frac{J_{\alpha}(\sqrt{x}) \sqrt{y} J_{\alpha}'(\sqrt{v})
            - J_{\alpha}(\sqrt{y}) \sqrt{x} J_{\alpha}'(\sqrt{x})}{2(x-y)}. \]
\end{remark}

\subsection{Overview of the rest of the paper}

As indicated in the beginning of the paper, our inspiration for this work
comes from the papers \cite{DKMVZ1,DKMVZ2} by Deift, Kriecherbauer, McLaughlin,
Venakides and Zhou. As in these papers, we use the Riemann--Hilbert
approach to orthogonal polynomials, introduced by Fokas, Its and
Kitaev \cite{FokasItsKitaev}, together with the steepest descent method for
Riemann--Hilbert problems of Deift and Zhou \cite{DeiftZhou}.

In this method a number of transformations $Y \mapsto T \mapsto S
\mapsto R$ are applied to the original Riemann--Hilbert problem in
order to arrive at a Riemann-Hilbert problem for $R$, which is
normalized at infinity, and whose jump matrices are uniformly
close to the identity matrix. Then $R$ is uniformly close to the
identity matrix, which gives the leading terms in the asymptotic
expansion for $Y$, and hence for the orthogonal polynomials. Using
the fact that the jump matrices in the Riemann--Hilbert problem
for $R$ have a full asymptotic expansion in powers of $1/n$, we
arrive at a full asymptotic expansion for $R$. Tracing back the
steps $Y \mapsto T \mapsto S \mapsto R$, we obtain the full
asymptotic expansions for the orthogonal polynomials and related
quantities given in Theorems \ref{theorempin},
\ref{theoremgamman}, \ref{theoremanbn},
\ref{theorem: monic orthogonal polynomial bulk: MJW},
amd \ref{theorem: monic orthogonal polynomial endpoint: MJW}

All steps are analogous to the corresponding steps in \cite{DKMVZ2}.
The differences with \cite{DKMVZ2}  are connected with the endpoints $\pm 1$.
The main differences are:

\begin{itemize}
\item In the setup of the Riemann--Hilbert problem for $Y$ we have to be careful
with the endpoints $\pm 1$. If $\alpha < 0$ or $\beta < 0$, then the
weight is unbounded at $1$ or $-1$, and the Riemann--Hilbert problem cannot
be stated in terms of continuous boundary values. We include
a growth condition in the Riemann--Hilbert problem in order to control
the possible growth of $Y$ near the
endpoints in these cases, see (\ref{RHPYd})--(\ref{RHPYe}).
\item The first transformation $Y \mapsto T$ uses the conformal mapping
$\varphi$ from $\mathbb C \setminus [-1,1]$
onto the exterior of the unit circle. This step is simpler than the
corresponding step in \cite{DKMVZ2}, since there the transformation
is based on the construction of a $g$-function, which depends on the
weight and the degree $n$. There is also a preliminary step $Y \mapsto U$
which scales the problem considered in \cite{DKMVZ2}, and which does not
play a role here.
\item The second transformation $T \mapsto S$ involves a factorization of
the jump matrix and a contour deformation. The factorization is somewhat
different from the one in \cite{DKMVZ2}. After contour deformation it
leads to jump matrices for $S$ with $w^{-1}$ in some of the entries.
These jump matrices therefore have a non-integrable singularity in case $\alpha \geq 1$
or $\beta \geq 1$. Consequently, the Riemann--Hilbert problem for $S$
cannot be stated in $L^p$-sense for some $p > 1$. Fortunately, the
growth conditions at the endpoints still guarantee a unique solution,
as we show.
\item Before we can do the third transformation $S \mapsto R$ we have
to construct a parametrix for $S$. The parametrix away from the
endpoints depends in our case on the weight $w$. This is where the
Szeg\H{o} function for $w$ comes in. In \cite{DKMVZ2} this particular
step did not depend on $w$.
\item The construction of the parametrix near the endpoints $\pm 1$
involves (modified) Bessel functions (as expected), while in \cite{DKMVZ2}
it involves Airy functions. The construction is similar. One extra
complication is that we have to take into account the growth condition
at $\pm 1$.
\item Having the parametrix for $S$ in all regions, we can do
the third transformation $S \mapsto R$ in the same way as in \cite{DKMVZ2}.
The asymptotic expansion of $R$ follows from the known asymptotics of
the modified Bessel functions at infinity, while of course it was
based on the asymptotics of the Airy functions in \cite{DKMVZ2}.
\end{itemize}

Throughout the paper, $w$ will be the weight (\ref{Definitiew}) on
$[-1,1]$ with $\alpha, \beta > -1$. The function $h$ is positive
on $[-1,1]$, and analytic with positive real part in the domain
\[ U = \{ z \in [-1,1] \mid d(z,[-1,1]) < r \} \]
with $0 < r < 1$.

The function $\varphi$ is the conformal map (\ref{DefinitieVarphi}),
$D(z)$ is the Szeg\H{o} function (\ref{Szegofunctie}) associated with $w$,
and $D_{\infty}$ is the limit of $D(z)$ for $z \to \infty$.
Other notation will be introduced when needed.

    \section{Riemann-Hilbert Problem}
\setcounter{equation}{0}

The Riemann--Hilbert approach starts from a characterization of
the monic orthogonal polynomial $\pi_n$ with respect to a weight
$w$ on $[-1,1]$ through the following Riemann--Hilbert problem (RHP)
for a $2 \times 2$ matrix valued function $Y(z)$.

\subsubsection*{RHP for \boldmath$Y$:}
\begin{enumerate}
    \item[(a)]
        $Y(z)$ is  analytic for $z\in\mathbb C \setminus [-1,1]$.
    \item[(b)]
        $Y$ possesses continuous boundary values for $x \in (-1,1)$
        denoted by $Y_{+}(x)$ and $Y_{-}(x)$, where $Y_{+}(x)$ and $Y_{-}(x)$
        denote the limiting values of $Y(z')$ as $z'$ approaches $x$ from
        above and below, respectively, and
        \begin{equation}\label{RHPYb}
            Y_+(x) = Y_-(x)
            \begin{pmatrix}
                1 & w(x) \\
                0 & 1
            \end{pmatrix},
            \qquad\mbox{for $x \in (-1,1)$.}
        \end{equation}
    \item[(c)]
        $Y(z)$ has the following asymptotic behavior at infinity:
        \begin{equation} \label{RHPYc}
            Y(z)= \left(I+ O \left( \frac{1}{z} \right)\right)
            \begin{pmatrix}
                z^{n} & 0 \\
                0 & z^{-n}
            \end{pmatrix}, \qquad \mbox{as $z\to\infty$.}
        \end{equation}
    \item[(d)]
        $Y(z)$ has the following behavior near $z=1$:
        \begin{equation}\label{RHPYd}
            Y(z)=\left\{
            \begin{array}{cl}
                O\begin{pmatrix}
                    1 & |z-1|^{\alpha} \\
                    1 & |z-1|^{\alpha}
                \end{pmatrix},
                &\mbox{if $\alpha<0$,} \\[2ex]
                O\begin{pmatrix}
                    1 & \log|z-1| \\
                    1 & \log|z-1|
                \end{pmatrix},
                &\mbox{if $\alpha=0$,} \\[2ex]
                O\begin{pmatrix}
                    1 & 1 \\
                    1 & 1
                \end{pmatrix},
                &\mbox{if $\alpha>0$,}
            \end{array}\right.
        \end{equation}
        as $z \to 1$, $z \in \mathbb C \setminus [-1,1]$.
    \item[(e)]
        $Y(z)$ has the following behavior near  $z=-1$:
        \begin{equation} \label{RHPYe}
            Y(z)=\left\{
            \begin{array}{cl}
                O\begin{pmatrix}
                    1 & |z+1|^{\beta} \\
                    1 & |z+1|^{\beta}
                \end{pmatrix}, &\mbox{if $\beta<0$,} \\[2ex]
                O\begin{pmatrix}
                    1 & \log|z+1| \\
                    1 & \log|z+1|
                \end{pmatrix},
                &\mbox{if $\beta=0$,} \\[2ex]
                O\begin{pmatrix}
                    1 & 1 \\
                    1 & 1
                \end{pmatrix},
                &\mbox{if $\beta>0$,}
            \end{array}\right.
        \end{equation}
        as $z \to -1$, $z \in \mathbb C\setminus [-1,1]$.
\end{enumerate}

\begin{remark}
    The $O$-terms in (\ref{RHPYd}) and (\ref{RHPYe}) are to be taken
    entrywise. So for example $Y(z)= O\begin{pmatrix} 1 & |z-1|^{\alpha} \\
    1 & |z-1|^{\alpha}  \end{pmatrix}$  means that
    $Y_{11}(z) = O(1)$, $Y_{12}(z) = O(|z-1|^{\alpha})$, etc.
\end{remark}

\begin{remark}
    The RHP for orthogonal polynomials is due to Fokas, Its, and Kitaev \cite{FokasItsKitaev}.
    See Deift \cite{Deift} for an excellent exposition.
    The above formulation of the RHP for weights on $[-1,1]$ differs from the
    formulation in \cite{Deift} in two respects.
    \begin{itemize}
        \item In \cite{Deift} uniqueness of the  solution of the RHP is
            formulated in the setting of $L^2$ boundary values. Here in (b) we
            assume continuous boundary values on the open interval $(-1,1)$,
            and we assume a certain growth condition at the endpoints $\pm 1$
            in (d) and (e). If $\alpha > - \frac{1}{2}$ and $\beta > -
            \frac{1}{2}$, then these growth conditions ensure that the RHP can
            also be taken in $L^2$ sense. This is no longer true if $\alpha
            \leq - \frac{1}{2}$ or $\beta \leq - \frac{1}{2}$. We will give an
            independent proof of the uniqueness of the solution of the RHP for
            $Y$ in the first lemma.
        \item
            The conditions (d) and (e) control the behavior at the endpoints $\pm 1$.
            There is no such condition in the RHP for orthogonal polynomials on the
            full real line, which are the object of research in \cite{Deift,DKMVZ1,DKMVZ2}.
            We note that in the following we will encounter even more singular
            behavior at the endpoints.
    \end{itemize}
\end{remark}

\begin{lemma}\label{lemmaRHPY}
    If a solution of the RHP for $Y$ exists then it is unique.
\end{lemma}
\begin{proof}
    The proof is as in \cite[p.44]{Deift} except that we have to pay
    special attention to the behavior at the endpoints $\pm 1$.

    Let $Y$ be a solution. Then $\det Y$ is analytic on
    $\mathbb{C}\setminus[-1,1]$ and by the jump condition
    (\ref{RHPYb}) we have for $x\in (-1,1)$,
    \[
        (\det Y)_{+}(x)=\det\left(Y_{+}(x)\right)=\det
        \left(Y_{-}(x)\right)\det
        \begin{pmatrix}
            1 & w(x) \\
            0 & 1
        \end{pmatrix}
        =(\det Y)_{-}(x),
    \]
    so that $\det Y$ is analytic in $\mathbb C \setminus\{-1,1\}$.

    From (\ref{RHPYd}) it follows that
    \[
        \det Y(z)=
        \left\{\begin{array}{cl}
            O\left(|z-1|^{\alpha}\right),& \mbox{if $\alpha<0$,} \\[1ex]
            O\left(\log|z-1|\right), & \mbox{if $\alpha=0$,} \\[1ex]
            O(1),& \mbox{if $\alpha>0$,}
        \end{array}\right.
    \]
    as $z \to 1$. Since $\alpha>-1$, the singularity of $\det Y$ at $1$ is removable.
    Similarly it follows from (\ref{RHPYe}) that the singularity at $-1$ is removable.
    Thus $\det Y$ is an entire function.

    From (\ref{RHPYc}) it follows that $\det Y(z) \to 1$, as $z\to\infty$.
    Therefore, by Liouville's theorem, we have $\det Y(z)=1$ for every $z$.
    In particular, we see that  $Y(z)$ is invertible for every $z$ and $Y^{-1}(z)$ is
    analytic for $z\in\mathbb{C}\setminus[-1,1]$.

    \medskip

    Now suppose that $\tilde Y$ is a second solution of the RHP for $Y$.
    Then $H(z) := \tilde Y(z)Y^{-1}(z)$ is defined and analytic
    for $z\in \mathbb C \setminus[-1,1]$. As in \cite[p.44]{Deift} we have that
    $H$ is analytic across the interval $(-1,1)$. We have, since $\det
    Y=1$,
    \begin{eqnarray*}
        H(z) &= & \tilde{Y}(z)Y^{-1}(z)=
        \begin{pmatrix}
            \tilde{Y}_{11}(z) & \tilde{Y}_{12}(z) \\
            \tilde{Y}_{21}(z) & \tilde{Y}_{22}(z)
        \end{pmatrix}
        \begin{pmatrix}
            Y_{22}(z) & -Y_{12}(z) \\
            -Y_{21}(z) & Y_{11}(z)
        \end{pmatrix} \\[1ex]
        & = &
        \begin{pmatrix}
            \tilde{Y}_{11}(z)Y_{22}(z)-\tilde{Y}_{12}(z)Y_{21}(z) &
            \tilde{Y}_{12}(z)Y_{11}(z)-\tilde{Y}_{11}(z)Y_{12}(z) \\
            \tilde{Y}_{21}(z)Y_{22}(z)-\tilde{Y}_{22}(z)Y_{21}(z) &
            \tilde{Y}_{22}(z)Y_{11}(z)-\tilde{Y}_{21}(z)Y_{12}(z)
        \end{pmatrix}.
    \end{eqnarray*}
    Since $Y$ and $\tilde{Y}$ both satisfy (\ref{RHPYd}) it follows
    for each entry $H_{ij}$ of $H$ that
    \[
        H_{ij}(z)=
        \left\{\begin{array}{cl}
            O\left(|z-1|^{\alpha}\right),& \mbox{if $\alpha<0$,} \\[1ex]
            O\left(\log|z-1|\right), & \mbox{if $\alpha=0$,} \\[1ex]
            O(1),& \mbox{ if $\alpha>0$,}
        \end{array}\right.
    \]
    as $z \to 1$.  Thus the point $1$ is a removable singularity for each $H_{ij}$.
    Similarly, the point $-1$ is a removable singularity for each $H_{ij}$.
    Therefore $H$ is analytic in $\mathbb C$. Since $H(z)\to I$, as $z\to\infty$,
    it follows again by Liouville's theorem that $H=I$ and so $\tilde{Y}=Y$.

    This proves the uniqueness of the solution of the RHP for $Y$.
\end{proof}

\begin{theorem}
    The matrix valued function $Y(z)$ given by
    \begin{equation} \label{RHPYsolution}
        Y(z) =
        \begin{pmatrix}
            \pi_n(z) & \frac{1}{2\pi i} \int_{-1}^1  \frac{\pi_n(x) w(x)}{x-z}dx \\[2ex]
            -2\pi i \gamma_{n-1}^2 \pi_{n-1}(z) & -\gamma_{n-1}^2 \int_{-1}^1 \frac{\pi_{n-1}(x)w(x)}{x-z} dx
        \end{pmatrix}
    \end{equation}
    is the unique solution of the RHP for $Y$.

    Recall that $\pi_n$ is the monic polynomial of degree $n$ orthogonal with
    respect to the weight $w$ and that $\gamma_n$ is the leading coefficient of
    the corresponding orthonormal polynomial.
\end{theorem}

\begin{proof}
    It is obvious that $Y(z)$ is analytic for $z\in\mathbb C \setminus
    [-1,1]$. The proofs that $Y$ satisfies (\ref{RHPYb}) and
    (\ref{RHPYc}) are as in \cite{Deift}. We now concentrate on the
    proof of (\ref{RHPYd}). It is clear that the $(1,1)$ and $(2,1)$
    entries of $Y$ are $O(1)$ as $z\to 1$. For the $(1,2)$ entry, we have to
    prove that for $z\to 1$,
    \begin{equation} \label{Y12a}
        Y_{12}(z)=
        \left\{\begin{array}{cl}
            O\left(|z-1|^{\alpha}\right), & \mbox{if $\alpha<0$,} \\[1ex]
            O(\log|z-1|), & \mbox{if $\alpha=0$,} \\[1ex]
            O(1),  & \mbox{if $\alpha>0$.}
        \end{array}\right.
    \end{equation}
    This is obvious if $\alpha>0$, since then the limit
    \begin{equation}\label{Y12}
        \lim_{z\to 1}Y_{12}(z)=\frac{1}{2\pi i}\int_{-1}^{1}\frac{\pi_{n}(x)w(x)}{x-1}dx
    \end{equation}
    exists. Note that the integral in the right-hand side of
    (\ref{Y12}) converges if $\alpha > 0$.

    \medskip

    Consider the case $\alpha=0$. Then
    \[
        Y_{12}(z) = \frac{1}{2\pi i}\int_{-1}^{0}\frac{\pi_{n}(x)(1+x)^{\beta}h(x)}{x-z}dx
        +\frac{1}{2\pi i}\int_{0}^{1}\frac{\pi_{n}(x)(1+x)^{\beta}h(x)}{x-z}dx.
    \]
    The first term is analytic in a neighborhood of 1. For the second
    term, we note that $\pi_{n}(x)(1+x)^{\beta}h(x)$ is analytic on
    $[0,1]$, hence it satisfies the H\"older condition on the
    interval $[0,1]$ so that by \cite[p.521]{AblowitzFokas} we have
    \[
        \frac{1}{2\pi i}\int_{0}^{1}\frac{\pi_{n}(x)(1+x)^{\beta}h(x)}{x-z}dx=c\log(z-1)+\Phi_{0}(z),
    \]
    as $z\to 1$, $z\in\mathbb{C}\setminus[0,1]$, with some constant
    $c$ and a bounded function $\Phi_{0}(z)$. Thus
    $Y_{12}(z)=O(\log|z-1|)$ as $z\to 1$, in case $\alpha=0$.

    \medskip

    \begin{figure}
        \center{ \resizebox{12cm}{!}{\includegraphics{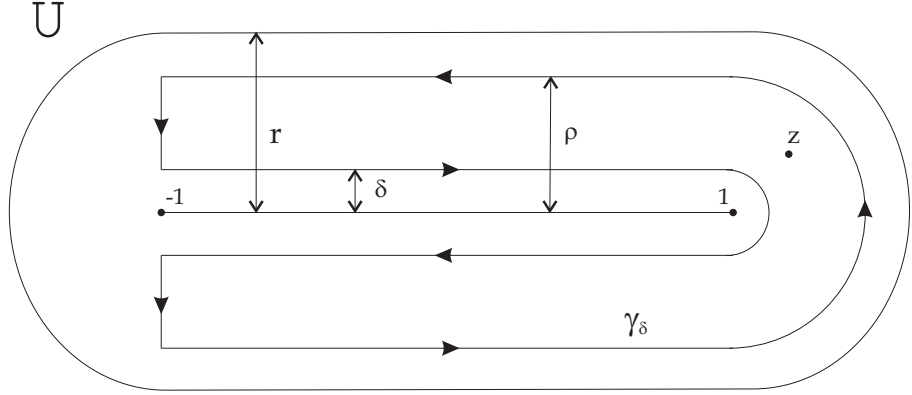}}
        \caption{The contour $\gamma_{\delta}$}\label{figuur1}}
    \end{figure}
    Finally we consider the case $\alpha<0$. Recall that the function $h$ is
    analytic in the open neighborhood $U$ of $[-1,1]$, see (\ref{DefinitieU}) and Figure \ref{figuur1}.
    Choose $0 < \rho < r$, so that the disk with radius $\rho$ and center $1$
    is contained in $U$.
    Let $|z - 1| < \rho$ with $z \not\in [-1,1]$. For $\delta>0$ sufficiently small, $z$
    lies inside the closed contour $\gamma_{\delta}$ also
    shown in  Figure \ref{figuur1}. We define the function
    $\pi_{n}(\zeta)(\zeta-1)^{\alpha}(\zeta+1)^{\beta}h(\zeta)$ with a
    cut along $(-\infty,1]$, so it is analytic in
    $U\setminus(-\infty,1]$. Therefore, by Cauchy's theorem,
    \[
        \frac{1}{2\pi i}\int_{\gamma_{\delta}}
        \frac{\pi_{n}(\zeta)(\zeta-1)^{\alpha}(\zeta+1)^{\beta}h(\zeta)}{\zeta-z}d\zeta=
        \pi_{n}(z)(z-1)^{\alpha}(z+1)^{\beta}h(z).
    \]
    Letting $\delta\to 0$, we find after an easy calculation that
    \[
        \Omega(z)+\frac{1}{2\pi i}\int_{-1}^{1}\frac{\pi_{n}(x)w(x)}{x-z}
        \left(e^{i\pi\alpha}-e^{-i\pi\alpha}\right)dx
        =\pi_{n}(z)(z-1)^{\alpha}(z+1)^{\beta}h(z),
    \]
    with $\Omega(z)$ analytic for $|z-1| < \rho$. Therefore
    \[
        Y_{12}(z)=\frac{\pi_{n}(z)(z+1)^{\beta}h(z)}{2i\sin\pi\alpha}(z-1)^{\alpha}
        -\frac{\Omega(z)}{2i\sin\pi\alpha},\qquad \mbox{ for } |z -1| < \rho, \ z \not\in [-1,1].
    \]
    Hence $Y_{12}(z)=O\left(|z-1|^{\alpha}\right)$ as  $z\to 1$, in case $\alpha < 0$.

    \medskip

    So (\ref{Y12a}) holds in all cases. Similarly it follows that
    $Y_{22}(z)$ has the correct order as $z\to 1$.
    Thus (\ref{RHPYd}) holds. The proof that (\ref{RHPYe}) holds
    follows in the same way.
\end{proof}

    \section{First transformation: $Y \mapsto T$}
\setcounter{equation}{0}

In this section we formulate an equivalent RHP. This new RHP is normalized
at infinity, and its jump matrix has rapidly oscillating entries.
Let $\sigma_3$ be the Pauli matrix
\[
    \sigma_3 =
    \begin{pmatrix}
        1 & 0 \\
        0 & -1
    \end{pmatrix}
\]
and define
\begin{equation} \label{TinY}
    T(z)=2^{n\sigma_{3}}Y(z)\varphi(z)^{-n\sigma_{3}},\qquad\mbox{for
    $z\in\mathbb{C} \setminus[-1,1]$,}
\end{equation}
with $Y$ the unique solution of the RHP for $Y$, given
by (\ref{RHPYsolution}).

Recall that $\varphi(z) = z + (z^2-1)^{1/2}$ is the conformal map
from $\mathbb C \setminus [-1,1]$ onto the exterior of the unit
disk. So it is clear from (\ref{TinY}) that $T$ is analytic in $\mathbb{C} \setminus
[-1,1]$.

We have for $x \in (-1,1)$,
\begin{equation}
\varphi_{+}(x)\varphi_{-}(x)=\left(x+i\sqrt{1-x^{2}}\right)\left(x-i\sqrt{1-x^{2}}\right)=1,
\end{equation}
and so, by (\ref{RHPYb}) and (\ref{TinY}), $T$ satisfies the
following jump relation on $(-1,1)$
\begin{eqnarray*}
    T_{+}(x) & = & 2^{n\sigma_{3}}Y_{+}(x)\varphi_{+}(x)^{-n\sigma_{3}} \\[1ex]
    & = & 2^{n\sigma_{3}}Y_{-}(x)
        \begin{pmatrix}
            1 & w(x)\\
            0 & 1
        \end{pmatrix} \varphi_{+}(x)^{-n\sigma_{3}} \\[1ex]
    & = & T_{-}(x)\varphi_{-}(x)^{n\sigma_{3}}
        \begin{pmatrix}
            1 & w(x)\\
            0 & 1
        \end{pmatrix} \varphi_{+}(x)^{-n\sigma_{3}} \\[1ex]
    & = & T_{-}(x)
        \begin{pmatrix}
            \varphi_{+}(x)^{-2n} & w(x)\\
            0 & \varphi_{-}(x)^{-2n}
        \end{pmatrix},\qquad\mbox{for $x\in(-1,1)$.}
\end{eqnarray*}
Using (\ref{RHPYc}), (\ref{TinY}), and the fact that $\varphi(z) = 2z + O(1/z)$
as $z\to\infty$, we have $T(z)=I+O(1/z)$ as $z\to\infty$. Finally,
$\varphi$ is bounded and bounded away from $0$ in a neighborhood
of $\pm 1$. Thus, $T$ has the same behavior as $Y$ around
the points $\pm 1$. So, we proved that $T$ is the solution
of the following RHP.

\subsubsection*{RHP for \boldmath$T$:}

\begin{enumerate}
    \item[(a)]
        $T(z)$ is analytic for $z\in\mathbb{C}\setminus[-1,1]$.
    \item[(b)]
        $T(z)$ satisfies the following jump relation on $(-1,1)$:
        \begin{equation}\label{RHPTb}
            T_{+}(x)=T_{-}(x)
            \begin{pmatrix}
                \varphi_{+}(x)^{-2n} & w(x) \\
                0 & \varphi_{-}(x)^{-2n}
            \end{pmatrix}, \qquad\mbox{for $x\in(-1,1)$.}
        \end{equation}
    \item[(c)]
        $T(z)$ has the following behavior at infinity:
        \begin{equation}\label{RHPTc}
            T(z) = I + O \left( \frac{1}{z} \right),
            \qquad \mbox{as $z \to \infty$.}
        \end{equation}
    \item[(d)]
        $T(z)$ has the following behavior as $z\to 1$:
        \begin{equation}\label{RHPTd}
            T(z)=\left\{
            \begin{array}{cl}
                O\begin{pmatrix}
                    1 & |z-1|^{\alpha} \\
                    1 & |z-1|^{\alpha}
                \end{pmatrix},& \mbox{if $\alpha<0$,} \\[2ex]
                O\begin{pmatrix}
                    1 & \log|z-1| \\
                    1 & \log|z-1|
                \end{pmatrix},& \mbox{if $\alpha=0$,} \\[2ex]
                O\begin{pmatrix}
                    1 & 1 \\
                    1 & 1
                \end{pmatrix}, &\mbox{if $\alpha>0$.}
            \end{array}\right.
        \end{equation}
    \item[(e)]
        $T(z)$ has the following behavior as $z\to -1$:
        \begin{equation} \label{RHPTe}
            T(z)=\left\{
            \begin{array}{cl}
                O\begin{pmatrix}
                    1 & |z+1|^{\beta} \\
                    1 & |z+1|^{\beta}
                \end{pmatrix},& \mbox{if $\beta<0$,} \\[2ex]
                O\begin{pmatrix}
                    1 & \log|z+1| \\
                    1 & \log|z+1|
                \end{pmatrix},& \mbox{if $\beta=0$,} \\[2ex]
                O\begin{pmatrix}
                    1 & 1 \\
                    1 & 1
                \end{pmatrix},& \mbox{if $\beta>0$.}
            \end{array}\right.
        \end{equation}
\end{enumerate}

The RHPs for $T$ and $Y$ are equivalent, since if $T$ is any
solution of the RHP for $T$ then it similarly follows that
\begin{equation}\label{YinT}
    Y(z) = 2^{-n\sigma_{3}}T(z)\varphi(z)^{n\sigma_{3}},
\end{equation}
solves the RHP for $Y$. Thus, by Lemma \ref{lemmaRHPY}, the RHP
for $T$ has a unique solution.

\begin{remark}
    Note that the transformation $Y \mapsto T$ has the same form
    (\ref{TinY}) for every weight $w$. What we mean is that we use the
    same function $\varphi$ for every $w$. This is in contrast to the
    case of varying weights on the real line \cite{DKMVZ1,DKMVZ2},
    where the first transformation is based on the
    construction of a so-called $g$-function, which is specific for
    the weight in question (it may even be $n$-dependent).

    In our situation the $g$-function is
    \[ g(z) = \log \varphi(z) = \log (z + (z^2-1)^{1/2}), \]
    which is defined and analytic on $\mathbb C \setminus (-\infty,1]$.
\end{remark}

\begin{remark}
    Since $\varphi_+(x)$ and $\varphi_-(x)$ have absolute value $1$
    for $x \in (-1,1)$, we see in (\ref{RHPTb}) that the diagonal
    entries in the jump matrix for $T$ have absolute value $1$, which
    are rapidly oscillating for large $n$.
\end{remark}

    \section{Second transformation: $T \mapsto S$}
\setcounter{equation}{0}

The second transformation is based on a factorization of the jump
matrix in (\ref{RHPTb}) as a product of three matrices. A simple
calculation, based on the fact that
$\varphi_{+}(x)\varphi_{-}(x)=1$ for $x\in (-1,1)$, shows that
\begin{eqnarray}\nonumber
    \lefteqn{
    \begin{pmatrix}
        \varphi_{+}(x)^{-2n} & w(x) \\
        0 & \varphi_{-}(x)^{-2n}
    \end{pmatrix}} \\[1ex]
    \label{opsplitsing*in*drie*matrices}
    & & =  \begin{pmatrix}
                1 & 0 \\
                w(x)^{-1}\varphi_{-}(x)^{-2n} & 1
            \end{pmatrix}
            \begin{pmatrix}
                0 & w(x) \\
                -w(x)^{-1} & 0
            \end{pmatrix}
            \begin{pmatrix}
                1 & 0 \\
                w(x)^{-1}\varphi_{+}(x)^{-2n} & 1
            \end{pmatrix}.
\end{eqnarray}

\begin{figure}
    \center{ \resizebox{12cm}{!}{\includegraphics{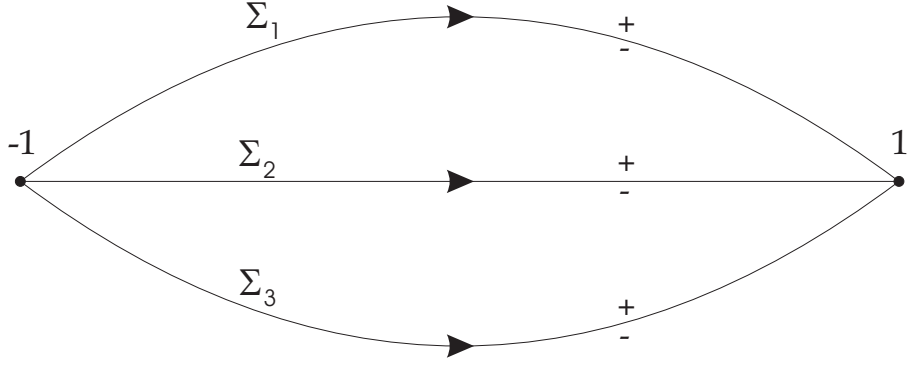}}
    \caption{The lens $\Sigma$}\label{figuur2}}
\end{figure}

As in \cite{Deift,DKMVZ1,DKMVZ2} we deform the RHP
for $T$ into a RHP with jumps on a lens shaped contour $\Sigma =
\Sigma_1 \cup \Sigma_2 \cup \Sigma_3$ as shown in Figure
\ref{figuur2}. Here $\Sigma_{1}$ and $\Sigma_{3}$ are the upper
and lower lips of the lens, respectively, and $\Sigma_{2} =
[-1,1]$. The lens is contained in $U$. We also write
\[
    \Sigma_j^o = \Sigma_j \setminus \{-1,1\}, \qquad \mbox{for $j=1,2,3$.}
\]

The weight $w$ has an analytic extension to
$U\setminus\left((-\infty,-1]\cup[1,\infty)\right)$,
also denoted by $w$, given by
\[
    w(z) = (1-z)^{\alpha} (1+z)^{\beta} h(z).
\]
The extended weight is non-zero in
$U\setminus\left((-\infty,-1]\cup[1,\infty)\right)$, since $\Re h > 0$
in $U$. Let $T$ be
the unique solution of the RHP for $T$ and define the matrix
valued function $S$ on $\mathbb{C} \setminus \Sigma$ by
\begin{equation} \label{SinT}
    S(z)=
    \left\{\begin{array}{cl}
        T(z), & \mbox{for $z$ outside the lens,} \\[2ex]
        T(z)
        \begin{pmatrix}
            1 & 0 \\
            -w(z)^{-1}\varphi(z)^{-2n} & 1
        \end{pmatrix}, & \mbox{for $z$ in the upper part of the lens,} \\[2ex]
        T(z)
        \begin{pmatrix}
            1 & 0 \\
            w(z)^{-1}\varphi(z)^{-2n} & 1
        \end{pmatrix}, & \mbox{for $z$ in the lower part of the lens.}
    \end{array} \right.
\end{equation}
Straightforward calculations then show that $S$ is a solution
of the following RHP.

\subsubsection*{RHP for \boldmath$S$:}

\begin{enumerate}
    \item[(a)]
        $S(z)$ is analytic for $z\in\mathbb{C}\setminus\Sigma$.
    \item[(b)]
        $S(z)$ satisfies the following jump relations on $\Sigma \setminus \{-1,1\}$:
        \begin{equation}\label{RHPSb1}
            S_{+}(z)=S_{-}(z)
            \begin{pmatrix}
                1 & 0 \\
                w(z)^{-1}\varphi(z)^{-2n} & 1
            \end{pmatrix}, \qquad \mbox{for $z\in\Sigma_{1}^o\cup\Sigma_{3}^o$,}
        \end{equation}
        \begin{equation}\label{RHPSb2}
            S_{+}(x)=S_{-}(x)
            \begin{pmatrix}
                0 & w(x) \\
                -w(x)^{-1} & 0
            \end{pmatrix}, \qquad \mbox{for $x\in (-1,1)$.}
        \end{equation}
    \item[(c)]
        $S(z)$ has the following behavior as $z\to\infty$:
        \begin{equation}\label{RHPSc}
            S(z)=I+O\left(\frac{1}{z}\right),\qquad\mbox{as $z\to\infty$.}
        \end{equation}
    \item[(d)]
        For $\alpha<0$, the matrix function $S(z)$ has the following behavior as $z\to 1$:
        \begin{equation}\label{RHPSd1}
            S(z)=
            O\begin{pmatrix}
                1 & |z-1|^{\alpha} \\
                1 & |z-1|^{\alpha}
            \end{pmatrix}, \qquad\mbox{as $z\to 1, z\in\mathbb{C}\setminus\Sigma$.}
        \end{equation}
        For $\alpha=0$, $S(z)$ has the following behavior as $z\to 1$:
        \begin{equation}\label{RHPSd2}
            S(z)=
            O\begin{pmatrix}
                \log|z-1| & \log|z-1| \\
                \log|z-1| & \log|z-1|
            \end{pmatrix}, \qquad \mbox{as $z\to 1, z \in \mathbb C \setminus\Sigma$.}
        \end{equation}
        For $\alpha>0$, $S(z)$ has the following behavior as $z\to 1$:
        \begin{equation}\label{RHPSd3}
            S(z)=\left\{\begin{array}{cl}
                O\begin{pmatrix}
                    1 & 1 \\
                    1 & 1
                \end{pmatrix},& \mbox{as $z\rightarrow 1$ outside the lens,} \\[2ex]
                O\begin{pmatrix}
                    |z-1|^{-\alpha} & 1 \\
                    |z-1|^{-\alpha} & 1
                \end{pmatrix}, & \mbox{as $z\to 1$ inside the lens.}
            \end{array}\right.
        \end{equation}
    \item[(e)]
        $S(z)$ has the same behavior near $-1$ if we replace in (\ref{RHPSd1}),
        (\ref{RHPSd2}), (\ref{RHPSd3}) $\alpha$ by $\beta$, $|z-1|$ by $|z+1|$ and take the
        limit $z \to -1$ instead of $z \to 1$.
\end{enumerate}

It is not completely obvious that the RHPs for $T$ and $S$ are
equivalent. We have that the (unique) solution of the RHP for $T$
gives rise to a solution of the RHP for $S$ by means of the
formulas (\ref{SinT}). Conversely, we may start from a solution
$\tilde{S}$ of the RHP for $S$ and define
\begin{equation} \label{TinS}
    \tilde{T}(z)=
    \left\{\begin{array}{cl}
        \tilde{S}(z), & \mbox{for $z$ outside the lens,} \\[2ex]
        \tilde{S}(z)
        \begin{pmatrix}
            1 & 0 \\
            w(z)^{-1}\varphi(z)^{-2n} & 1
        \end{pmatrix}, & \mbox{for $z$ in the upper part of the lens,} \\[2ex]
        \tilde{S}(z)
        \begin{pmatrix}
            1 & 0 \\
            -w(z)^{-1}\varphi(z)^{-2n} & 1
        \end{pmatrix}, & \mbox{for $z$ in the lower part of the lens.}
    \end{array}\right.
\end{equation}
We show that then $\tilde{T}$ is a solution of the RHP for $T$,
which in turn implies, since the solution of the RHP for $T$ is
unique, that the solution of the RHP for $S$ is unique as well.

\begin{lemma}\label{TheoremSecondTransform}
    Let $\tilde{S}$ be a solution of the RHP for $S$ and define
    $\tilde{T}$ by the formulas {\rm (\ref{TinS})}. Then $\tilde{T}$ is a
    solution of the RHP for $T$.
\end{lemma}

\begin{proof}
    It is straightforward to show that $\tilde{T}$ satisfies parts
    (a), (b) and (c) of the RHP for $T$. We need to show that (d) and
    (e) also hold. Let $T$ be the unique solution of the RHP for $T$.
    Then $H(z):=\tilde{T}(z) T^{-1}(z)$ is defined and analytic for
    $z\in\mathbb{C}\setminus[-1,1]$. Since
    $\tilde{T}$ and $T$ have the same jumps on $(-1,1)$, it follows
    that $H$ is analytic across $(-1,1)$. We have,
    since $\det T=1$,
    \begin{equation}\label{HifvT}
        H(z)=
        \begin{pmatrix}
            \tilde{T}_{11}(z)T_{22}(z)-\tilde{T}_{12}(z)T_{21}(z) &
            \tilde{T}_{12}(z)T_{11}(z)-\tilde{T}_{11}(z)T_{12}(z)\\[2ex]
            \tilde{T}_{21}(z)T_{22}(z)-\tilde{T}_{22}(z)T_{21}(z) &
            \tilde{T}_{22}(z)T_{11}(z)-\tilde{T}_{21}(z)T_{12}(z)
        \end{pmatrix}.
    \end{equation}
    From the definition (\ref{TinS}) of $\tilde{T}(z)$ and
    (\ref{RHPSd1})--(\ref{RHPSd3}) it follows that
    \[
        \tilde{T}(z)=
        O\begin{pmatrix}
            1 & |z-1|^{\alpha}\\
            1 & |z-1|^{\alpha}
        \end{pmatrix}, \qquad\mbox{as $z \to 1$,}
    \]
    in case $\alpha < 0$, that
    \[
        \tilde{T}(z)=
        O\begin{pmatrix}
            \log|z-1| & \log|z-1|\\
            \log|z-1| & \log|z-1|
        \end{pmatrix}, \qquad \mbox{as $z\to 1$,}
    \]
    in case $\alpha = 0$, and that
    \[
        \tilde{T}(z)=
        \left\{\begin{array}{cl}
            O\begin{pmatrix}
                1 & 1\\
                1 & 1
            \end{pmatrix}, & \mbox{as $z\to 1$ outside the lens,} \\[2ex]
            O\begin{pmatrix}
                |z-1|^{-\alpha} & 1 \\
                |z-1|^{-\alpha} & 1
            \end{pmatrix}, & \mbox{as $z\to 1$ inside the lens.}
        \end{array}\right.
    \]
    in case $\alpha > 0$.

    Combining this with (\ref{RHPTd}) and (\ref{HifvT}), we obtain
    that
    \begin{equation} \label{Halpha1}
        H(z)=
        O\begin{pmatrix}
            |z-1|^{\alpha} & |z-1|^{\alpha} \\
            |z-1|^{\alpha} & |z-1|^{\alpha}
        \end{pmatrix}, \qquad\mbox{as $z\to 1$,}
    \end{equation}
    in case $\alpha < 0$, that
    \begin{equation} \label{Halpha2}
        H(z)=
        O\begin{pmatrix}
            \log|z-1| & (\log|z-1|)^{2} \\
            \log|z-1| & (\log|z-1|)^{2}
        \end{pmatrix}, \qquad\mbox{as $z\to 1$,}
    \end{equation}
    in case $\alpha = 0$, and that
    \begin{equation} \label{Halpha3}
        H(z)=
        \left\{\begin{array}{cl}
        O\begin{pmatrix}
            1 & 1\\
            1 & 1
        \end{pmatrix}, & \mbox{as  $z\to 1$  outside the lens,}\\[2ex]
        O\begin{pmatrix}
            |z-1|^{-\alpha} & |z-1|^{-\alpha}\\
            |z-1|^{-\alpha} & |z-1|^{-\alpha}
        \end{pmatrix},
        & \mbox{as $z\to 1$ inside the lens,}
        \end{array}\right.
    \end{equation}
    in case $\alpha > 0$.

    In all cases it now follows that $H$ has a removable singularity at $1$.
    This is clear if $\alpha \leq 0$ from (\ref{Halpha1}) and (\ref{Halpha2}),
    since $\alpha > -1$. If $\alpha > 0$, then it follows from (\ref{Halpha3})
    that $H$ remains bounded if we approach $1$ from outside the lens. Then
    the entries of $H$ cannot have a pole at $1$, since in the case of a pole,
    the function becomes unbounded if we approach $1$ from all directions.
    Moreover, we also get from (\ref{Halpha3}) that $(z-1)^m H(z)$ is bounded near $1$,
    for any integer $m > \alpha$. Then the entries of $H$ cannot have an
    essential singularity either, so that $1$ is a removable singularity of $H$.

    Similarly, it follows that in all cases, $H$ has a removable
    singularity at $-1$. Therefore, $H$ is entire. It follows from
    (\ref{RHPTc}), (\ref{RHPSc}), and (\ref{TinS}) that $H(z)\to\ I$ as
    $z\to\infty$. So, by Liouville's theorem, $H(z)=I$. This gives
    $\tilde{T} = T$, so that indeed $\tilde{T}$ is the (unique)
    solution of the RHP for $T$.
\end{proof}

\begin{remark}
    For $\alpha \geq 1$, the jump matrices (\ref{RHPSb1}) and
    (\ref{RHPSb2}) for $S$, as well as $S$ itself, have a
    non-integrable singularity near $z = 1$, since the factor
    $w^{-1}(z)$ behaves like $(1-z)^{-\alpha}$ as $z \to 1$. It then
    follows that the RHP for $S$ does not have a solution in any $L^p$
    sense with $p > 1$, if $\alpha \geq 1$.
\end{remark}

\begin{remark}
    Since $|\varphi(z)| > 1$ for every $z \in \mathbb C \setminus
    [-1,1]$, we see from (\ref{RHPSb1}) that the jump matrix for
    $S$ on $\Sigma_1^o \cup \Sigma_3^o$,
    \[
        \begin{pmatrix}
            1 & 0 \\
            w(z)^{-1}\varphi(z)^{-2n} & 1
        \end{pmatrix}
    \]
    tends to the identity matrix as $n \to \infty$.
    So the effect of the transformation $T \mapsto S$ has been to
    transform the oscillating entries on the diagonal in the jump
    matrix for $T$ into exponentially decaying off-diagonal entries in
    the jump matrix for $S$. This lies at the heart of the Deift--Zhou
    steepest descent method for Riemann--Hilbert problems.
\end{remark}

    \section{The parametrix for the outside region}
\setcounter{equation}{0}

Since the jump matrix in (\ref{RHPSb1}) tends to the identity
matrix as $n \to \infty$, we expect that the leading order
asymptotics are determined by a solution $N$ of the following
RHP.

\subsubsection*{RHP for \boldmath$N$:}

\begin{enumerate}
    \item[(a)]
        $N(z)$ is analytic for $z \in \mathbb C \setminus [-1,1]$.
    \item[(b)]
        $N(z)$ satisfies the following jump relation on the interval $(-1,1)$:
        \begin{equation} \label{RHPNb}
            N_+(x) = N_-(x)
            \begin{pmatrix}
                0 & w(x) \\
                -w(x)^{-1} & 0
            \end{pmatrix}, \qquad \mbox{for $x \in (-1,1)$.}
        \end{equation}
    \item[(c)]
        $N(z)$ has the following behavior as $z\to\infty$:
        \begin{equation} \label{RHPNc}
            N(z) = I + O\left(\frac{1}{z} \right), \qquad \mbox{as $z \to\infty$.}
        \end{equation}
\end{enumerate}

We present a solution of the RHP for $N$ in terms of the Szeg\H{o}
function $D(z) = D(z;w)$ associated with $w$, see
(\ref{Szegofunctie}). The Szeg\H{o} function is the scalar function
$D : \mathbb C \setminus [-1,1] \to \mathbb C$ that satisfies the
following scalar RHP.

\subsubsection*{RHP for \boldmath$D$:}

\begin{enumerate}
\item[(a)]
    $D(z)$ has no zeros and is analytic for $z \in \mathbb C \setminus [-1,1]$.
\item[(b)]
    $D_+(x) D_-(x) = w(x)$ \mbox{ for $x \in (-1,1)$.}
\item[(c)]
    $\lim\limits_{z \to \infty} D(z) = D_{\infty} \in (0, \infty)$.
\end{enumerate}
It is well-known and easy to check that $D$ as given by (\ref{Szegofunctie}) is indeed
a solution of this RHP.

\begin{remark}
    If $D_1$ and $D_2$ are the Szeg\H{o} functions for the weights $w_1$
    and $w_2$, respectively, then the product $D_1D_2$ is the Szeg\H{o}
    function for $w_1w_2$.
    Since the Szeg\H{o} function for the pure Jacobi weight
    $w_{\alpha,\beta}(x)=(1-x)^\alpha(1+x)^\beta$ is equal to
    \begin{equation}\label{SzegoD1}
        D(z;w_{\alpha,\beta})=\frac{(z-1)^{\alpha/2}(z+1)^{\beta/2}}{\varphi(z)^{(\alpha+\beta)/2}},
        \qquad\mbox{for $z\in\mathbb{C}\setminus[-1,1]$,}
    \end{equation}
    it follows that for $z\in\mathbb{C}\setminus[-1,1]$,
    \begin{eqnarray}
        \nonumber
        D(z;w) & = & D(z;w_{\alpha,\beta}) D(z;h) \\[1ex]
        \label{SzegoD}
        & = & \frac{(z-1)^{\alpha/2}(z+1)^{\beta/2}}{\varphi(z)^{(\alpha+\beta)/2}}\exp
            \left(\frac{(z^{2}-1)^{1/2}}{2\pi}\int_{-1}^{1}\frac{\log h(x)}{\sqrt{1-x^{2}}}\frac{dx}{z-x}\right).
    \end{eqnarray}
\end{remark}

\begin{proposition}
A solution of the RHP for $N$ is given by
\begin{equation}\label{RHPNsolution}
    N(z) = D_{\infty}^{\sigma_3}
    \begin{pmatrix}
        \frac{a(z) + a(z)^{-1}}{2} & \frac{a(z) - a(z)^{-1}}{2i} \\[1ex]
        \frac{a(z) - a(z)^{-1}}{-2i} & \frac{a(z) + a(z)^{-1}}{2}
    \end{pmatrix} D(z)^{-\sigma_3}.
\end{equation}
where $D(z) = D(z;w)$ is the Szeg\H{o} function associated with $w$,
and
\begin{equation} \label{az}
a(z) = \frac{(z-1)^{1/4}}{(z+1)^{1/4}}.
\end{equation}
\end{proposition}

\begin{proof}
If we seek a solution $N$ of the RHP for $N$ of the form
\begin{equation} \label{NinN1}
    N(z) = D_{\infty}^{\sigma_3} N^{(1)}(z) D(z)^{-\sigma_3},
\end{equation}
then $N^{(1)}(z)$ is analytic in $\mathbb{C}\setminus[-1,1]$, and
\[
    N^{(1)}(z)=I+O\left(\frac{1}{z}\right),\qquad\mbox{as $z\to\infty$.}
\]
Furthermore, using (\ref{RHPNb}) and $D_+(x)D_-(x)=w(x)$ for $x\in(-1,1)$,
we have
\[
    N_+^{(1)}(x) = N_-^{(1)}(x)
    \begin{pmatrix}
        0 & 1 \\
        -1 & 0
    \end{pmatrix}, \qquad \mbox{for $x \in (-1,1)$.}
\]
It is well-known (see \cite{Deift,DKMVZ1,DKMVZ2}) that the RHP
for $N^{(1)}$ is solved by
\[
    N^{(1)}(z) =
    \begin{pmatrix}
        1 & 1 \\
        i & -i
    \end{pmatrix}
    a(z)^{\sigma_3}
    \begin{pmatrix}
        1 & 1 \\
        i & -i
    \end{pmatrix}^{-1} =
    \begin{pmatrix}
        \frac{a(z) + a(z)^{-1}}{2} & \frac{a(z) - a(z)^{-1}}{2i} \\[1ex]
        \frac{a(z) - a(z)^{-1}}{-2i} & \frac{a(z) + a(z)^{-1}}{2}
    \end{pmatrix},
\]
where $a(z)$ is given by (\ref{az}). This proves (\ref{RHPNsolution}).
\end{proof}

\section{Parametrix near the endpoints}
\setcounter{equation}{0}

The jump matrices for $S$ and $N$ are not uniformly close to each
other near the endpoints $\pm 1$. This is reflected in the fact
that $S N^{-1}$ is not bounded near the endpoints. We need a local
analysis around these points. We consider the endpoint $z=1$ in
detail, the analysis for $z =-1$ being similar. In a small but
fixed neighborhood $U_{\delta} = \{z \mid |z-1| < \delta\}$ of $z
= 1$, with $0 < \delta< r$ (so that the closure of $U_{\delta}$ lies in
$U$) we construct a parametrix $P$,  that satisfies exactly the
same jump relations as $S$, that matches with $N$ on the boundary
$\partial U_{\delta}$ of $U_{\delta}$, and that has the same
behavior as $S(z)$ near $z=1$. We thus seek a matrix valued function
$P$ that satisfies the following RHP.

\subsubsection*{RHP for \boldmath$P$:}

\begin{enumerate}
\item[(a)]
    $P(z)$ is defined and analytic for $z \in U_{\delta_0}\setminus\Sigma$
    for some $\delta_0 > \delta$.
\item[(b)]
    $P(z)$ satisfies the following jump relations on
    $U_{\delta} \cap \Sigma^o$:
    \begin{equation}\label{RHPPb1}
        P_+(z) = P_-(z)
        \begin{pmatrix}
            1 & 0 \\
            w(z)^{-1}\varphi (z)^{-2n} & 1
        \end{pmatrix}, \qquad\mbox{for $z\in U_{\delta} \cap \left(\Sigma_{1}^o
        \cup\Sigma_{3}^o\right)$,}
    \end{equation}
    \begin{equation}\label{RHPPb2}
        P_+(x) = P_-(x)
        \begin{pmatrix}
            0 & w(x) \\
            -w(x)^{-1} & 0
        \end{pmatrix}, \qquad\mbox{for $x \in U_{\delta} \cap \Sigma_{2}^o$.}
    \end{equation}
    \item[(c)]
        On $\partial U_{\delta}$ we have, as $n \to \infty$
        \begin{equation} \label{RHPPc}
            P(z) N^{-1}(z) = I + O \left( \frac{1}{n} \right),
            \qquad \mbox{uniformly for $z \in \partial U_{\delta}\setminus\Sigma$.}
        \end{equation}
    \item[(d)]
        For $\alpha<0$, the matrix function $P(z)$ has the following behavior as $z\to 1$:
        \begin{equation}\label{RHPPd1}
            P(z)=
            O\begin{pmatrix}
                1 & |z-1|^{\alpha} \\
                1 & |z-1|^{\alpha}
            \end{pmatrix},\qquad \mbox{as $z\to 1 ,z\in U_{\delta}\setminus\Sigma$.}
        \end{equation}
        For $\alpha=0$, $P(z)$ has the following behavior as $z\to 1$:
        \begin{equation}\label{RHPPd2}
            P(z)=
            O\begin{pmatrix}
                \log|z-1| & \log|z-1| \\
                \log|z-1| & \log|z-1|
            \end{pmatrix}, \qquad \mbox{as $z\to 1, z \in U_{\delta}\setminus\Sigma$.}
        \end{equation}
        For $\alpha>0$, $P(z)$ has the following behavior as $z\to 1$:
        \begin{equation}\label{RHPPd3}
            P(z)=
            \left\{\begin{array}{cl}
                O\begin{pmatrix}
                    1 & 1 \\
                    1 & 1
                \end{pmatrix},& \mbox{as $z\to 1$ outside the lens,}\\[2ex]
                O\begin{pmatrix}
                    |z-1|^{-\alpha} & 1 \\
                    |z-1|^{-\alpha} & 1
                \end{pmatrix}, & \mbox{as $z\to 1$ inside the lens.}
            \end{array}\right.
        \end{equation}
\end{enumerate}

\medskip

To find $P$, we concentrate first on parts (a), (b), and (d).
Later we will also consider the matching condition (\ref{RHPPc}).

Recall that $h$ is analytic with positive real part on $U$. Hence, the function
\begin{equation} \label{defW(z)}
    W(z) = \left( (z-1)^{\alpha} (z+1)^{\beta} h(z) \right)^{1/2}
\end{equation}
is defined and analytic for $z \in U \setminus (-\infty, 1]$.
The branch of the square root is chosen which is positive for $z > 1$.
Note that
\begin{equation} \label{W^2(z)}
    W^2(z) = \left\{ \begin{array}{rl}
     e^{\alpha \pi i} w(z), & \qquad \mbox{for } \Im z > 0, \\[10pt]
     e^{-\alpha \pi i} w(z), & \qquad \mbox{for } \Im z < 0.
     \end{array} \right.
\end{equation}
Thus in particular, $W^2_{\pm}(x) =  e^{\pm \alpha \pi i} w(x)$
for $x \in (-1, 1)$, and so
\begin{equation}\label{W^2(z)bis}
     W_{+}(x) W_{-}(x) = w(x), \qquad \mbox{for } x \in (-1,1).
\end{equation}
Thus $W$ satisfies the same jump relation as the Szeg\H{o} function
$D(z)$. Note however that $W$ does not extend to an analytic
function on $\mathbb C\setminus [-1,1]$.

We seek $P$ in the form
\begin{equation} \label{SformS1}
    P(z) = E_n(z) P^{(1)}(z) W(z)^{-\sigma_3} \varphi(z)^{-n\sigma_{3}},
\end{equation}
where the invertible matrix valued function $E_n$  is analytic in
a neighborhood of $U_{\delta}$. By (\ref{RHPPb1}) and (\ref{SformS1}) the matrix
function $P^{(1)}$ should satisfy on $U_{\delta} \cap
(\Sigma_{1}^o \cup \Sigma_{3}^o)$ the jump relation
\begin{eqnarray*}
    P^{(1)}_+(z) & = & P^{(1)}_-(z)\varphi(z)^{-n\sigma_3} W(z)^{-\sigma_3}
        \begin{pmatrix}
            1 & 0 \\
            w(z)^{-1}\varphi(z)^{-2n} & 1
        \end{pmatrix}
        W(z)^{\sigma_3}\varphi(z)^{n\sigma_3} \\[1ex]
    & = & P^{(1)}_-(z)
        \begin{pmatrix}
            1 & 0 \\
            w(z)^{-1} W^2(z) & 1
        \end{pmatrix} \\[1ex]
    & = & P^{(1)}_-(z)
        \begin{pmatrix}
            1 & 0 \\
            e^{\pm \alpha \pi i} & 1
        \end{pmatrix}
\end{eqnarray*}
where in the last step we used (\ref{W^2(z)}). In $e^{\pm \alpha
\pi i}$ the $+$ holds on $\Sigma_1^o$ and the $-$ on $\Sigma_3^o$.
For $x \in (1-\delta,1)$, we find using (\ref{RHPPb2}),
(\ref{W^2(z)bis}), (\ref{SformS1}), and the fact that $\varphi_{+}(x)\varphi_{-}(x)
= 1$,
\begin{eqnarray*}
    P^{(1)}_+(x) & = & P^{(1)}_-(x)\varphi_-(x)^{-n\sigma_3} W_-(x)^{-\sigma_3}
        \begin{pmatrix}
            0 & w(x) \\
            -w(x)^{-1}  & 0
        \end{pmatrix}
        W_+(x)^{\sigma_3}\varphi_+(x)^{n\sigma_3} \\[1ex]
    & = & P^{(1)}_-(x)
        \begin{pmatrix}
            0 & w(x) W_-(x)^{-1} W_+(x)^{-1}  \\
            - w(x)^{-1}  W_-(x) W_+(x) & 0
        \end{pmatrix} \\[1ex]
    & = & P^{(1)}_-(x)
        \begin{pmatrix}
            0 & 1 \\
            -1  & 0
        \end{pmatrix}.
\end{eqnarray*}
We therefore see that we must look for a matrix valued function
$P^{(1)}$ that satisfies the following RHP.

\subsubsection*{RHP for \boldmath$P^{(1)}$:}
\begin{enumerate}
    \item[(a)]
        $P^{(1)}(z)$ is defined and analytic for $z\in U_{\delta_0}\setminus\Sigma$ for some $\delta_0 > \delta$.
    \item[(b)]
        $P^{(1)}(z)$ satisfies the following jump relations:
        \begin{equation}\label{RHPP1b}
            \left\{\begin{array}{ll}
                P^{(1)}_{+}(z)=P^{(1)}_{-}(z)
                \begin{pmatrix}
                    1 & 0 \\
                    e^{\alpha\pi i} & 1
                \end{pmatrix}, & \mbox{for $z\in U_{\delta}\cap\Sigma_{1}^o$,} \\[2ex]
                P^{(1)}_{+}(x)=P^{(1)}_{-}(x)
                \begin{pmatrix}
                    0 & 1 \\
                    -1 & 0
                \end{pmatrix}, & \mbox{for $x\in U_{\delta}\cap\Sigma_{2}^o$,} \\[2ex]
                P^{(1)}_{+}(z)=P^{(1)}_{-}(z)
                \begin{pmatrix}
                    1 & 0 \\
                    e^{-\alpha\pi i} & 1
                \end{pmatrix}, & \mbox{for $z\in U_{\delta}\cap\Sigma_{3}^o$.}
            \end{array} \right.
        \end{equation}
    \item[(c)]
        For $\alpha<0$, the matrix function $P^{(1)}(z)$ has the following behavior as $z\to 1$:
        \[
            P^{(1)}(z)=
            O\begin{pmatrix}
                |z-1|^{\alpha/2} & |z-1|^{\alpha/2}\\
                |z-1|^{\alpha/2} & |z-1|^{\alpha/2}
            \end{pmatrix},
            \qquad\mbox{as $z\to 1, z\in U_{\delta}\setminus\Sigma$.}
        \]
        For $\alpha=0$, $P^{(1)}(z)$ has the following behavior as $z\to 1$:
        \[
            P^{(1)}(z)=
            O\begin{pmatrix}
                \log|z-1| & \log|z-1|\\
                \log|z-1| & \log|z-1|
            \end{pmatrix},
            \qquad \mbox{as $z\to 1, z \in U_{\delta}\setminus\Sigma$.}
        \]
        For $\alpha>0$, $P^{(1)}(z)$ has the following behavior as $z\to 1$:
        \[
            P^{(1)}(z)=
            \left\{\begin{array}{cl}
                O\begin{pmatrix}
                    |z-1|^{\alpha/2} & |z-1|^{-\alpha/2}\\
                    |z-1|^{\alpha/2} & |z-1|^{-\alpha/2}
                \end{pmatrix}, & \mbox{as $z\to 1$ outside the lens,}\\[2ex]
                O\begin{pmatrix}
                    |z-1|^{-\alpha/2} & |z-1|^{-\alpha/2}\\
                    |z-1|^{-\alpha/2} & |z-1|^{-\alpha/2}
                \end{pmatrix}, & \mbox{as $z\to 1$ inside the lens.}
            \end{array}\right.
        \]
\end{enumerate}

\begin{remark}
    Note that the condition (c) in the RHP for $P^{(1)}$ follows from
    the condition (d) in the RHP for $P$, since
    \[
        P^{(1)}(z) = E_n(z)^{-1} P(z) W(z)^{\sigma_3} \varphi(z)^{-n\sigma_3}
    \]
    where $\varphi$ is bounded and bounded away from $0$ near $z = 1$,
    and $W(z)$ behaves like $c (z-1)^{\alpha/2}$, with a positive constant $c$,
    as $z \to 1$.
\end{remark}

\begin{remark}
    Observe that the jump matrices for $P^{(1)}$ in (\ref{RHPP1b}) are
    constants on the three parts $\Sigma_1^o$, $\Sigma_2^o$ and $\Sigma_3^o$.
\end{remark}

To construct $P^{(1)}$ we first define a function $g$ as
\begin{equation} \label{defg(z)}
    g(z)=\log\varphi(z)=\log\left(z+(z^{2}-1)^{1/2}\right),\qquad\mbox{for
    $z\in\mathbb{C}\setminus(-\infty,1]$.}
\end{equation}
This function is analytic in $\mathbb{C}\setminus (-\infty,1]$ and
for $x\in(-1,1)$ we have, since $\varphi_{+}(x)\varphi_{-}(x)=1$,
that $g_{+}(x)=-g_{-}(x)$. Therefore $g^{2}$ is also analytic
across $(-1,1)$. Since $g^{2}(1)=0$, the singularity at $1$ is removable,
so that the function $f$ defined by
\begin{equation} \label{deff(z)}
    f(z)=g^{2}(z)/4, \qquad \mbox{for $z \in \mathbb{C}\setminus(-\infty,-1]$,}
\end{equation}
is analytic in $\mathbb C \setminus (-\infty,-1]$.
The behavior of $g(z)$ near 1 is
\begin{equation}\label{Behaviorg}
    g(z)=\sqrt{2}(z-1)^{1/2}-\frac{\sqrt2}{12}(z-1)^{3/2}+O\left((z-1)^{5/2}\right),
    \qquad\mbox{as $z\to 1$,}
\end{equation}
and the behavior of $f(z)$ near 1 is
\begin{equation} \label{Behaviorf}
    f(z)=\frac{1}{2}(z-1)-\frac{1}{12}(z-1)^{2}+O\left((z-1)^{3}\right),
    \qquad\mbox{as $z\to 1$.}
\end{equation}

Because $f$ is analytic on $\mathbb{C}\setminus (-\infty,-1]$
and $f'(1)=1/2\neq 0$, $f$ is a one-to-one conformal mapping on a
neighborhood of $1$. A more detailed analysis of $f$ reveals that
$f$ is one-to-one on the disk $|z-1| < 2$ around $1$.
In any case, since $0 < \delta < r < 1$, we have that
$f$ is one-to-one on $U_{\delta}$. Since $f(z)$ is real for $z$ real,
we then have that $f$ maps $U_{\delta}\cap\mathbb{C}_{\pm}$ one-to-one
onto $f(U_{\delta})\cap\mathbb{C}_{\pm}$,
where $\mathbb{C}_{+}=\{z \mid \Im z>0\}$ and $\mathbb{C}_{-}=\{z \mid \Im z<0\}$.

\medskip

\begin{figure}
    \center{\resizebox{8cm}{!}{\includegraphics{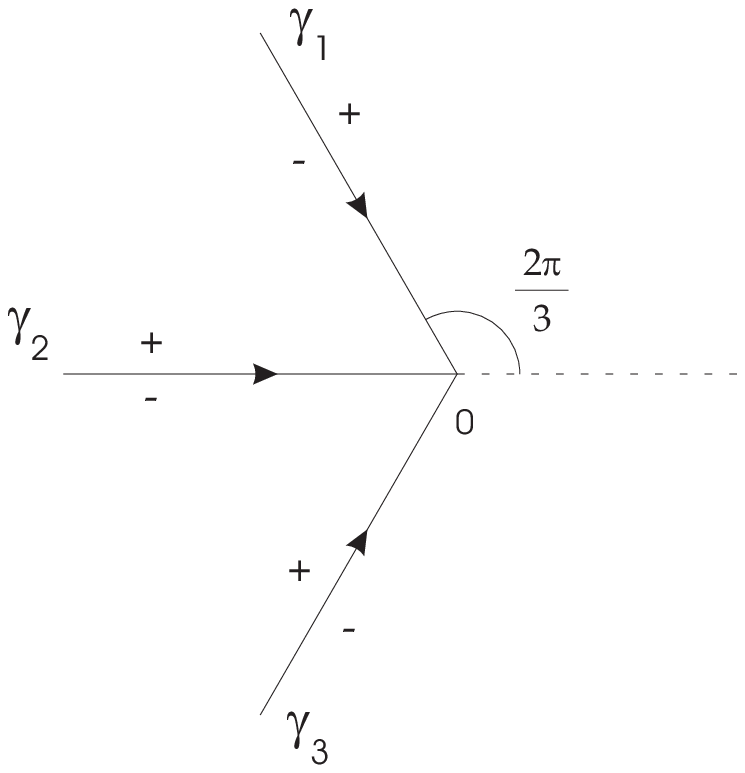}}
    \caption{The contour $\Sigma_{\Psi}$}\label{figuur3}}
\end{figure}

We will use the mapping $\zeta = n^2f(z)$ to transfer the RHP in
the $z$-plane to a RHP in the $\zeta$-plane.
Consider the following RHP in the complex $\zeta$-plane, on a contour $\Sigma_{\Psi}$
consisting of three infinite rays,
\[
    \gamma_1 :  \arg \zeta = 2\pi/3, \qquad
    \gamma_2 : \arg \zeta = \pi, \qquad
    \gamma_3 : \arg \zeta = -2\pi/3,
\]
oriented as in Figure \ref{figuur3}.

\subsubsection*{RHP for \boldmath$\Psi$:}

\begin{enumerate}
    \item[(a)]
        $\Psi(\zeta)$ is analytic for $\zeta\in\mathbb{C}\setminus\Sigma_{\Psi}$,
    \item[(b)]
        $\Psi(\zeta)$ satisfies the following jump relations:
        \begin{eqnarray}
            \label{RHPPSIb1}
            \Psi_+(\zeta) &=& \Psi_-(\zeta)
            \begin{pmatrix}
                1 & 0 \\
                e^{\alpha \pi i}  & 1
            \end{pmatrix}, \qquad \mbox{for $\zeta \in \gamma_{1}$,} \\[2ex]
            \label{RHPPSIb2}
            \Psi_+(\zeta) &=& \Psi_-(\zeta)
            \begin{pmatrix}
                0 & 1 \\
                -1  & 0
            \end{pmatrix}, \qquad \mbox{for $\zeta \in \gamma_{2}$,} \\[2ex]
            \label{RHPPSIb3}
            \Psi_+(\zeta) &=& \Psi_-(\zeta)
            \begin{pmatrix}
                1 & 0 \\
                e^{-\alpha \pi i}  & 1
            \end{pmatrix}, \qquad \mbox{for $\zeta \in \gamma_{3}$.}
        \end{eqnarray}
    \item[(c)]
        For $\alpha<0$, the matrix function $\Psi(\zeta)$ has the
        following behavior as $\zeta\to 0$:
        \begin{equation} \label{RHPPSIc1}
            \Psi(\zeta)=
            O\begin{pmatrix}
                |\zeta|^{\alpha/2} & |\zeta|^{\alpha/2}\\
                |\zeta|^{\alpha/2} & |\zeta|^{\alpha/2}
            \end{pmatrix}, \qquad \mbox{as $\zeta\to 0$.}
        \end{equation}
        For $\alpha=0$, $\Psi(\zeta)$ has the following behavior as $\zeta\to 0$:
        \begin{equation} \label{RHPPSIc2}
            \Psi(\zeta)=
                O\begin{pmatrix}
                    \log|\zeta| & \log|\zeta|\\
                    \log|\zeta| & \log|\zeta|
                \end{pmatrix}, \qquad \mbox{as $\zeta \to 0$.}
        \end{equation}
        For $\alpha>0$, $\Psi(\zeta)$ has the following behavior as $\zeta\to 0$:
        \begin{equation} \label{RHPPSIc3}
            \Psi(\zeta)=
            \left\{\begin{array}{cl}
                O\begin{pmatrix}
                    |\zeta|^{\alpha/2} & |\zeta|^{-\alpha/2}\\
                    |\zeta|^{\alpha/2} & |\zeta|^{-\alpha/2}
                \end{pmatrix},
                & \mbox{as $\zeta\to 0$ in $|\arg\zeta|<2\pi/3$},\\[2ex]
                O\begin{pmatrix}
                    |\zeta|^{-\alpha/2} & |\zeta|^{-\alpha/2}\\
                    |\zeta|^{-\alpha/2} & |\zeta|^{-\alpha/2}
                \end{pmatrix}, & \mbox{as $\zeta\to 0$ in $2\pi/3 < |\arg\zeta| <\pi$}.
            \end{array}\right.
        \end{equation}
\end{enumerate}

\begin{figure}
    \center{ \resizebox{12cm}{!}{\includegraphics{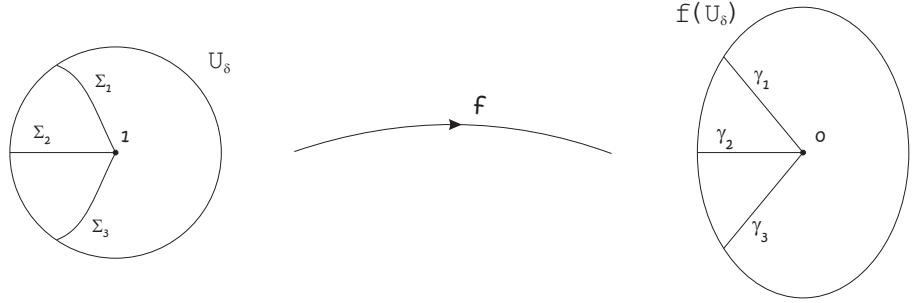}}
    \caption{The definition of $\Sigma_1$, $\Sigma_{2}$ and
    $\Sigma_{3}$ in $U_{\delta}$}\label{figuur4}}
\end{figure}

In the above, we have chosen the angle $2\pi/3$ rather
arbitrarily. It could be replaced by any angle in $(0, \pi)$.
However, having chosen the angle, we {\em define} the contours
$\Sigma_{1}\cap U_{\delta}$ and $\Sigma_{3}\cap U_{\delta}$ in
$U_{\delta}$ as the preimages of the rays $\gamma_{1}$ and
$\gamma_{3}$ under the mapping $\zeta=n^2 f(z)$, see Figure
\ref{figuur4}. Note that this definition does not depend on $n$.
Then it follows immediately that if $\Psi(\zeta)$ is a solution of
the RHP for $\Psi$ we can define
\begin{equation} \label{P1inPsi}
    P^{(1)}(z) = \Psi(n^{2}f(z)),
\end{equation}
and $P^{(1)}$ will satisfy the RHP for $P^{(1)}$.

\medskip

We construct a solution $\Psi$ out of (modified) Bessel functions
of order $\alpha$. We refer to Chapter 9 of
\cite{AbramowitzStegun} for definitions and properties of the
modified Bessel functions $I_{\alpha}$ and $K_{\alpha}$, and the
Hankel functions $H_{\alpha}^{(1)}$ and $H_{\alpha}^{(2)}$, as
well as their relation with the ordinary Bessel functions
$J_{\alpha}$ and $Y_{\alpha}$. We define $\Psi(\zeta)$ for $|\arg
\zeta| < 2 \pi/3$ as
\begin{equation}\label{RHPPSIsolution1}
    \Psi(\zeta) =
    \begin{pmatrix}
        I_{\alpha} (2 \zeta^{1/2}) & \frac{i}{\pi} K_{\alpha}(2 \zeta^{1/2}) \\[1ex]
        2\pi i \zeta^{1/2} I_{\alpha}'(2\zeta^{1/2}) & -2 \zeta^{1/2} K_{\alpha}'(2\zeta^{1/2})
    \end{pmatrix}.
\end{equation}
For $2\pi/3 < \arg \zeta < \pi$ we define it as
\begin{equation}\label{RHPPSIsolution2}
    \Psi(\zeta) =
    \begin{pmatrix}
        \frac{1}{2} H_{\alpha}^{(1)}(2(-\zeta)^{1/2}) &
        \frac{1}{2} H_{\alpha}^{(2)}(2(-\zeta)^{1/2}) \\[1ex]
        \pi \zeta^{1/2} \left(H_{\alpha}^{(1)}\right)'(2(-\zeta)^{1/2}) &
        \pi \zeta^{1/2} \left(H_{\alpha}^{(2)}\right)'(2(-\zeta)^{1/2})
    \end{pmatrix} e^{\frac{1}{2}\alpha\pi i \sigma_3 }.
\end{equation}
And finally for $- \pi < \arg \zeta  < -2\pi/3$ it is defined as
\begin{equation}\label{RHPPSIsolution3}
    \Psi(\zeta) =
    \begin{pmatrix}
        \frac{1}{2} H_{\alpha}^{(2)}(2(-\zeta)^{1/2}) &
        -\frac{1}{2} H_{\alpha}^{(1)}(2 (-\zeta)^{1/2}) \\[1ex]
        -\pi \zeta^{1/2} \left(H_{\alpha}^{(2)}\right)'(2 (- \zeta)^{1/2}) &
        \pi \zeta^{1/2} \left(H_{\alpha}^{(1)}\right)'(2 (-\zeta)^{1/2})
    \end{pmatrix}
    e^{-\frac{1}{2} \alpha \pi i \sigma_3}.
\end{equation}

\begin{theorem}
    The matrix valued function $\Psi(\zeta)$, defined  in
    {\rm(\ref{RHPPSIsolution1})--(\ref{RHPPSIsolution3})}, is a
    solution of the RHP for $\Psi$.
\end{theorem}

\begin{proof}
    (a) Since the functions $I_{\alpha}$, $K_{\alpha}$,
    $H_{\alpha}^{(1)}$, and $H_{\alpha}^{(2)}$ (as well as their
    derivatives) are defined and analytic in the complex plane with a
    branch cut along the negative real axis, the matrix valued
    function $\Psi$ defined by (\ref{RHPPSIsolution1})--(\ref{RHPPSIsolution3})
    is analytic in the respective regions.

    \medskip

    (b) We need to verify that the jump conditions
    (\ref{RHPPSIb1})--(\ref{RHPPSIb3}) are satisfied. We first verify
    the jump condition (\ref{RHPPSIb2}) on the negative real axis. The
    Hankel functions evaluated at $2(-\zeta)^{1/2}$ that appear in
    (\ref{RHPPSIsolution2}) and (\ref{RHPPSIsolution3}) are analytic
    across the negative $\zeta$-axis. Only the factor $\zeta^{1/2}$
    appearing in the second rows of (\ref{RHPPSIsolution2}) and
    (\ref{RHPPSIsolution3}) has different boundary values on the $+$
    and $-$ sides of the negative real axis. Taking this into account,
    we easily see that the jump condition (\ref{RHPPSIb2}) is satisfied.

    Next, we check (\ref{RHPPSIb1}). For $\zeta\in \gamma_{1}$ we use
    (\ref{RHPPSIsolution1}) to evaluate $\Psi_{+}(\zeta)$ and
    (\ref{RHPPSIsolution2}) to evaluate $\Psi_{-}(\zeta)$. From
    (\ref{RHPPSIsolution2}) the $(1,1)$ entry on the right of
    (\ref{RHPPSIb1}) is
    \[
        e^{\frac{1}{2} \alpha \pi i} \frac{1}{2} \left( H_{\alpha}^{(1)}(2(-\zeta)^{1/2}) +
        H_{\alpha}^{(2)}(2(-\zeta)^{1/2}) \right).
    \]
    Using $(-\zeta)^{1/2}=\zeta^{1/2}\exp(-\frac{1}{2}\pi i)$ for
    $\zeta\in\gamma_1$ and formulas 9.1.3 and 9.1.4 of
    \cite{AbramowitzStegun}, we find that this is equal to
    \begin{equation} \label{eq529}
        e^{\frac{1}{2}\alpha\pi i} J_{\alpha}(2\zeta^{1/2}e^{-\frac{1}{2}\pi i}).
    \end{equation}
    If we now use the fact that, cf. formulas 9.1.35 and 9.6.3 of
    \cite{AbramowitzStegun},
    \[
        I_{\alpha}(z) = e^{\frac{1}{2} \alpha \pi i} J_{\alpha}(z
        e^{-\frac{1}{2} \pi i}),\qquad\mbox{for }-\pi<\arg z<\frac{\pi}{2},
    \]
    we see that (\ref{eq529}) is equal to $I_{\alpha}(2\zeta^{1/2})$,
    and this is the $(1,1)$ entry on the left of (\ref{RHPPSIb1}). So
    the $(1,1)$ entries of both sides of (\ref{RHPPSIb1}) agree.

    The $(1,2)$ entry on the left of (\ref{RHPPSIb1}) is $(i/\pi)
    K_{\alpha}(2\zeta^{1/2})$. Using formula 9.6.4 of
    \cite{AbramowitzStegun}, we rewrite $K_{\alpha}$ in terms of the
    Hankel function $H_{\alpha}^{(2)}$ as follows
    \[
        \frac{i}{\pi} K_{\alpha}(2 \zeta^{1/2})
        = \frac{1}{2}\ e^{-\frac{1}{2} \alpha \pi i} H_{\alpha}^{(2)}(2 (-\zeta)^{1/2}).
    \]
    This is exactly the $(1,2)$ entry on the right of
    (\ref{RHPPSIb2}). So the $(1,2)$ entries of both sides (\ref{RHPPSIb1}) also
    agree.

    Now we know that the first rows of the left-hand side and the
    right-hand side of (\ref{RHPPSIb1}) are the same. To see that the
    second rows agree as well, we observe that in both
    (\ref{RHPPSIsolution1}) and (\ref{RHPPSIsolution2}) (and also
    (\ref{RHPPSIsolution3})) the second row is equal to the derivative
    with respect to $\zeta$ of the first row multiplied by the factor
    $2\pi i \zeta$. Since the jump matrix does not depend on $\zeta$,
    the equality of the first rows yields the equality of the second
    rows as well. Thus (\ref{RHPPSIb1}) holds. Similar considerations
    show that (\ref{RHPPSIb3}) holds as well.

    \medskip

    (c) The behavior near $0$ stated in (\ref{RHPPSIc1})--(\ref{RHPPSIc3})
    follows easily from formulas 9.6.4, 9.6.6,
    9.6.7, 9.6.8, 9.6.9, 9.6.26 and 9.6.27 of \cite{AbramowitzStegun}.
\end{proof}

Having $\Psi$ we define $P^{(1)}(z) = \Psi(n^2 f(z))$ as in
(\ref{P1inPsi}) and $P^{(1)}$ is a solution of the RHP for
$P^{(1)}$. Then we define as in (\ref{SformS1})
\begin{equation}\label{RHPPsolution}
    P(z)=E_{n}(z)\Psi(n^{2}f(z))W(z)^{-\sigma_{3}}\varphi(z)^{-n\sigma_{3}},
\end{equation}
with $E_{n}$ an analytic matrix valued function in $U_{\delta_0}$.
We do not have $E_n$ yet, but no matter how it is chosen, $P$
satisfies parts (a), (b) and (d) of the RHP for $P$. The
multiplication on the left by the analytic factor $E_n$ has no
influence on parts (a), (b), and (d). Observe that for part (d)
it is important that the matrices in the $O$-terms in
(\ref{RHPPd1})--(\ref{RHPPd3}) are constant along the columns.

\medskip

We use the freedom we have in choosing $E_{n}$ to ensure that
(\ref{RHPPsolution}) satisfies the matching condition
(\ref{RHPPc}) of the RHP for $P$. To that end we need the
asymptotic behavior of the Bessel functions at infinity. Inserting
the asymptotics of the modified Bessel functions $I_{\alpha}$ and
$K_{\alpha}$ as given in formulas 9.7.1--9.7.4 of
\cite{AbramowitzStegun} into (\ref{RHPPSIsolution1}), we find that
uniformly as $\zeta \to \infty$ in the sector $|\arg \zeta| <
2\pi/3$,
\begin{eqnarray}\nonumber
    \Psi(\zeta) & = &
    \begin{pmatrix}
        \frac{1}{2\sqrt{\pi}} \zeta^{-1/4} e^{2\zeta^{1/2}} \left(1 + O(\zeta^{-1/2})\right) &
        \frac{i}{2\sqrt{\pi}} \zeta^{-1/4} e^{-2\zeta^{1/2}} \left(1 + O(\zeta^{-1/2})\right) \\[1ex]
        \sqrt{\pi} i \zeta^{1/4} e^{2\zeta^{1/2}} \left(1+ O(\zeta^{-1/2})\right) &
        \sqrt{\pi} \zeta^{1/4} e^{-2\zeta^{1/2}} \left(1 + O(\zeta^{-1/2})\right)
    \end{pmatrix}\\[1ex]
    \label{AsymptoticsPSI}
    & = & \left(2\pi \zeta^{1/2} \right)^{-\sigma_3/2} \frac{1}{\sqrt{2}}
    \begin{pmatrix}
        1 + O(\zeta^{-1/2}) & i + O(\zeta^{-1/2}) \\[1ex]
        i + O(\zeta^{-1/2}) & 1 + O(\zeta^{-1/2})
    \end{pmatrix}
    e^{2\zeta^{1/2} \sigma_3}.
\end{eqnarray}
Inserting the asymptotics of the Hankel functions given in
formulas 9.2.3, 9.2.4, 9.2.13, and 9.2.14 of
\cite{AbramowitzStegun} into (\ref{RHPPSIsolution2}) and
(\ref{RHPPSIsolution3}), we find that the asymptotics of
$\Psi(\zeta)$ in the other two sectors are given by {\em exactly
the same} formula (\ref{AsymptoticsPSI}). Thus
(\ref{AsymptoticsPSI}) holds uniformly as $\zeta\to \infty$,
$\zeta\in\mathbb{C}\setminus\Sigma_{\Psi}$.
Taking $\zeta = n^2f(z)$, we then have
\begin{equation} \label{eq534}
    \Psi(n^2 f(z)) = \left(2\pi n \right)^{-\sigma_3/2} f(z)^{-\sigma_3/4} \frac{1}{\sqrt{2}}
    \begin{pmatrix}
        1 + O(\frac{1}{n}) & i + O(\frac{1}{n}) \\[1ex]
        i + O(\frac{1}{n}) & 1 + O(\frac{1}{n})
    \end{pmatrix}
    e^{2n f^{1/2}(z) \sigma_3},
\end{equation}
as $n\to\infty$, uniformly for $z \in \partial U_{\delta}\setminus\Sigma$.
Since $2f^{1/2} = g = \log \phi$ we have $\exp(2n f^{1/2}(z)\sigma_3)
=\varphi(z)^{n\sigma_{3}}$. Combining (\ref{RHPPsolution}) and
(\ref{eq534}) we then get
\begin{equation} \label{eq535}
    P(z) = E_n(z) \left(2\pi n \right)^{-\sigma_3/2} f(z)^{-\sigma_3/4} \frac{1}{\sqrt{2}}
    \begin{pmatrix}
        1 + O(\frac{1}{n}) & i + O(\frac{1}{n}) \\[1ex]
        i + O(\frac{1}{n}) & 1 + O(\frac{1}{n})
    \end{pmatrix}
    W(z)^{-\sigma_3},
\end{equation}
as $n\to\infty$, uniformly for $z \in \partial U_{\delta} \setminus
\Sigma$. Since we want $P$ to match with $N$ on $\partial U_{\delta}$,
we now {\em define} $E_n$ as
\begin{equation} \label{eq536}
    E_n(z) =  N(z) W(z)^{\sigma_3} \frac{1}{\sqrt{2}}
    \begin{pmatrix}
        1 & -i \\
        -i & 1
    \end{pmatrix}
    f(z)^{\sigma_3/4} \left(2\pi n \right)^{\sigma_3/2}.
\end{equation}
Then $E_n$ is defined and analytic in $U \setminus (-\infty,1]$.
With this $E_n$, we see from (\ref{eq535}) that
\begin{equation}\label{eq537}
    P(z) N^{-1}(z) = N(z) W(z)^{\sigma_3} \left( I + O\left(\frac{1}{n}\right) \right)
    W(z)^{-\sigma_3} N^{-1}(z),
\end{equation}
as $n\to\infty$, uniformly for $z \in \partial U_{\delta}
\setminus \Sigma$. Since $W$ as well as all entries of $N$ are
uniformly bounded and uniformly bounded away from $0$ on $\partial
U_{\delta}$, we get from (\ref{eq537}) that the matching condition
(\ref{RHPPc}) is satisfied. It is also clear that $E_n$ is
invertible, since by (\ref{eq536}) $E_n$ is a product of five
matrices, all with determinant one, so that
\begin{equation} \label{detEn}
    \det E_n(z) = 1.
\end{equation}
\medskip

Now everything is fine, except for the fact that $E_n$ from
(\ref{eq536}) is defined and analytic on $U \setminus
(-\infty, 1]$, but we need it to be analytic in a full
neighborhood of $1$. We check that $E_n$ has no
jumps across the interval $(-1,1)$, and in addition that it
has a removable singularity at $1$. To that end, we first
determine the behavior of $W/D$ near the point $z=1$. The lemma
states more than what is needed for the proof that $E_n$
is analytic. The extra information will be useful later on.

\begin{lemma} \label{BehaviorWoverD}
    For $z \in U_{\delta}$, we have
    \[
        \left(\frac{W(z)}{D(z)}\right)^2
        =\varphi(z)^{\alpha+\beta}\exp\left((z^2-1)^{1/2}
        \sum_{n=0}^{\infty}c_{n}(z-1)^{n}\right),
    \]
    with coefficients $c_n$ given by {\rm(\ref{Definitiecn})}.
    In particular, we have that
    \[ \frac{W(z)}{D(z)} = 1 + O(|z-1|^{1/2}), \qquad \mbox{as } z \to 1. \]
\end{lemma}

\begin{proof}
    We note that by (\ref{SzegoD})
    \begin{equation}\label{SzegoDbis}
        D^2(z)=\frac{(z-1)^{\alpha}(z+1)^{\beta}}{\varphi(z)^{\alpha+\beta}}
        \exp\left( \frac{(z^{2}-1)^{1/2}}{\pi}\int_{-1}^{1}\frac{\log
        h(x)}{\sqrt{1-x^{2}}}\frac{dx}{z-x}\right).
    \end{equation}
    Let $\gamma$ be a closed contour in $U$ going around the interval
    $[-1,1]$ in the positive direction and also encircling the point $z \in U_{\delta}$.
    After an easy residue calculation, we then find
    \[
        \frac{1}{\pi} \int_{-1}^{1}\frac{\log h(x)}{\sqrt{1-x^2}}\frac{dx}{z-x}=
        \frac{\log h(z)}{(z^2-1)^{1/2}} - \frac{1}{2\pi i}\int_{\gamma}\frac{\log
        h(\zeta)}{(\zeta^2-1)^{1/2}}\frac{d\zeta}{\zeta-z}.
    \]
    Inserting this in (\ref{SzegoDbis}), we get
    \[
        D^2(z) = \frac{(z-1)^{\alpha}(z+1)^{\beta}}{\varphi(z)^{\alpha+\beta}}
            h(z) \exp \left( - (z^2-1)^{1/2} \frac{1}{2\pi i} \int_{\gamma}\frac{\log
        h(\zeta)}{(\zeta^2-1)^{1/2}}\frac{d\zeta}{\zeta-z}\right).
    \]
    Since $W^2(z) = (z-1)^{\alpha} (z+1)^{\beta} h(z)$, it then follows that
    \[
        \left(\frac{W(z)}{D(z)}\right)^2
        =\varphi(z)^{\alpha+\beta}
        \exp\left((z^2-1)^{1/2}
        \frac{1}{2\pi i}\int_{\gamma}\frac{\log h(\zeta)}{(\zeta^2-1)^{1/2}}
        \frac{d\zeta}{\zeta-z}\right).
    \]
    Since for $z \in U_{\delta}$,
    \[
        \frac{1}{2\pi i}\int_{\gamma}\frac{\log
        h(\zeta)}{(\zeta^2-1)^{1/2}}\frac{d\zeta}{\zeta-z}=\sum_{n=0}^{\infty}c_{n}(z-1)^{n}
    \]
    with coefficients $c_n$ given by formula (\ref{Definitiecn}),
    the lemma is proved.
\end{proof}

\begin{proposition}
    The matrix valued function $E_{n}(z)$ defined in {\rm (\ref{eq536})} is
    analytic in $U \setminus (-\infty,-1]$.
\end{proposition}

\begin{proof}
    We first show that $E_n$ has no jump across the interval
    $(-1,1)$. Let $x \in (-1,1)$.
    Since $N$ has the jump (\ref{RHPNb}) and $W_+(x) W_-(x) = w(x)$,
    we have
    \begin{eqnarray}
        \nonumber
        N_+(x) W_+(x)^{\sigma_3} & = & N_-(x)
            \begin{pmatrix}
                0 & w(x) \\
                -w(x)^{-1} & 0
            \end{pmatrix}
            \left( \frac{w(x)}{W_-(x)} \right)^{\sigma_3} \\[1ex]
        \nonumber
        & = & N_-(x)
        \begin{pmatrix}
            0 & W_-(x) \\
            -W_-(x)^{-1} & 0
        \end{pmatrix} \\[1ex]
        \label{JumpNWsigma3}
        & = & N_-(x) W_-^{\sigma_3}(x)
        \begin{pmatrix}
            0 & 1 \\
            -1 & 0
        \end{pmatrix}.
    \end{eqnarray}
    Furthermore, we have that $f(x)$ is negative, and
    taking fourth roots, we get $f_+(x)^{1/4} = i f_-(x)^{1/4}$,
    so that
    \begin{equation} \label{fplusfmin}
        f_+(x)^{\sigma_3/4} =
        \begin{pmatrix}
            i & 0 \\
            0 & -i
        \end{pmatrix}
        f_-(x)^{\sigma_3/4}.
    \end{equation}
    By (\ref{eq536}), (\ref{JumpNWsigma3}), and (\ref{fplusfmin})
    we then have for $x \in (-1,1)$,
    \begin{eqnarray*}
        (E_{n})_+(x) & = & N_+(x)W_+(x)^{\sigma_{3}}\frac{1}{\sqrt{2}}
            \begin{pmatrix}
                1 & -i\\
                -i & 1
            \end{pmatrix}
            f_+(x)^{\sigma_3/4} (2\pi n)^{\sigma_3/2}\\[1ex]
        & = & N_-(x) W_-(x)^{\sigma_3}
            \begin{pmatrix}
                0 & 1 \\
                -1 & 0
            \end{pmatrix}
            \frac{1}{\sqrt{2}}
            \begin{pmatrix}
                1 & -i \\
                -i & 1
            \end{pmatrix}
            \begin{pmatrix}
                i & 0 \\
                0 & -i
            \end{pmatrix}
            f_-(x)^{\sigma_3/4} (2\pi n)^{\sigma_3/2}\\[1ex]
        & = & N_{-}(x)W_{-}(x)^{\sigma_3} \frac{1}{\sqrt{2}}
            \begin{pmatrix}
                1 & -i\\
                -i & 1
            \end{pmatrix}
            f_-(x)^{\sigma_3/4} (2\pi n)^{\sigma_3/2}\\[1ex]
        & = & (E_n)_-(x).
    \end{eqnarray*}
    Therefore $E_n$ is analytic across $(-1,1)$.

    \medskip

    Next, we have to show that $E_n$ has a removable singularity at
    $1$. Using Lemma \ref{BehaviorWoverD} and (\ref{RHPNsolution}) we see that
    \begin{eqnarray}
        \nonumber
        N(z)W(z)^{\sigma_3} &=& D_{\infty}^{\sigma_3}
            \begin{pmatrix}
                \frac{a(z) + a(z)^{-1}}{2} & \frac{a(z) - a(z)^{-1}}{2i} \\[1ex]
                \frac{a(z) - a(z)^{-1}}{-2i} & \frac{a(z) + a(z)^{-1}}{2}
            \end{pmatrix}
            \left(\frac{W(z)}{D(z)} \right)^{\sigma_3} \\[1ex]
        \label{BehaviorNWSigma3}
        & = &
            O\begin{pmatrix}
                |z-1|^{-1/4} & |z-1|^{-1/4} \\
                |z-1|^{-1/4} & |z-1|^{-1/4}
            \end{pmatrix}
            O\begin{pmatrix}
                1 & 0 \\
                0 & 1
            \end{pmatrix}, \qquad \mbox{as $z\to 1$.}
    \end{eqnarray}
    Since $f$ has a simple zero at $1$ we also have
    \begin{equation}\label{BehaviorfSigma3over4}
        f(z)^{\sigma_{3}/4}=
        O\begin{pmatrix}
            |z-1|^{1/4} & 0 \\
            0 & |z-1|^{-1/4}
        \end{pmatrix}, \qquad \mbox{as $z\to 1$.}
    \end{equation}
    Thus by (\ref{eq536}), (\ref{BehaviorNWSigma3}), and
    (\ref{BehaviorfSigma3over4}) we have that all entries of $E_n(z)$
    have at most a square-root singularity at $z=1$. Since we already
    know that $E_n$ is analytic in a punctured neighborhood of $1$,
    the isolated singularity at $1$ is removable.
    Thus, $E_n$ is analytic in $U \setminus (-\infty,-1]$.
\end{proof}

The construction of the parametrix $P$ near $z=1$ is now
completed. A similar construction based on Bessel functions of
order $\beta$ yields a parametrix $\tilde{P}$ defined in a
neighborhood $\tilde{U}_{\delta}=\{z \mid |z+1|<\delta\}$ of
$z=-1$ that satisfies the following RHP.

\subsubsection*{RHP for \boldmath$\tilde{P}$:}

\begin{enumerate}
    \item[(a)]
        $\tilde{P}(z)$ is defined and analytic for $z \in \tilde{U}_{\delta_0} \setminus
        \Sigma$ for some $\delta_0 > \delta$.
    \item[(b)]
        $\tilde{P}(z)$ satisfies the following jump relations on
        $\tilde{U}_{\delta} \cap \Sigma^o$:
        \begin{equation}\label{RHPPTildeb1}
            \tilde{P}_{+}(z)=\tilde{P}_{-}(z)
            \begin{pmatrix}
                1 & 0\\
                w(z)^{-1}\varphi(z)^{-2n} & 1
            \end{pmatrix}, \qquad \mbox{for $z\in\tilde{U}_{\delta} \cap
            \left(\Sigma_{1}^o\cup\Sigma_{3}^o\right)$,}
        \end{equation}
        \begin{equation}\label{RHPPTildeb2}
            \tilde{P}_+(x) = \tilde{P}_-(x)
            \begin{pmatrix}
                0 & w(x) \\
                -w(x)^{-1} & 0
            \end{pmatrix}, \qquad \mbox{for $x \in \tilde{U}_{\delta} \cap \Sigma_2^o$.}
        \end{equation}
    \item[(c)]
        On $\partial\tilde{U}_{\delta}$ we have, as $n\to\infty$
        \begin{equation}\label{RHPPTildec}
            \tilde{P}(z) N^{-1}(z) = I + O \left( \frac{1}{n} \right),
            \qquad \mbox{uniformly for $ z\in\partial\tilde{U}_{\delta}\setminus\Sigma$.}
        \end{equation}
    \item[(d)]
        For $\beta<0$, the matrix function $\tilde{P}(z)$ has the following behavior as $z\to -1$:
        \[
            \tilde{P}(z)=
            O\begin{pmatrix}
                1 & |z+1|^{\beta} \\
                1 & |z+1|^{\beta}
            \end{pmatrix}, \qquad \mbox{as $z\to -1, z\in\tilde{U}_{\delta}\setminus\Sigma$.}
        \]
        For $\beta=0$, $\tilde{P}(z)$ has the following behavior as $z\to -1$:
        \[
            \tilde{P}(z)=
            O\begin{pmatrix}
                \log|z+1| & \log|z+1| \\
                \log|z+1| & \log|z+1|
            \end{pmatrix}, \qquad \mbox{as $z\to -1, z\in\tilde{U}_{\delta}\setminus \Sigma$.}
        \]
        For $\beta>0$, $\tilde{P}(z)$ has the following behavior as $z\to -1$:
        \[
            \tilde{P}(z)=
            \left\{\begin{array}{cl}
                O\begin{pmatrix}
                    1 & 1 \\
                    1 & 1
                \end{pmatrix}, & \mbox{as $z\to -1$ outside the lens,}\\[2ex]
                O\begin{pmatrix}
                    |z+1|^{-\beta} & 1 \\
                    |z+1|^{-\beta} & 1
                \end{pmatrix}, & \mbox{as $z\to -1$ inside the lens.}
            \end{array} \right.
        \]
\end{enumerate}

Details of the construction of $\tilde{P}$ are as follows. Instead
of working with $\varphi$ it is more convenient to introduce
\begin{equation} \label{eq545}
    \tilde{\varphi}(z) = \varphi(-z) = - \varphi(z),
\end{equation}
which is defined and analytic in $\mathbb C \setminus [-1,1]$.
Replacing $\varphi$ by $\tilde{\varphi}$ does not change the jump
matrix in (\ref{RHPPTildeb1}). Since $\varphi_+(x)\varphi_-(x)=1$
for $x\in(-1,1)$ it is obvious that we also have
\begin{equation}\label{eq545b}
    \tilde{\varphi}_+(x)\tilde{\varphi}_-(x)=1,\qquad\mbox{for }x\in(-1,1).
\end{equation}
The function
\begin{equation} \label{eq546}
    \tilde{W}(z) = \left((1-z)^{\alpha}(-1-z)^{\beta} h(z) \right)^{1/2},
    \qquad \mbox{for $z \in U  \setminus [-1,\infty)$,}
\end{equation}
is analytic in $U \setminus [-1,\infty)$, and satisfies
\begin{equation} \label{eq547}
    \tilde{W}^2(z) =
    \left\{ \begin{array}{ll}
        e^{-\beta \pi i} w(z), & \qquad\mbox{for $\Im z > 0$,} \\[1ex]
        e^{\beta \pi i} w(z), &  \qquad\mbox{for $\Im z < 0$.}
    \end{array} \right.
\end{equation}
Thus in particular, $\tilde{W}^{2}_{\pm}(x)=e^{\mp\beta\pi i}w(x)$
for $x\in(-1,1)$, and so
\begin{equation}\label{eq547b}
    \tilde{W}_{+}(x)\tilde{W}_{-}(x)=w(x),\qquad\mbox{for } x\in(-1,1).
\end{equation}

Writing
\begin{equation} \label{eq548}
   \tilde{P}(z) = \tilde{E}_n(z) \tilde{P}^{(1)}(z) \tilde{W}(z)^{-\sigma_3}
    \tilde{\varphi}(z)^{-n\sigma_3},
\end{equation}
with an invertible analytic factor $\tilde{E}_n$, we see, using
(\ref{RHPPTildeb1}), (\ref{RHPPTildeb2}), (\ref{eq545b}),
(\ref{eq547}) and (\ref{eq547b}), that $\tilde{P}^{(1)}$ should
satisfy the following RHP.

\subsubsection*{RHP for \boldmath$\tilde{P}^{(1)}$:}

\begin{enumerate}
    \item[(a)]
        $\tilde{P}^{(1)}(z)$ is defined and analytic for $z\in \tilde{U}_{\delta_0}\setminus\Sigma$
        for some $\delta_0 > \delta$.
    \item[(b)]
        $\tilde{P}^{(1)}(z)$ satisfies the following jump relations:
        \begin{equation}\label{RHPPTilde2b}
            \left\{\begin{array}{ll}
                \tilde{P}^{(1)}_{+}(z)=\tilde{P}^{(1)}_{-}(z)
                \begin{pmatrix}
                    1 & 0 \\
                    e^{-\beta\pi i} & 1
                \end{pmatrix},& \mbox{for $z\in\tilde{U}_{\delta}\cap\Sigma_{1}^o$,} \\[2ex]
                \tilde{P}^{(1)}_{+}(x)=\tilde{P}^{(1)}_{-}(x)
                \begin{pmatrix}
                    0 & 1 \\
                    -1 & 0
                \end{pmatrix}, & \mbox{for $x\in \tilde{U}_{\delta}\cap\Sigma_{2}^o$,} \\[2ex]
                \tilde{P}^{(1)}_{+}(z)=\tilde{P}^{(1)}_{-}(z)
                \begin{pmatrix}
                    1 & 0\\
                    e^{\beta\pi i} & 1
                \end{pmatrix}, & \mbox{for $z\in \tilde{U}_{\delta}\cap\Sigma_{3}^o$.}
            \end{array} \right.
        \end{equation}
    \item[(c)]
        For $\beta<0$, the matrix function $\tilde{P}^{(1)}(z)$ has the following behavior as $z\to -1$:
        \[
            \tilde{P}^{(1)}(z)=
            O\begin{pmatrix}
                |z+1|^{\beta/2} & |z+1|^{\beta/2} \\
                |z+1|^{\beta/2} & |z+1|^{\beta/2}
            \end{pmatrix}, \qquad\mbox{as $z\to -1, z\in \tilde{U}_{\delta}\setminus\Sigma$.}
        \]
        For $\beta=0$, $\tilde{P}^{(1)}(z)$ has the following behavior as $z\to -1$:
        \[
            \tilde{P}^{(1)}(z)=
            O\begin{pmatrix}
                \log|z+1| & \log|z+1| \\
                \log|z+1| & \log|z+1|
            \end{pmatrix}, \qquad\mbox{as $z\to -1, z\in\tilde{U}_{\delta}\setminus\Sigma$.}
        \]
        For $\beta>0$, $\tilde{P}^{(1)}(z)$ has the following behavior as $z\to -1$:
        \[
            \tilde{P}^{(1)}(z)=
            \left\{\begin{array}{cl}
                O\begin{pmatrix}
                    |z+1|^{\beta/2} & |z+1|^{-\beta/2} \\
                    |z+1|^{\beta/2} & |z+1|^{-\beta/2}
                \end{pmatrix}, & \mbox{as $z\to -1$ outside the lens,}\\[2ex]
                O\begin{pmatrix}
                    |z+1|^{-\beta/2} & |z+1|^{-\beta/2}\\
                    |z+1|^{-\beta/2} & |z+1|^{-\beta/2}
                \end{pmatrix}, & \mbox{as $z\to -1$ inside the lens.}
            \end{array}\right.
        \]
\end{enumerate}
The jumps (\ref{RHPPTilde2b}) of $\tilde{P}^{(1)}$ are analogous
to the jumps (\ref{RHPP1b}) of $P^{(1)}$. To construct
$\tilde{P}^{(1)}$ we first define a function $\tilde{g}$ as
\begin{equation}\label{defgTilde(z)}
    \tilde{g}(z)=\log\tilde{\varphi}(z),\qquad\mbox{for $z\in\mathbb{C}\setminus [-1,\infty)$.}
\end{equation}
and a function $\tilde{f}$ by
\begin{equation} \label{deffTilde(z)}
    \tilde{f}(z) = \tilde{g}^2(z)/4,
    \qquad \mbox{for $z \in \tilde{U}_{\delta} \setminus [-1,\infty)$.}
\end{equation}
Then $\zeta = \tilde{f}(z)$ is a one-to-one conformal  mapping of
$\tilde{U}_{\delta}$ onto a neighborhood of $\zeta = 0$. However
the ordering is reversed. Thus if $z < -1$, then $\zeta > 0$,
while if $z > -1$, then $\zeta < 0$. The contours $\Sigma_1$ and
$\Sigma_3$ are chosen such that $\tilde{f}(\Sigma_1 \cap \tilde{U}_{\delta})$
is part of the ray $\gamma_{3}$, and $\tilde{f}(\Sigma_{3} \cap \tilde{U}_{\delta})$
is part of the ray $\gamma_{1}$ in the $\zeta$-plane.

\begin{figure}
    \center{ \resizebox{8cm}{!}{\includegraphics{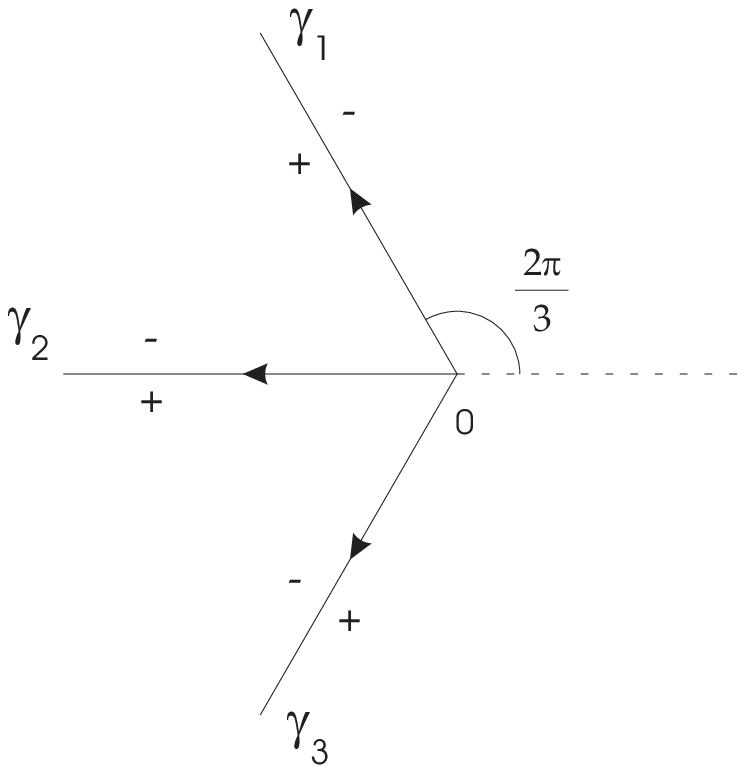}}
    \caption{The contour $\Sigma_{\tilde{\Psi}}$}\label{figuur5}}
\end{figure}

Then we need a fixed matrix function $\tilde{\Psi}$ defined in the
complex $\zeta$-plane that satisfies the following RHP
on the contour $\Sigma_{\tilde{\Psi}}$ shown in Figure \ref{figuur5}.

\subsubsection*{RHP for \boldmath$\tilde{\Psi}$:}

\begin{enumerate}
    \item[(a)]
        $\tilde{\Psi}(\zeta)$ is analytic for $\zeta \in
        \mathbb{C}\setminus\Sigma_{\tilde{\Psi}}$.
    \item[(b)]
        $\tilde{\Psi}(\zeta)$ satisfies the following jump relations:
        \begin{eqnarray*}
            \tilde{\Psi}_{+}(\zeta) &=& \tilde{\Psi}_{-}(\zeta)
            \begin{pmatrix}
                1 & 0 \\
                e^{\beta \pi i}  & 1
            \end{pmatrix}, \qquad\mbox{for $\zeta\in\gamma_{1}$,} \\[2ex]
            \tilde{\Psi}_+(\zeta) &=& \tilde{\Psi}_-(\zeta)
            \begin{pmatrix}
                0 & 1 \\
                -1  & 0
            \end{pmatrix}, \qquad\mbox{for $\zeta \in \gamma_{2}$,} \\[2ex]
            \tilde{\Psi}_+(\zeta) &=& \tilde{\Psi}_-(\zeta)
            \begin{pmatrix}
                1 & 0 \\
                e^{-\beta \pi i}  & 1
            \end{pmatrix},
            \qquad\mbox{for $\zeta\in\gamma_{3}$.}
        \end{eqnarray*}
    \item[(c)]
        $\tilde{\Psi}(\zeta)$ has the same behavior as $\Psi(\zeta)$
        near $\zeta = 0$, if we replace $\alpha$ by $\beta$.
\end{enumerate}
When comparing this RHP with (\ref{RHPPSIb1})--(\ref{RHPPSIb3}) we
have to realize that the {\em directions on the contours are
reversed}. The reversal of the directions is taken care of by
conjugating with the Pauli matrix $\sigma_3$. Thus we take
\begin{equation} \label{eq557}
    \tilde{\Psi}(\zeta) = \sigma_3 \Psi(\zeta; \alpha \to \beta) \sigma_3,
\end{equation}
where $\Psi(\cdot; \alpha \to \beta)$ is the matrix valued function defined by
(\ref{RHPPSIsolution1})--(\ref{RHPPSIsolution3}), but with Bessel
functions of order $\beta$ instead of order $\alpha$.
Then $\tilde{P}^{(1)}(z)=\tilde{\Psi}(n^2\tilde{f}(z))$ is a
solution of the RHP for $\tilde{P}^{(1)}$ and thus
\begin{equation} \label{eq558}
    \tilde{P}(z) = \tilde{E}_n(z) \tilde{\Psi}(n^2 \tilde{f}(z))
    \tilde{\varphi}(z)^{-n\sigma_3} \tilde{W}(z)^{-\sigma_3},
\end{equation}
with $\tilde{E}_n$ an analytic matrix valued function, satisfies
the parts (a), (b), and (d) of the RHP for $\tilde{P}$. As before
we construct the analytic factor $\tilde{E}_n$ so that the
matching condition (c) is satisfied as well. We define
$\tilde{E}_n$ by
\begin{equation} \label{eq559}
    \tilde{E}_n(z) =  N(z) \tilde{W}(z)^{\sigma_3} \frac{1}{\sqrt{2}}
    \begin{pmatrix}
        1 & i \\
        i & 1
    \end{pmatrix}
    \tilde{f}(z)^{\sigma_3/4} \left(2\pi n \right)^{\sigma_3/2},
    \quad \mbox{for $z \in \tilde{U}_{\delta}$.}
\end{equation}
Similar reasoning as before shows that $\tilde{E}_n$ is analytic
in $U \setminus [1,\infty)$, and that the matching condition is satisfied.

This completes the construction of the local parametrix
$\tilde{P}$ near $-1$. Note the slight difference (apart from the
tildes) between $\tilde{E}_n$ defined in (\ref{eq559}) and $E_n$
defined in (\ref{eq536}).

\medskip

For later use, we remark that the following analogue of Lemma
\ref{BehaviorWoverD} holds. Its proof is similar.

\begin{lemma}\label{BehaviorWTildeoverD}
    For $z \in \tilde{U}_{\delta}$, we have
    \[
        \left(\frac{\tilde{W}(z)}{D(z)}\right)^2=\tilde{\varphi}(z)^{\alpha+\beta}
        \exp\left((z^2-1)^{1/2} \sum_{n=0}^{\infty}d_{n}(z+1)^{n}\right),
    \]
    with coefficients $d_n$ given by {\rm (\ref{Definitiedn})}.
\end{lemma}

    \section{Third transformation $S \mapsto R$}
\setcounter{equation}{0}

In the final transformation of our original problem we define the
matrix valued function $R$ by
\begin{eqnarray}
    \label{DefRBuiten}
    R(z) &=& S(z) N^{-1}(z),\qquad \mbox{for $z\in\mathbb{C}\setminus(\overline{U}_{\delta}\cup\overline{\tilde{U}}_{\delta}
        \cup \Sigma)$,} \\[1ex]
    \label{DefRU}
    R(z) &=& S(z) P^{-1}(z),\qquad \mbox{for $z \in U_{\delta} \setminus \Sigma$,} \\[1ex]
    \label{DefRUTilde}
    R(z) &=& S(z) \tilde{P}^{-1}(z),\qquad\mbox{for $z \in \tilde{U}_{\delta} \setminus\Sigma$.}
\end{eqnarray}

\begin{remark}
    The inverses of the matrices $N$, $P$, and $\tilde{P}$
    used in (\ref{DefRBuiten})--(\ref{DefRUTilde}) exist, since
    the determinants of these matrices are equal to $1$.
    For $N$, this is easy to see from its definition (\ref{RHPNsolution}).
    To see this for $P$, we note that by (\ref{RHPPsolution}) and (\ref{detEn})
    it is sufficient to show that $\det \Psi(\zeta) = 1$ for
    $\zeta \in \mathbb C\setminus \Sigma_{\Psi}$.
    It is clear that $\det \Psi(\zeta)$ is analytic in
    $\mathbb C \setminus \Sigma_{\Psi}$.
    Since the jump matrices in (\ref{RHPPSIb2})--(\ref{RHPPSIb3}) have determinant
    one, it follows that $\det \Psi(\zeta)$ has an analytic continuation to
    $\mathbb{C}\setminus\{0\}$.
    By (\ref{RHPPSIc1})--(\ref{RHPPSIc3}), the behavior of
    $\det \Psi(\zeta)$ as $\zeta \to 0$ is
    \[
        \det\Psi(\zeta)=
        \left\{\begin{array}{cl}
            O\left(|\zeta|^{\alpha}\right), & \mbox{if $\alpha<0$,} \\[2ex]
            O\left((\log|\zeta|)^2 \right), & \mbox{if $\alpha=0$,} \\[2ex]
            O(1), & \mbox{if $\alpha>0$ and $|\arg(\zeta)| < 2\pi/3$,} \\[1ex]
            O(|\zeta|^{-\alpha}) & \mbox{if $\alpha>0$ and $2\pi/3 < |\arg(\zeta)| < \pi$.}
        \end{array} \right.
    \]
    Then the singularity at $0$ is removable, which is obvious in case $\alpha \leq 0$,
    and which follows as in the proof of Lemma \ref{TheoremSecondTransform} in the
    paragraph after formula (\ref{Halpha3}) in case $\alpha > 0$.
    Thus $\det\Psi(\zeta)$ is entire.
    From the asymptotics at infinity, given in (\ref{AsymptoticsPSI}),
    we see that $\det \Psi(\zeta) \to 1$ as $\zeta\to \infty$. By Liouville's theorem
    we then find that $\det\Psi(\zeta)=1$ for every $\zeta$. As already noted
    this implies that $\det P(z) = 1$.
    In the same way, it follows that $\det \tilde{P}(z) = 1$.
\end{remark}

\begin{figure}
    \center{ \resizebox{12cm}{!}{\includegraphics{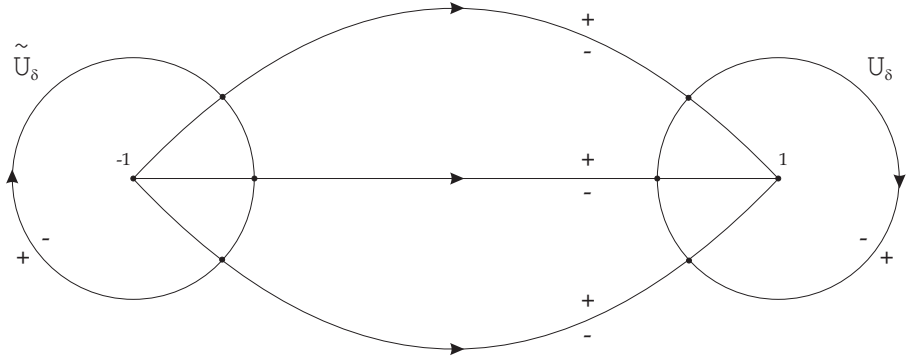}}
    \caption{System of contours for ${R}$}\label{figuur6}}
\end{figure}

\begin{figure}
    \center{ \resizebox{12cm}{!}{\includegraphics{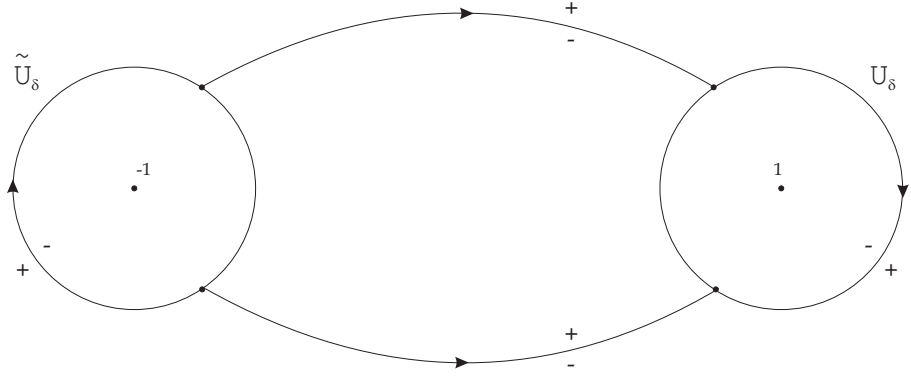}}
    \caption{Reduced system of contours $\Sigma_{R}$}\label{figuur7}}
\end{figure}

Now $R$ given by (\ref{DefRBuiten})--(\ref{DefRUTilde})
is defined and analytic in $\mathbb C \setminus
(\Sigma \cup \partial U_{\delta} \cup \partial \tilde{U}_{\delta})$.
Then $R$ has jumps on a system of contours shown in Figure
\ref{figuur6}.
However, since by construction the jumps of $S$ and $N$ on the
interval $(-1+\delta, 1-\delta)$ agree, it follows from
(\ref{DefRBuiten}) that $R$ has the identity as jump matrix on
this interval, and so has an analytic continuation to $(-1+\delta,
1-\delta)$. Also, since the jumps of $S$ and $P$ on
$U_{\delta}\cap\Sigma^o$ are exactly the same, it follows from
(\ref{DefRU}) that $R$ has an analytic continuation across these
contours as well. Similarly, it follows from (\ref{DefRUTilde})
that $R$ is analytic across $\tilde{U}_{\delta}\cap\Sigma^o$. What
remains are the jumps on a reduced system of contours, which we
will call $\Sigma_{R}$, see Figure \ref{figuur7}, and possible
isolated singularities at $\pm 1$. We prove that these
singularities are removable.

If $\alpha<0$, we use the fact that $\det P(z) = 1$, and
(\ref{RHPPd1}) to conclude that
\[
    P^{-1}(z) =
    O\begin{pmatrix}
        |z-1|^{\alpha} & |z-1|^{\alpha} \\
        1 & 1
    \end{pmatrix}, \qquad \mbox{as $z \to 1$.}
\]
Then from (\ref{RHPSd1}) and (\ref{DefRU}) it follows that
\[
    R(z)=
    O\begin{pmatrix}
        |z-1|^{\alpha} & |z-1|^{\alpha} \\
        |z-1|^{\alpha} & |z-1|^{\alpha}
    \end{pmatrix}, \qquad \mbox{as $z\to 1$.}
\]
Since $R$ is analytic in $U_{\delta} \setminus \{1\}$ and $\alpha > -1$,
it follows  that the singularity at $1$ is removable.

If $\alpha=0$, it follows in a similar way from (\ref{RHPSd2}),
(\ref{RHPPd2}), and (\ref{DefRU}) that
\[
    R(z)=
    O\begin{pmatrix}
        \left(\log|z-1|\right)^{2} & \left(\log|z-1|\right)^{2} \\
        \left(\log|z-1|\right)^{2} & \left(\log|z-1|\right)^{2}
    \end{pmatrix}, \qquad \mbox{as $z\to 1$,}
\]
from which it also follows that the singularity at $1$ is removable.

And finally, if $\alpha>0$, it follows from (\ref{RHPSd3}),
(\ref{RHPPd3}), and (\ref{DefRU}) that
\[
    R(z)=
    \left\{\begin{array}{cl}
        O\begin{pmatrix}
            1 & 1 \\
            1 & 1
        \end{pmatrix}, & \mbox{as $z\to 1$ outside the lens,} \\[2ex]
        O\begin{pmatrix}
            |z-1|^{-\alpha} & |z-1|^{-\alpha} \\
            |z-1|^{-\alpha} & |z-1|^{-\alpha}
        \end{pmatrix}, & \mbox{as $z\to 1$ inside the lens.}
    \end{array}\right.
\]
From this it also follows that the singularity at $1$ is removable,
cf.\ the argument given after
formula (\ref{Halpha3}) in Lemma \ref{TheoremSecondTransform}.

\medskip

So in all cases we proved that the singularity at $1$ is removable.
Similarly it follows  that the singularity at $-1$ is removable.
Hence $R$ is analytic in $\mathbb{C}\setminus\Sigma_{R}$, where
\[\Sigma_R = \partial U_{\delta} \cup \partial \tilde{U}_{\delta}
    \cup \left((\Sigma_1 \cup \Sigma_3) \setminus (U_{\delta} \cup \tilde{U}_{\delta})\right). \]
The contour $\Sigma_R$ is oriented as indicated in Figure \ref{figuur7}.

\medskip

We next calculate the jumps of $R$. For $z \in \Sigma_R \setminus
(\partial U_{\delta} \cup \partial \tilde{U}_{\delta})$, we use (\ref{RHPSb1}) and
(\ref{DefRBuiten}) to obtain
\begin{eqnarray*}
    R_+(z) & = & S_+(z) N^{-1}(z) \ = \ S_-(z)
    \begin{pmatrix}
        1 & 0 \\
        w(z)^{-1} \varphi(z)^{-2n} & 1
    \end{pmatrix} N^{-1}(z)  \\[1ex]
    & = & R_-(z) N(z)
    \begin{pmatrix}
        1 & 0 \\
        w(z)^{-1}\varphi(z)^{-2n} & 1
    \end{pmatrix} N^{-1}(z).
\end{eqnarray*}
On the boundary of $U_{\delta}$ we have $R_{+}(z)=S(z)N^{-1}(z)$
and $R_-(z) = S(z) P^{-1}(z)$ so that
\[
    R_{+}(z)=R_{-}(z)P(z)N^{-1}(z), \qquad \mbox{for $z\in\partial U_{\delta}$.}
\]
Similarly, on the boundary of $\tilde{U}_{\delta}$ we have
\[
    R_{+}(z)=R_{-}(z)\tilde{P}(z)N^{-1}(z),\qquad\mbox{for $z\in\partial\tilde{U}_{\delta}$.}
\]
So we see that $R$ is a solution of the following RHP:

\subsubsection*{RHP for \boldmath$R$:}

\begin{enumerate}
    \item[(a)]
        $R(z)$ is analytic for $z \in \mathbb C \setminus \Sigma_R$.
    \item[(b)]
        $R(z)$ satisfies the following jump relations on $\Sigma_R$:
        \begin{eqnarray}\label{RHPRb1}
            R_{+}(z) &=& R_{-}(z)P(z) N^{-1}(z), \qquad \mbox{for $z\in \partial U_{\delta}$,} \\[2ex]
            \label{RHPRb2}
            R_{+}(z) &=& R_{-}(z)\tilde{P}(z) N^{-1}(z), \qquad \mbox{for $z \in \partial \tilde{U}_{\delta}$,} \\[2ex]
            \label{RHPRb3}
            R_{+}(z) &=& R_{-}(z)N(z)
            \begin{pmatrix}
                1 & 0 \\
                w(z)^{-1}\varphi(z)^{-2n} & 1
            \end{pmatrix}N^{-1}(z), \\[1ex]
            & & \qquad \qquad \qquad \qquad \qquad \nonumber
            \mbox{ for } z \in \Sigma_R \setminus (\partial U_{\delta} \cup
            \partial \tilde{U}_{\delta}).
        \end{eqnarray}
    \item[(c)]
        $R(z) =  I + \begin{displaystyle}O \left( \frac{1}{z}
        \right)\end{displaystyle}, \qquad \mbox{as } z \to \infty$.
\end{enumerate}

\begin{remark}
    Since $\log|\varphi(z)| > c > 0$ on the upper and lower lips of
    $\Sigma_{R}$ and $N(z)$ does not depend on $n$, the jump matrix in
    (\ref{RHPRb3}) tends to the identity matrix at an exponential rate
    \begin{equation} \label{AsympR1}
        R_+(z) = R_-(z) \left(I + O\left(e^{-2cn}\right)\right), \qquad \mbox{as $n\to\infty$}
    \end{equation}
    uniformly on $\Sigma_R \setminus (\partial U_{\delta} \cup \partial \tilde{U}_{\delta})$.
    By (\ref{RHPPc}) and (\ref{RHPRb1}) we have on the boundary of
    $U_{\delta}$ as $n\to\infty$ that
    \begin{equation} \label{AsympR2}
        R_{+}(z) = R_{-}(z) \left(I +
        O\left(\frac{1}{n}\right)\right),\qquad \mbox{as $n\to\infty$,}
    \end{equation}
    uniformly on $\partial U_{\delta}$. At the two points where
    $\partial U_{\delta}$ intersects with the arcs $\Sigma_1$ and $\Sigma_3$,
    the boundary values are
    taken within a single connected component of $\mathbb C \setminus \Sigma_R$.
    Since $R$ is continuous up to the boundary from within each component,
    (\ref{AsympR2}) holds everywhere on $\partial U_{\delta}$. Similarly by
    (\ref{RHPPTildec}) and (\ref{RHPRb2})
    \begin{equation} \label{AsympR3}
        R_{+}(z) = R_{-}(z) \left(I + O\left(\frac{1}{n} \right)
        \right),\qquad \mbox{as $n\to\infty$,}
    \end{equation}
    uniformly on $\partial \tilde{U}_{\delta}$.
    From
    (\ref{AsympR1})--(\ref{AsympR3}) it follows as in \cite{DKMVZ1} that
    \begin{equation} \label{AsympR}
        R(z) = I + O\left(\frac{1}{n} \right),
        \qquad \mbox{as $n\to\infty$ uniformly on $\mathbb C \setminus \Sigma_R$,}
    \end{equation}
    The relation (\ref{AsympR}) is the main term in
    the asymptotics for $R$ and it is enough to give the leading term in the
    asymptotics of the orthogonal polynomials. In the next section we
    will give a complete asymptotic expansion for $R$.
\end{remark}

    \section{Asymptotics for \boldmath$R(z)$}
\setcounter{equation}{0}

\subsection{Asymptotic expansion for $\Delta$}
We write the jump matrix for $R$ in the form $I + \Delta(s)$, so
that by (\ref{RHPRb1})--(\ref{RHPRb3})
\begin{equation} \label{DefDelta}
    I + \Delta(s) =
    \left\{\begin{array}{cl}
        P(s) N^{-1}(s), & \mbox{for $s\in \partial U_{\delta}$,} \\[2ex]
        \tilde{P}(s) N^{-1}(s), & \mbox{for $s \in \partial \tilde{U}_{\delta}$,} \\[2ex]
        N(s)
        \begin{pmatrix}
            1 & 0\\
            w^{-1}(s)\varphi(s)^{-2n} & 1
        \end{pmatrix}N^{-1}(s), & \mbox{for $s \in \Sigma_R \setminus
        (\partial U_{\delta}\cup\partial \tilde{U}_{\delta})$.}
    \end{array} \right.
\end{equation}
We show that $\Delta$ has  an asymptotic expansion in powers of $1/n$
of the form
\begin{equation} \label{AsympOntDelta}
    \Delta(s) \sim \sum_{k=1}^{\infty} \frac{\Delta_k(s)}{n^k}
    \qquad \mbox{as }n \to \infty,
\end{equation}
uniformly for $s \in \Sigma_R$.
On  $\Sigma_R \setminus (\partial U_{\delta} \cup \partial \tilde{U}_{\delta})$, the jump
matrix is the identity matrix plus an exponentially small term,
so that for every $k \in \mathbb N$,
\begin{equation}\label{Deltak=0}
    \Delta_k(s)=0, \qquad\mbox{for $s\in\Sigma_{R}\setminus\left(\partial
    U_{\delta}\cup\partial\tilde{U}_{\delta}\right)$.}
\end{equation}

We now determine $\Delta_k(s)$ for $s\in\partial U_{\delta}$.
The modified Bessel functions $I_{\alpha}(\zeta)$, $K_{\alpha}(\zeta)$,
and their derivatives, have an asymptotic expansion as $\zeta \to \infty$,
see formulas 9.7.1--9.7.4 of \cite{AbramowitzStegun}.
Using this expansion in (\ref{RHPPSIsolution1}), we find that
\begin{eqnarray}
    \nonumber
    \Psi(\zeta)  & = & (2\pi)^{-\sigma_3/2} \zeta^{-\sigma_3/4} \frac{1}{\sqrt{2}}
        \begin{pmatrix}
            1 & i \\
            i & 1
        \end{pmatrix} \\[1ex]
    \label{DeltaAsymptotics}
    & & \left[ I + \sum_{k=1}^{\infty} \frac{(\alpha,k-1)}{4^k \zeta^{k/2}}
    \begin{pmatrix}
        \frac{(-1)^k}{k} (\alpha^2 + \frac{1}{2}k-\frac{1}{4}) & -(k-\frac{1}{2})i \\[1ex]
        (-1)^k (k-\frac{1}{2}) i & \frac{1}{k} (\alpha^2 + \frac{1}{2}k - \frac{1}{4})
    \end{pmatrix}\right] e^{2 \zeta^{1/2}\sigma_3}
\end{eqnarray}
as $\zeta \to \infty$, uniformly in the sector $|\arg \zeta| < 2 \pi/3$.
Here $(\alpha,0)=1$ and
\[
    (\alpha,k) = \frac{(4\alpha^2-1)(4\alpha^2-9) \cdots (4\alpha^2 - (2k-1)^2)}{2^{2k} k!},
\]
see also e.g.\ \cite[\S4.9,\S4.12]{AAR}.
Inserting the asymptotic expansions of the Hankel functions into
(\ref{RHPPSIsolution2}) and (\ref{RHPPSIsolution3}), we obtain the
same formula (\ref{DeltaAsymptotics}) as $\zeta \to \infty$ in the
other sectors. Thus (\ref{DeltaAsymptotics}) is valid uniformly as
$\zeta \to \infty$.
Then by a calculation using (\ref{defg(z)}), (\ref{deff(z)}),
(\ref{RHPPsolution}), (\ref{eq536}), and (\ref{DefDelta}), we find
that  uniformly for $s \in \partial U_{\delta}$,
\begin{eqnarray*}
    \Delta(s) & = & P(s) N^{-1}(s) - I \\[2ex]
    & \sim & N(s) W(s)^{\sigma_3} \left[ \sum_{k=1}^{\infty} \frac{(\alpha,k-1)}{n^k 2^k g(s)^k}
    \begin{pmatrix}
        \frac{(-1)^k}{k} (\alpha^2 + \frac{1}{2}k-\frac{1}{4}) & -(k-\frac{1}{2})i \\[1ex]
        (-1)^k (k-\frac{1}{2}) i & \frac{1}{k} (\alpha^2 + \frac{1}{2}k - \frac{1}{4})
        \end{pmatrix}\right] \\[1ex]
    & & \ W(s)^{-\sigma_3} N^{-1}(s).
\end{eqnarray*}
Thus,
\begin{eqnarray}
    \nonumber
    \Delta_k(s) & = & \frac{(\alpha,k-1)}{2^k g(s)^k} N(s) W(s)^{\sigma_3}
    \begin{pmatrix}
        \frac{(-1)^k}{k} (\alpha^2 + \frac{1}{2}k-\frac{1}{4}) & -(k-\frac{1}{2})i \\[1ex]
        (-1)^k (k-\frac{1}{2}) i & \frac{1}{k} (\alpha^2 + \frac{1}{2}k - \frac{1}{4})
    \end{pmatrix} \\[1ex]
    \label{DefinitieDeltak}
    & & \qquad W(s)^{-\sigma_3} N^{-1}(s),
    \qquad \mbox{for $s\in\partial U_{\delta}$,}
\end{eqnarray}
Analogously, using (\ref{defgTilde(z)})--(\ref{eq559}), and
(\ref{DefDelta}), we find
\begin{eqnarray}
    \nonumber
    \Delta_k(s) & = & \frac{(\beta,k-1)}{2^k \tilde{g}(s)^k} N(s) \tilde{W}(s)^{\sigma_3}
    \begin{pmatrix}
        \frac{(-1)^k}{k} (\beta^2 + \frac{1}{2}k-\frac{1}{4}) & (k-\frac{1}{2})i \\[1ex]
        (-1)^{k+1} (k-\frac{1}{2}) i & \frac{1}{k} (\alpha^2 + \frac{1}{2}k - \frac{1}{4})
    \end{pmatrix} \\[1ex]
    \label{DefinitieDeltakTilde}
    & & \qquad\tilde{W}(s)^{-\sigma_3} N^{-1}(s),
    \qquad \mbox{for $s\in\partial \tilde{U}_{\delta}$.}
\end{eqnarray}

To summarize, we have shown
\begin{lemma}
    $\Delta(s)$ has an asymptotic expansion in powers of $1/n$ of the
    form {\rm(\ref{AsympOntDelta})}, which is valid uniformly for $s \in
    \Sigma_R$. For $s \in \Sigma_R \setminus (\partial U_{\delta} \cup
    \partial \tilde{U}_{\delta})$, we have $\Delta_k(s) = 0$ for every $k$. For
    $s \in \partial U_{\delta}$ and $s \in \partial
    \tilde{U}_{\delta}$, we have that $\Delta_k$ is given by
    {\rm(\ref{DefinitieDeltak})--(\ref{DefinitieDeltakTilde})},
    respectively.
\end{lemma}

For later use, we need
\begin{lemma}
    \begin{enumerate}
        \item[\rm (a)] For every $k \in \mathbb N$, the restriction of $\Delta_k$
            to $\partial U_{\delta}$ has a meromorphic continuation to a neighborhood
            $U_{\delta_0}$ of $1$ for some $\delta_0 > \delta$. The continuation is analytic,
            except at $1$, where it has a pole of order at most $[(k+1)/2]$.
        \item[\rm (b)] For every $k \in \mathbb N$, the restriction of $\Delta_k$
            to $\partial \tilde{U}_{\delta}$ has a meromorphic continuation to a neighborhood
            $\tilde{U}_{\delta_0}$ of $-1$ for some $\delta_0 > \delta$. The continuation is analytic,
            except at $-1$, where it has a pole of order at most $[(k+1)/2]$.
    \end{enumerate}
\end{lemma}

\begin{proof}
    (a) The formula (\ref{DefinitieDeltak}) shows that $\Delta_k$ is
        analytic on $U \setminus (-\infty,1]$.
        We recall, see (\ref{JumpNWsigma3}), that for $x \in (-1,1)$,
        \[
            N_+(x)W_+(x)^{\sigma_3} = N_-(x)W_-(x)^{\sigma_3}
            \begin{pmatrix}
                0 & 1 \\
                -1 & 0
            \end{pmatrix}.
        \]
        Since $g_+(x) = - g_-(x)$, we then have for $x \in (-1,1)$,
        \begin{eqnarray*}
            \Delta_{k,+}(x) & = & (-1)^k \frac{(\alpha,k-1)}{2^k g_-(x)^k} N_-(x)W_-(x)^{\sigma_3} \\[1ex]
            & & \qquad
            \begin{pmatrix}
                0 & 1 \\
                -1 & 0
            \end{pmatrix}
            \begin{pmatrix}
                \frac{(-1)^k}{k} (\alpha^2 + \frac{1}{2}k-\frac{1}{4}) & -(k-\frac{1}{2})i \\[1ex]
                (-1)^{k} (k-\frac{1}{2}) i & \frac{1}{k} (\alpha^2 + \frac{1}{2}k - \frac{1}{4})
            \end{pmatrix}
            \begin{pmatrix}
                0 & -1\\
                1 & 0
            \end{pmatrix}
            \\[1ex]
            & & \qquad
            W_-(x)^{-\sigma_3} N_-^{-1}(x) \\[2ex]
            & = & (-1)^k \frac{(\alpha,k-1)}{2^k g_-(x)^k} N_-(x)W_-(x)^{\sigma_3} \\[1ex]
            & & \qquad (-1)^k
            \begin{pmatrix}
                \frac{(-1)^k}{k} (\alpha^2 + \frac{1}{2}k-\frac{1}{4}) & -(k-\frac{1}{2})i \\[1ex]
                (-1)^{k} (k-\frac{1}{2}) i & \frac{1}{k} (\alpha^2 + \frac{1}{2}k - \frac{1}{4})
            \end{pmatrix} \\[1ex]
            & & \qquad
            W_-(x)^{-\sigma_3} N_-^{-1}(x) \\[2ex]
            & = & \Delta_{k,-}(x).
        \end{eqnarray*}
        Thus $\Delta_k$ is analytic across $(-1,1)$.

        To determine the behavior near $1$, we note that $g(z) = O(|z-1|^{1/2})$ as $z \to 1$, and
        \[
            N(z)W(z)^{\sigma_3} =
            O\begin{pmatrix}
                |z-1|^{-1/4} & |z-1|^{-1/4} \\
                |z-1|^{-1/4} & |z-1|^{-1/4}
            \end{pmatrix},
        \]
        as $z \to 1$, see (\ref{BehaviorNWSigma3}).
        Thus we have by formula (\ref{DefinitieDeltak}),
        \[
            \Delta_k(z) =
            O\begin{pmatrix}
                |z-1|^{-\frac{k}{2} - \frac{1}{2}} & |z-1|^{-\frac{k}{2} - \frac{1}{2}} \\
                |z-1|^{-\frac{k}{2} - \frac{1}{2}} & |z-1|^{-\frac{k}{2} - \frac{1}{2}}
            \end{pmatrix}, \qquad \mbox{as $z \to 1$.}
        \]
        Since $z=1$ is an isolated singularity of $\Delta_k(z)$, it then follows
        that $z=1$ is a pole of order at most $[(k+1)/2]$. This proves part (a).

        \medskip

        Part (b) is proved similarly.
\end{proof}

\subsection{Asymptotic expansion for $R$}

As in Theorem 7.10 of \cite{DKMVZ1} we obtain from
(\ref{AsympOntDelta}) that $R$ has an asymptotic expansion of the
form
\begin{equation} \label{AsympOntR}
    R(z) \sim I + \sum_{k=1}^{\infty} \frac{R_k(z)}{n^k}, \qquad
    \mbox{as $n\to\infty$,}
\end{equation}
which is valid uniformly for
$z \in \mathbb C \setminus (\partial U_{\delta} \cup \partial \tilde{U}_{\delta})$.
We have that
\begin{equation} \label{AnalyticRk}
    R_k(z) \mbox{ is analytic for }
    z \in \mathbb C \setminus (\partial U_{\delta} \cup \partial \tilde{U}_{\delta}),
\end{equation}
and
\begin{equation} \label{AsympRk}
    R_k(z) = O\left( \frac{1}{z} \right), \qquad \mbox{as $z \to \infty$.}
\end{equation}

The expansion (\ref{AsympOntR}) is valid uniformly near infinity in the following sense.
\begin{lemma}\label{LemmaAsympR}
    For every $l$ there exists a constant $C > 0$, such that for every
    $z$ with  $|z|\geq 2$,
    \begin{equation} \label{AsymptOntRopOneindig}
        \left\|R(z)-I-\sum_{k=1}^{l}\frac{R_k(z)}{n^k}\right\|\leq\frac{C}{|z|n^{l+1}}.
    \end{equation}
    Here $\| \cdot \|$ denotes any matrix norm.
\end{lemma}

\begin{proof}
    From $R_+ = R_- + R_-\Delta$ we obtain the following integral representation for $R$
    (cf.\ the discussion in the proof of Theorem 7.8 in \cite{DKMVZ1}),
    \begin{equation}\label{IntegraalVoorstellingR}
        R(z)=I+\frac{1}{2\pi
        i}\int_{\Sigma_R}\frac{R_-(s)\Delta(s)}{s-z}\ ds, \qquad
        \mbox{for $z\in\mathbb{C}\setminus\Sigma_R$.}
    \end{equation}
    Expanding the jump relation $R_+=R_-(I+\Delta)$ using (\ref{AsympOntDelta})
    and (\ref{AsympOntR}), and collecting powers of $1/n$,
    we find for  $k=1,2,\ldots$,
    \begin{equation}\label{SprongRk}
        R_{k,+}(s)-R_{k,-}(s)=\sum_{j=1}^{k}R_{k-j,-}(s)\Delta_j(s), \qquad
        \mbox{for $s\in\partial U_{\delta} \cup \partial
        \tilde{U}_{\delta}$,}
    \end{equation}
    where $R_{0,-}=I$. From the Sokhotskii-Plemelj formula, (\ref{Deltak=0}),
    and (\ref{AsympRk}) it then follows that
    \begin{equation}\label{IntegraalVoorstellingRk}
        R_k(z)=\frac{1}{2\pi i}\int_{\Sigma_R}\sum_{j=1}^{k}R_{k-j,-}(s)\Delta_j(s)\frac{ds}{s-z},
        \qquad \mbox{for $z\in\mathbb{C}\setminus\Sigma_R$.}
    \end{equation}
    Combining (\ref{IntegraalVoorstellingR}) and
    (\ref{IntegraalVoorstellingRk}) we then have for $z \in \mathbb C \setminus \Sigma_R$,
    \[
        R(z)-I-\sum_{k=1}^{l}\frac{R_k(z)}{n^k}=\frac{1}{2\pi
        i}\int_{\Sigma_R}\left(R_-(s)\Delta(s)-\sum_{k=1}^l\sum_{j=1}^k\frac{R_{k-j,-}(s)\Delta_j
        (s)}{n^k}\right)\frac{ds}{s-z}.
    \]
    One easily verifies that we can rewrite the term in the parenthesis as follows
    \begin{eqnarray}
        \nonumber
        \lefteqn{ R_-(s)\Delta(s)-\sum_{k=1}^l\sum_{j=1}^k\frac{R_{k-j,-}(s)\Delta_j
            (s)}{n^k}\ = \ (R_{-}(s)-I)\left(\Delta(s)-\sum_{k=1}^{l-1}\frac{\Delta_k(s)}{n^k}\right)} \\[1ex]
        \label{Parenthesis}
        & & + \left(\Delta(s)-\sum_{k=1}^l \frac{\Delta_k(s)}{n^k}\right)
            + \sum_{k=1}^{l-1}\frac{\Delta_k(s)}{n^k}\left(R_{-}(s)-I-\sum_{j=1}^{l-k}\frac{R_{j,-}(s)}{n^{j}}\right).
    \end{eqnarray}
    From (\ref{AsympOntDelta}), (\ref{AsympOntR}), and the fact that each $\Delta_j$
    is bounded on $\Sigma_R$, it follows that each term in (\ref{Parenthesis}) is
    $O\left(n^{-(l+1)}\right)$ as $n\to\infty$, uniformly for
    $s\in\Sigma_R$. Hence the norm of the left-hand side of (\ref{Parenthesis})
    is $\leq \ds \frac{c}{n^{l+1}}$ for some constant $c$ uniformly for $s \in \Sigma_R$,
    and so
    \[
    \left\| R(z)-I-\sum_{k=1}^{l}\frac{R_k(z)}{n^k} \right\|
    \leq \frac{c}{n^{l+1}} \frac{1}{2\pi} \int_{\Sigma_R} \frac{|ds|}{|s-z|}.
    \]
    Since the contour $\Sigma_R$ is bounded and $|s-z|$ is uniformly bounded
    away from $0$ for $s \in \Sigma_R$ and $|z| \geq 2$, we obtain (\ref{AsymptOntRopOneindig})
    for some constant $C$. This completes the proof of the lemma.
\end{proof}

The functions $R_k$ can be computed from the known form
(\ref{AsympOntDelta}), (\ref{DefinitieDeltak}) and
(\ref{DefinitieDeltakTilde}) of the jump matrix. We will do the
calculations explicitly for the first two functions $R_1(z)$ and
$R_2(z)$. It will then be clear how to proceed for the further
terms.

\subsubsection{Determination of $R_1$}
Expanding the jump relation $R_+ = R_-(I + \Delta)$ up to order
$1/n$ we find, cf. (\ref{SprongRk})
\begin{equation} \label{JumpR1}
    R_{1,+}(s) - R_{1,-}(s) = \Delta_1(s),
    \qquad \mbox{for } s \in \partial U_{\delta}\cup\partial
    \tilde{U}_{\delta}.
\end{equation}
Then (\ref{JumpR1}) and (\ref{AnalyticRk})--(\ref{AsympRk}) with $k=1$,
give an additive RHP for $R_1$. This RHP can be solved by the
Sokhotskii-Plemelj formulas, but in this case we can state the
solution explicitly as follows.
From Lemma 8.1 we know that $\Delta_1(z)$ is analytic in a neighborhood of
$z=-1$ and $z=1$ except for simple poles at $\pm 1$. So
\[  \Delta_1(z)=\frac{A^{(1)}}{z-1}+O(1),\quad\mbox{as } z\to 1, \qquad
    \Delta_1(z)=\frac{B^{(1)}}{z+1}+O(1),\quad\mbox{as } z\to -1,
\]
for certain constant matrices $A^{(1)}$ and $B^{(1)}$.
By inspection we then see that
\begin{equation}\label{DefinitieR1}
R_1(z)= \left\{\begin{array}{ll}
     \ds \frac{A^{(1)}}{z-1}+\frac{B^{(1)}}{z+1},
     & \mbox{ for }z\in\mathbb{C}\setminus
     \left(\overline{U}_\delta\cup\overline{\tilde{U}}_\delta\right),\\[2ex]
     \ds \frac{A^{(1)}}{z-1}+\frac{B^{(1)}}{z+1}-\Delta_1(z),
     & \mbox{ for } z\in  U_\delta\cup\tilde{U}_\delta.
\end{array}\right.
\end{equation}
solves the additive RHP (\ref{AnalyticRk}), (\ref{AsympRk}), (\ref{JumpR1}).

\medskip

The constant matrix $A^{(1)}$ can be calculated
from the formula (\ref{DefinitieDeltak}) with $k=1$,
and the known behavior of $g$, $N$, and $W$ near $1$.
The behavior of $g$ is given by (\ref{Behaviorg}), and for
$N$ and $W$ we have as in  (\ref{BehaviorNWSigma3})
\begin{equation} \label{BehaviorNWsigma3}
    N(z) W(z)^{\sigma_3} = D_{\infty}^{\sigma_3}
   \begin{pmatrix}
    \frac{a(z) + a(z)^{-1}}{2} &
        \frac{a(z) - a(z)^{-1}}{2i} \\[1ex]
    \frac{a(z) - a(z)^{-1}}{-2i} &
    \frac{a(z) + a(z)^{-1}}{2}
   \end{pmatrix} \left(\frac{W(z)}{D(z)} \right)^{\sigma_3},
\end{equation}
where $W/D$ has an expansion near $1$ given in Lemma
\ref{BehaviorWoverD}. The functions $a(z) = \frac{(z-1)^{1/4}}{(z+1)^{1/4}}$
and $\varphi(z) = z + (z^2-1)^{1/2}$ that
appear in (\ref{BehaviorNWsigma3}) and Lemma \ref{BehaviorWoverD}
also have expansions near $1$.  We then obtain
\begin{equation}
    A^{(1)} = \label{DefinitieA1}
    \frac{4\alpha^2-1}{16}\ D_{\infty}^{\sigma_3}
    \begin{pmatrix}
        -1 & i \\
        i & 1
    \end{pmatrix}D_{\infty}^{-\sigma_3},
\end{equation}

Similarly, $B^{(1)}$ follows from (\ref{DefinitieDeltakTilde}),
Lemma \ref{BehaviorWTildeoverD}, and
the behavior of $\tilde{g}$, $N$, and $\tilde{W}$ near $-1$. The
result is that
\begin{equation}
    B^{(1)} = \label{DefinitieB1}
    \frac{4\beta^2-1}{16}\ D_{\infty}^{\sigma_3}
    \begin{pmatrix}
        1 & i\\
        i & -1
    \end{pmatrix}
    D_{\infty}^{-\sigma_3}.
\end{equation}
This completes the determination of $R_1$.

\begin{remark}
    To calculate $A^{(1)}$  we only need the leading terms in
    the expansions of $g$, $N$,  and $W$.  So in particular it is
    enough to use $W(z)/D(z) = 1 + O(|z-1|^{1/2})$ as $z \to 1$,
    see Lemma \ref{BehaviorWoverD}.
    Similarly, for $B^{(1)}$ it is enough to know
    $\tilde{W}(z)/D(z) = 1 + O(|z+1|^{1/2})$ as $z \to -1$,
    see Lemma \ref{BehaviorWTildeoverD}.
    This means that the analytic factor $h$ in the weight only
    plays a role through the constant $D_{\infty}$ in the
    determination of $R_1$.
\end{remark}

\subsubsection{Determination of $R_2$}

Expanding the jump relation $R_+ = R_-(I + \Delta_R)$ up to order
$1/n^2$ we find, cf. (\ref{SprongRk}),
\begin{equation} \label{JumpR2}
    R_{2,+}(s) - R_{2,-}(s) = R_{1,-}(s)\Delta_{1}(s)+\Delta_2(s),
    \qquad\mbox{for } s\in \partial U_{\delta}\cup \partial \tilde{U}_{\delta}.
\end{equation}
We view (\ref{JumpR2}) as a jump relation for $R_2$
on $\partial U_{\delta} \cup \partial \tilde{U}_{\delta}$ with additive jump
matrix
\[
    R_{1,-}(s)\Delta_{1}(s)+\Delta_2(s).
\]
Then (\ref{JumpR2}) and (\ref{AnalyticRk})--(\ref{AsympRk}) with $k=2$,
give an additive RHP for $R_2$.

Since $R_{1,-}(s)$ is the boundary value of the restriction of $R_1$
to the disks $U_{\delta}$ and $\tilde{U}_{\delta}$, and since
$\Delta_1$ and $\Delta_2$ are analytic in a neighborhood of $-1$ and $1$,
except for simple poles at $\pm 1$, see Lemma 8.2,
we have
\[
    R_{1}(z)\Delta_1(z)+\Delta_2(z)=\frac{A^{(2)}}{z-1}+O(1),\qquad\mbox{as }z\to 1,
\]
and
\[
    R_{1}(z)\Delta_1(z)+\Delta_2(z)=\frac{B^{(2)}}{z+1}+O(1),\qquad\mbox{as }z\to -1,
\]
for certain constant matrices $A^{(2)}$ and $B^{(2)}$.
As in the determination of $R_1$, we then find
\begin{equation} \label{DefinitieR2}
    R_2(z) =
    \left\{ \begin{array}{ll}
        \ds \frac{A^{(2)}}{z-1}+\frac{B^{(2)}}{z+1},
        & \mbox{for } z\in\mathbb{C}\setminus \left(\overline{U}_\delta\cup
            \overline{\tilde{U}}_\delta \right),\\[2ex]
        \ds \frac{A^{(2)}}{z-1}+\frac{B^{(2)}}{z+1}-R_{1}(z)\Delta_1(z)-\Delta_2(z),
        & \mbox{for } z\in  U_\delta\cup\tilde{U}_\delta.
    \end{array} \right.
\end{equation}

The calculation of $A^{(2)}$ is somewhat more complicated
than the determination of $A^{(1)}$.
It involves the two term expansion of the functions $g$, $W/D$, $a$,
and $\varphi$ near $z = 1$. In particular we need two terms in
the expansion in Lemma \ref{BehaviorWoverD}, which means that the
constant $c_0$ comes in. Using Maple we are able to
compute $A^{(2)}$, and the result is
\begin{equation}\label{DefinitieA2}
A^{(2)} = \frac{4\alpha^2-1}{256}
D_{\infty}^{\sigma_{3}}\begin{pmatrix} A_{2}(\alpha,\beta,c_0) &
iB_{2}(\alpha,\beta,c_0)\\ iC_{2}(\alpha,\beta,c_0)
& D_{2}(\alpha,\beta,c_0)
\end{pmatrix}D_{\infty}^{-\sigma_{3}},
\end{equation}
with $c_0$ given by formula (\ref{Definitiecn}) and
\begin{eqnarray*}
A_{2}(\alpha,\beta,c_0)&=& 8\alpha+8\beta+8c_0-4\beta^2+1,\\[1ex]
D_{2}(\alpha,\beta,c_0)&=& -8\alpha-8\beta-8c_0-4\beta^2+1,\\[1ex]
B_{2}(\alpha,\beta,c_0)&=&
    -8\alpha-8\beta-8c_0+4\alpha^2+ 4\beta^2-10,\\[1ex]
C_{2}(\alpha,\beta,c_0)&=&
    -8\alpha-8\beta-8c_0-4\alpha^2- 4\beta^2+10.
\end{eqnarray*}
Similarly, the result for $B^{(2)}$ is
\begin{equation}\label{DefinitieB2}
B^{(2)}=\frac{4\beta^2-1}{256}D_{\infty}^{\sigma_{3}}\begin{pmatrix}
-A_{2}(\beta,\alpha,d_0) & iB_{2}(\beta,\alpha,d_0)\\
iC_{2}(\beta,\alpha,d_0) & -D_{2}(\beta,\alpha,d_0)
\end{pmatrix}D_{\infty}^{-\sigma_{3}},
\end{equation}
with $d_0$ given by formula (\ref{Definitiedn}).

\begin{remark}
    To calculate $R_2$ we need two terms in the
    expansion of $W/D$ near $1$ and $\tilde{W}/D$ near $-1$, cf.
    Lemmas \ref{BehaviorWoverD} and
    \ref{BehaviorWTildeoverD}. That's why we see the influence of
    the analytic factor $h$ in the form  of the coefficients $c_0$ and $d_0$.
\end{remark}

\subsubsection{Result for general $R_k$}

In a similar (although increasingly more complicated) way
we can determine $R_k$ for every $k$. The calculations become
unreasonable to do by hand. Assisted by the computer algebra package
Maple, we find for $R_3$,
\begin{equation} \label{DefinitieR3}
R_3(z)=\frac{A^{(3)}}{z-1}+\frac{B^{(3)}}{z+1}+
\frac{C^{(3)}}{(z-1)^{2}}+\frac{D^{(3)}}{(z+1)^{2}},\qquad
\mbox{for } z\in\mathbb{C}\setminus
\left(\overline{U}_{\delta}\cup \overline{\tilde{U}}_{\delta}\right).
\end{equation}
where $A^{(3)}$ and $B^{(3)}$ are given by
\[
    A^{(3)} = \frac{4\alpha^2-1}{8192} D_{\infty}^{\sigma_{3}}
    \begin{pmatrix}
        A_{3}(\alpha,\beta,c_0,d_0) & iB_{3}(\alpha,\beta,c_0,d_0) \\
        iC_{3}(\alpha,\beta,c_0,d_0) & D_{3}(\alpha,\beta,c_0,d_0)
    \end{pmatrix}
    D_{\infty}^{-\sigma_{3}},
\]
\[
    B^{(3)}= \frac{4 \beta^2-1}{8192} D_{\infty}^{\sigma_{3}}
    \begin{pmatrix}
        -A_{3}(\beta,\alpha,d_{0},c_{0}) & iB_{3}(\beta,\alpha,d_0,c_0) \\
        iC_{3}(\beta,\alpha,d_0,c_0) & -D_{3}(\beta,\alpha,d_0,c_0)
    \end{pmatrix}
    D_{\infty}^{-\sigma_{3}}
\]
with $c_0$ and $d_0$ given by the formulas (\ref{Definitiecn}) and
(\ref{Definitiedn}) respectively, and
\begin{eqnarray}
    \nonumber
    A_{3}(\alpha,\beta,c_0,d_0)&=& 16(4\beta^2-1)(c_0+d_0+2\alpha+2\beta)-4\alpha^2(4\beta^2-1) \\[1ex] \nonumber
        && -128(\alpha^2+\beta^2+c_0^{2})-256c_0(\alpha+\beta)-256\alpha\beta+12\beta^2-16\beta^4-2,\\[2ex]
    \nonumber
    D_{3}(\alpha,\beta,c_0,d_0)&=& 16(4\beta^2-1)(c_0+d_0+2\alpha+2\beta)+4\alpha^2(4\beta^2-1) \\[1ex] \nonumber
        && +128(\alpha^2+\beta^2+c_0^{2})+256c_0(\alpha+\beta)+256\alpha\beta-12\beta^2+16\beta^4+2,
\end{eqnarray}
and
\begin{eqnarray} \nonumber
    B_{3}(\alpha,\beta,c_0,d_0)&=& 384(\alpha+\beta)^2+384c_0(c_0+2\alpha+2\beta)
        -388\alpha^2-252\beta^2+64\alpha^4 \\[1ex] \nonumber
        && +48\beta^4+144\alpha^2\beta^2+528-384(\alpha+\beta)(\alpha^2+\beta^2)+960(\alpha+\beta) \\[1ex] \nonumber
        && -192\beta^2(c_0+d_0)-384c_0\alpha^2+912c_0+48d_0, \\[2ex]
    \nonumber
    C_{3}(\alpha,\beta,c_0,d_0)&=& 384(\alpha+\beta)^2+384c_0(c_0+2\alpha+2\beta)
        -388\alpha^2-252\beta^2+64\alpha^4 \\[1ex] \nonumber
        && +48\beta^4+144\alpha^2\beta^2+528+384(\alpha+\beta)(\alpha^2+\beta^2)-960(\alpha+\beta) \\[1ex] \nonumber
        && +192\beta^2(c_0+d_0)+384c_0\alpha^2-912c_0-48d_0.
\end{eqnarray}
For $C^{(3)}$ and $D^{(3)}$, we find
\begin{eqnarray}
    \nonumber
    C^{(3)} &= & \frac{(4\alpha^2-1)(4\alpha^2-9)(4\alpha^2-25)}{12288}
    \begin{pmatrix}
        -1 & i \\
        i & 1
    \end{pmatrix}, \\[1ex]
    \nonumber
    D^{(3)} & = & \frac{(4\beta^2-1)(4\beta^2-9)(4\beta^2-25)}{12288}
    \begin{pmatrix}
        -1 & -i \\
        -i & 1
    \end{pmatrix}.
\end{eqnarray}

For general $k$, we get that
$R_k(z)$ in the region $\mathbb C \setminus
\left(\overline{U}_{\delta}\cup \overline{\tilde{U}}_{\delta}\right)$
is a rational function with poles at $1$ and $-1$. The residues in
$1$ and $-1$ are denoted by $A^{(k)}$ and $B^{(k)}$ respectively. Then
\begin{equation} \label{DefinitieRk}
    R_k(z) = \frac{A^{(k)}}{z-1} + \frac{B^{(k)}}{z+1} +
        O\left(\frac{1}{z^2}\right), \qquad \mbox{as $z \to \infty$.}
\end{equation}

    \section{Proof of Theorems \ref{theorempin}, \ref{theoremgamman}, and \ref{theoremanbn}}
\setcounter{equation}{0}

In this section we go back to the orthogonal
polynomials and prove the Theorems  \ref{theorempin}, \ref{theoremgamman},
and \ref{theoremanbn}. By tracing back the steps
\[
    Y \mapsto T \mapsto S \mapsto R
\]
and using the asymptotic expansion for $R$ obtained in the
previous section, we obtain an asymptotic expansion for $Y$. We
will restrict ourselves to the asymptotics in the exterior region
$\mathbb C \setminus [-1,1]$. The asymptotics on the interval
$(-1,1)$ and near the endpoints will be considered in the next
section. Let $K$ be a closed set in
$\overline{\mathbb{C}}\setminus[-1,1]$. Then we can choose the
lens and the regions $U_{\delta}$, $\tilde{U}_{\delta}$ disjoint
from $K$. From (\ref{YinT}), (\ref{SinT}), and (\ref{DefRBuiten}),
we then obtain
\begin{equation} \label{YinRN}
    Y(z) = 2^{-n\sigma_3} R(z) N(z) \varphi(z)^{n\sigma_3}, \qquad
    \mbox{ for } z \in K.
\end{equation}

\subsection{Asymptotics for $\pi_n$}

We start with the proof of Theorem \ref{theorempin} concerning the
asymptotics of the monic orthogonal polynomials $\pi_n(z)$ for $z \in K$.

\begin{varproof} \textbf{of Theorem \ref{theorempin}.}
Since $\pi_n = Y_{11}$ we obtain from (\ref{YinRN})
\begin{equation}
    \pi_{n}(z) = 2^{-n}\varphi(z)^{n}\left(R_{11}(z)N_{11}(z)+R_{12}(z)N_{21}(z)\right).
\end{equation}
If we use (\ref{RHPNsolution}) and (\ref{AsympOntR}),
we get uniformly for $z$ in $K$,
\begin{eqnarray*}
    \frac{2^n \pi_{n}(z)}{\varphi(z)^{n}} & \sim &
    \frac{D_{\infty}}{D(z)}\frac{a(z) + a(z)^{-1}}{2}
        \left[1 + \sum_{k=1}^{\infty} \frac{(R_k)_{11}(z)}{n^k} \right]
    \\[1ex]
    \nonumber
    &&
        \qquad +\, \frac{1}{D_{\infty}D(z)} \frac{a(z)-a(z)^{-1}}{-2i}
        \left[ \sum_{k=1}^{\infty} \frac{(R_k)_{12}(z)}{n^k} \right].
\end{eqnarray*}
Since
\[ \frac{a(z) + a(z)^{-1}}{2} = \frac{\varphi(z)^{1/2}}{\sqrt{2} (z^2-1)^{1/4}}, \]
we obtain Theorem \ref{theorempin} with
\begin{equation} \label{DefinitiePhik}
    \Pi_k(z) = (R_k)_{11}(z) + \frac{i}{D_{\infty}^2} \frac{a(z)-a(z)^{-1}}{a(z) + a(z)^{-1}}
    (R_k)_{12}(z).
\end{equation}
Using (\ref{DefinitieR1}), (\ref{DefinitieA1}), and
(\ref{DefinitieB1}) and simplifying the expressions, we see that
$\Pi_1$ is given by (\ref{DefinitiePi1}). From
(\ref{DefinitieR2}), (\ref{DefinitieA2}), and (\ref{DefinitieB2}),
we obtain the formula (\ref{DefinitiePi2}) for $\Pi_2$.
\end{varproof}

Note that  the factor $D_{\infty}^{-2}$ in the formula
(\ref{DefinitiePhik}) is cancelled by the factor $D_{\infty}^2$
which is included in the $(1,2)$ entry of $R_k$
(see formulas (\ref{DefinitieR1}), (\ref{DefinitieA1}), (\ref{DefinitieB1})
for the case $k=1$), so that $D_{\infty}$ does not appear in the final expression
for $\Pi_k$. This is due to the fact that we are dealing with monic
polynomials in Theorem \ref{theorempin}. The functions $\Pi_k(z)$ are rational functions of $\varphi(z)$
with coefficients that are polynomial expressions in $\alpha$, $\beta$, and
the coefficients $c_k$ and $d_k$.

\subsection{Asymptotics for $\gamma_n$}

We next consider the asymptotic behavior of the leading
coefficient $\gamma_n$ of the orthonormal polynomial.

\begin{varproof} \textbf{of Theorem \ref{theoremgamman}.}
To emphasize the $n$-dependence, let us use $Y^{(n)}$
to denote the solution (\ref{RHPYsolution}) of the RHP for $Y$.
Then we see from (\ref{RHPYsolution}) that
\[
    Y_{21}^{(n+1)}(z)= -2\pi i\gamma_n^2\pi_n(z).
\]
Since $\pi_n$ is a monic polynomial of degree $n$, we therefore have that
\[
    \gamma_{n}^{2}= - \frac{1}{2\pi i} \lim_{z\to\infty}z^{-n}Y_{21}^{(n+1)}(z).
\]
This expression appeared earlier in \cite[Theorem 3.1]{DKMVZ2} and
\cite[Section 3.2]{Deift}.
Using (\ref{YinRN}) (with $n$ replaced by $n+1$)
we find that for $z$ in the exterior region,
\begin{equation}\label{Verg1AsympGamman}
    \gamma_{n}^{2}= -\frac{1}{2\pi i}
    \lim_{z\to\infty} \left[2^{n+1}\left(\frac{\varphi(z)}{z}\right)^{n+1}
    \left(zR^{(n+1)}_{21}(z)N_{11}(z)+zR^{(n+1)}_{22}(z)N_{21}(z)\right)\right],
\end{equation}
where $R^{(n+1)}$ denotes the $R$ matrix corresponding to $n+1$.
Recall that $N$ does not depend on $n$.

We note that
\begin{equation}\label{Verg2AsympGamman}
    \frac{\varphi(z)}{z} =2 +O\left( \frac{1}{z^{2}} \right), \qquad\mbox{as }z\to\infty,
\end{equation}
and
\begin{equation}\label{Verg3AsympGamman}
    N_{11}(z) = 1 + O\left(\frac{1}{z}\right), \qquad
    R^{(n+1)}_{22}(z) = 1 + O\left(\frac{1}{z}\right), \qquad \mbox{as $z \to \infty$.}
\end{equation}
From (\ref{RHPNsolution}) it follows that
\begin{equation}
    N_{21}(z) = D_{\infty}^{-1} D(z)^{-1} \frac{a(z) - a(z)^{-1}}{-2i}
    \label{Verg4AsympGamman}
    = \frac{1}{2i D_{\infty}^2 z} + O\left(\frac{1}{z^2}\right),
    \quad \mbox{as } z \to \infty,
\end{equation}
and from (\ref{AsympOntR})
\[
    R_{21}^{(n+1)}(z) \sim \sum_{k=1}^{\infty}
    \frac{(R_k)_{21}(z)}{(n+1)^k}, \qquad \mbox{as } n\to\infty.
\]
Using this and Lemma \ref{LemmaAsympR} and (\ref{DefinitieRk}) we find
\begin{equation}\label{Verg5AsympGamman}
    \lim_{z \to \infty} z R_{21}(z)
    \sim \sum_{k=1}^{\infty} \frac{A_{21}^{(k)} +
    B_{21}^{(k)}}{(n+1)^k}, \qquad \mbox{as }n\to\infty.
\end{equation}
Combining (\ref{Verg1AsympGamman})--(\ref{Verg5AsympGamman}), we
obtain that $\gamma_n^2$ has an asymptotic expansion
\[
    \gamma_{n}^{2} \sim
    \frac{2^{2n}}{\pi D_{\infty}^2}\left(1+2i D_{\infty}^2
    \sum_{k=1}^{\infty} \frac{A_{21}^{(k)} + B_{21}^{(k)}}{(n+1)^k}\right),
    \qquad \mbox{as } n \to \infty.
\]
Taking square roots, and rearranging terms we obtain the asymptotic
expansion (\ref{AsympOntGamma}) for $\gamma_n$.
Using the exact expressions for $A^{(1)}$, $B^{(1)}$, $A^{(2)}$,
and $B^{(2)}$ given in (\ref{DefinitieA1}), (\ref{DefinitieB1}),
(\ref{DefinitieA2}), and (\ref{DefinitieB2}), we obtain the explicit expressions
(\ref{DefinitieGamma1}) and (\ref{DefinitieGamma2}) for $\Gamma_1$ and $\Gamma_2$.
\end{varproof}

\subsection{Asymptotics for the recurrence coefficients}
\setcounter{equation}{0}

Recall the three term recurrence relation
\begin{equation}\label{RecursieRelatiePi}
    \pi_{n+1}(z)=(z-b_{n})\pi_{n}(z)-a_{n}^{2}\pi_{n-1}(z)
\end{equation}
satisfied by the monic orthogonal polynomials.
In this section we give the proof of Theorem \ref{theoremanbn}
concerning the asymptotic behavior for
the recurrence coefficients $a_{n}$ and $b_{n}$.

\begin{varproof} \textbf{of Theorem \ref{theoremanbn}.}
First, we
determine formulas for $a_{n}$ and $b_{n}$ in terms of $Y$.
We use an argument taken from \cite[Section 3.3]{Deift}.
Let $Y^{(n)}$ and $Y^{(n+1)}$ be the solutions of the RHP for
$Y$ corresponding to $n$ and $n+1$, respectively. Since $Y^{(n+1)}$ and
$Y^{(n)}$ are analytic in
$\mathbb{C}\setminus[-1,1]$ and satisfy the same jump relation on
the interval $(-1,1)$, it follows that
\begin{equation} \label{DefHn}
    H^{(n)} : = Y^{(n+1)}\left(Y^{(n)}\right)^{-1}
\end{equation}
is analytic in $\mathbb{C}\setminus\{-1,1\}$. In the same way as was done
in the proof of Lemma \ref{lemmaRHPY}, we can prove that $H^{(n)}$ has removable
singularities at the points $\pm 1$. Therefore $H^{(n)}$ is entire.
From the asymptotic behavior (\ref{RHPYc}), the definition (\ref{DefHn}),
and the fact that $\det Y^{(n)} \equiv 1$, it follows that
\[
    H^{(n)}(z)=
    \begin{pmatrix}
        z+O(1) & O(1) \\
        O(1) & O\left(z^{-1}\right)
    \end{pmatrix}, \qquad\mbox{as $z\to\infty$.}
\]
Using Liouville's theorem we then find
that
\begin{equation} \label{Definitiernsntn}
    H^{(n)}(z) = \begin{pmatrix} z-r_n & s_n \\ t_n & 0 \end{pmatrix}
\end{equation}
for certain constants $r_n$, $s_n$, and $t_n$.
Thus
\begin{equation}\label{RecursieRelatieY}
Y^{(n+1)}(z)=
\begin{pmatrix}
    z-r_{n} & s_{n}\\
    t_{n} & 0
\end{pmatrix}Y^{(n)}(z).
\end{equation}
The $(1,1)$ entry of this matrix relation gives
\[
    \pi_{n+1}(z)=(z-r_{n})\pi_{n}(z)-2\pi i\gamma_{n-1}^{2}s_{n}\pi_{n-1}(z),
\]
which is the recurrence relation  for the orthogonal polynomials.
Comparing with (\ref{RecursieRelatiePi}), we
see that $b_n = r_n$ and $a_n^2 = 2\pi i\gamma_{n-1}^2 s_n$,
The $(2,1)$ entry of
(\ref{RecursieRelatieY}) shows that $-2\pi i\gamma_{n}^{2}
=t_{n}$. Thus
\begin{equation} \label{bnaninHn}
    b_n = r_n, \qquad a_n^2 = -s_n t_{n-1}.
\end{equation}

Using (\ref{DefHn}), (\ref{Definitiernsntn}), (\ref{bnaninHn})
and the fact that $\det Y^{(n)}(z) = 1$, we find that
\[
    b_n=  z-Y_{11}^{(n+1)}(z)Y_{22}^{(n)}(z)+Y_{12}^{(n+1)}(z)Y_{21}^{(n)}(z),
\]
and
\[
    a_{n}^{2}=  \left(Y_{11}^{(n+1)}(z) Y_{12}^{(n)}(z) -Y_{12}^{(n+1)}(z) Y_{11}^{(n)}(z)\right)
    \left(Y_{21}^{(n)}(z)Y_{22}^{(n-1)}(z)-Y_{22}^{(n)}(z) Y_{21}^{(n-1)}(z) \right),
\]
for every $z \in \mathbb C \setminus [-1,1]$. From the asymptotic
behavior (\ref{RHPYc}) of $Y$ we then find, upon letting $z \to \infty$,
\begin{equation}\label{bnFormule1}
    b_{n} = \lim_{z\to\infty}\left(z-Y_{11}^{(n+1)}(z)Y_{22}^{(n)}(z)\right),
\end{equation}
and
\begin{equation}\label{anFormule1}
    a_{n}^{2}=\lim_{z\to\infty}z^{2}Y_{12}^{(n)}(z)Y_{21}^{(n)}(z).
\end{equation}

We use (\ref{bnFormule1}) and (\ref{anFormule1}) to
express the recurrence coefficients in terms of $R$. From (\ref{YinRN}) and (\ref{anFormule1}), we get for the
recurrence coefficient $a_{n}$
\[ a_{n}^{2}=\lim_{z\to\infty}z^{2}\left(R_{21}^{(n)} (z)N_{11}(z)+R_{22}^{(n)}(z)N_{21}(z)\right)
        \left(R_{11}^{(n)}(z)N_{12}(z)+R_{12}^{(n)}(z)N_{22}(z)\right).
\]
From the behavior of $R^{(n)}$ and $N$ at infinity (in
particular (\ref{RHPNsolution})) one then easily sees that
\begin{equation}\label{anFormule2}
    a_{n}^{2}=\lim_{z\to\infty}\left(zR_{21}^{(n)}(z)+\frac{1}{2i D_{\infty}^2}\right)
        \left(-\frac{D_{\infty}^2}{2i}+zR_{12}^{(n)}(z)\right).
\end{equation}
For $b_n$ we find in the same way
\begin{equation}\label{bnFormule2}
    b_n = \lim_{z\to\infty}\left(z- z R_{11}^{(n+1)}(z)R_{22}^{(n)}(z)\right).
\end{equation}

In the asymptotic expansion for $R$, we may take the limit $z \to
\infty$, by Lemma \ref{LemmaAsympR}. We find
\begin{equation} \label{anFormule3}
a_n^2 \sim \left(\frac{1}{2i D_{\infty}^2} +
    \sum_{k=1}^{\infty} \frac{A_{21}^{(k)} + B_{21}^{(k)}}{n^k}\right)
    \left(- \frac{D_{\infty}^2}{2i} +
    \sum_{k=1}^{\infty} \frac{A_{12}^{(k)} + B_{12}^{(k)}}{n^k} \right),
\end{equation}
and
\begin{eqnarray} \nonumber
b_n & \sim & \lim_{z \to \infty} \left[ z - z
    \left(1 + \sum_{k=1}^{\infty} \frac{(R_{k})_{11}(z)}{(n+1)^k} \right)
    \left(1 + \sum_{l=1}^{\infty} \frac{(R_l)_{22}(z)}{n^l} \right) \right] \\[1ex]
    & = & \nonumber
    \lim_{z \to \infty} -z \left( \sum_{k=1}^{\infty} \frac{(R_k)_{11}(z)}{(n+1)^k}
    + \sum_{l=1}^{\infty} \frac{(R_l)_{22}(z)}{n^l} \right) \\[1ex]
    & = & - \sum_{k=1}^{\infty} \left( \frac{A^{(k)}_{11} + B^{(k)}_{11}}{(n+1)^k}
        + \frac{A^{(k)}_{22} + B^{(k)}_{22}}{n^k} \right).
        \label{bnFormule3}
\end{eqnarray}
Expanding (\ref{anFormule3}), we find a complete asymptotic expansion
for $a_n^2$ in powers of $1/n$. The constant term is $1/4$, and
the term with $1/n$ has coefficient
\[ \frac{1}{2i} \left( D_{\infty}^{-2} (A_{12}^{(1)} + B_{12}^{(1)})
    - D_{\infty}^2 (A_{21}^{(1)} + B_{21}^{(1)}) \right), \]
which by the formulas (\ref{DefinitieA1}) and (\ref{DefinitieB1})
for $A^{(1)}$ and $B^{(1)}$ is equal to zero.
Thus
\[ a_n  = \frac{1}{2} + O\left(\frac{1}{n^2}\right). \]
The full expansion of $a_n$ follows from (\ref{anFormule3})
and gives us the expansion (\ref{asympan}) with explicitly computable
coefficients $A_k$.

For $b_n$, we find from (\ref{bnFormule3}) that the $1/n$ term
in the expansion  has the coefficient
\[ - \left(A_{11}^{(1)} + B_{11}^{(1)} + A_{22}^{(1)} + B_{22}^{(1)} \right), \]
which is zero by (\ref{DefinitieA1}) and (\ref{DefinitieB1}).
Thus
\[ b_n = O\left(\frac{1}{n^2}\right). \]
The full expansion (\ref{asympbn}) of $b_n$ follows by rearranging terms in (\ref{bnFormule3}).
\end{varproof}

\begin{remark}
Note that (\ref{anFormule1}) can also be obtained directly from
(\ref{RHPYsolution}), since
\[ Y_{21}(z) = -2\pi i \gamma_{n-1}^2 z^{n-1} \left(1 + O\left(\frac{1}{z}\right)\right) \]
and
\[ Y_{12}(z) = \frac{1}{2\pi i} \int_{-1}^1 \frac{\pi_n(x) w(x)}{x-z} dx
    = - \frac{1}{2\pi i} \gamma_n^{-2} z^{-n-1} \left(1+ O\left(\frac{1}{z}\right)\right). \]
Then (\ref{anFormule1}) follows easily, because
\[ a_n = \frac{\gamma_{n-1}}{\gamma_n}. \]

The formula (\ref{anFormule1}) is due to Deift et al.\ \cite[Theorem 3.1]{DKMVZ2},
see also \cite[Section 3.2]{Deift}.
However the formula (\ref{bnFormule1}) is different from the one given in
these references, which is due to the fact that the orthogonality weight is
varying with $n$ in \cite{DKMVZ2}.
\end{remark}

\section{Proof of Theorems \ref{theorem: monic orthogonal polynomial bulk: MJW}
    and \ref{theorem: monic orthogonal polynomial endpoint: MJW}}

In this section we determine the asymptotics for the monic
orthogonal polynomials $\pi_n$ on the interval of orthogonality
$[-1,1]$. We use the behavior of the first column of $Y$ (which
contains the orthogonal polynomials) on $[-1,1]$, which has
already been determined, see \cite{KuijlaarsVanlessen}. In
\cite{KuijlaarsVanlessen} this behavior has been used to derive
universal behavior for eigenvalue correlations, associated to the
modified Jacobi unitary ensemble, in the bulk and at the hard edge
of the spectrum.

\medskip

The behavior of the first column of $Y$ in the bulk of the
interval $[-1,1]$ is given by \cite[Lemma
3.1]{KuijlaarsVanlessen}. Namely,
\begin{eqnarray} \label{lemma: asymptotics bulk: MJW: equation}
    \nonumber
    \lefteqn{
    \begin{pmatrix}
        Y_{11}(x) \\
        Y_{21}(x)
    \end{pmatrix}
    =
    \frac{1}{\sqrt{w(x)}}\, 2^{-n\sigma_3}
    R(x)D_\infty^{\sigma_3}} \\[1ex]
    &&
    \qquad\qquad\qquad\times\,
        \begin{pmatrix}
            \frac{a_+(x)+a_+(x)^{-1}}{2} &
            \frac{a_+(x)-a_+(x)^{-1}}{2i} \\[1ex]
            \frac{a_+(x)-a_+(x)^{-1}}{-2i} & \frac{a_+(x)+a_+(x)^{-1}}{2}
        \end{pmatrix} e^{i(n\arccos x+\psi(x))\sigma_3}
    \begin{pmatrix}
        1 \\
        1
    \end{pmatrix},
\end{eqnarray}
for $x\in(-1+\delta,1-\delta)$, where $\psi$ is given by
(\ref{definition: psi: MJW}).

\begin{varproof}\textbf{of Theorem \ref{theorem: monic orthogonal polynomial bulk: MJW}.}
    Let $K$ be a compact subset of $(-1,1)$. We can choose the
    radius $\delta$ of the disks $U_\delta$ and $\tilde U_\delta$ so small that
    $K\subset (-1+\delta,1-\delta)$. From (\ref{lemma: asymptotics bulk: MJW:
    equation}) and the facts that
    \begin{equation}\label{formulas a}
        \frac{a(z)+a(z)^{-1}}{2}=\frac{\varphi(z)^{1/2}}{\sqrt 2
        (z^2-1)^{1/4}}, \quad \mbox{and} \quad
        \frac{a(z)-a(z)^{-1}}{2i}=i\frac{\varphi(z)^{-1/2}}{\sqrt 2
        (z^2-1)^{1/4}},
    \end{equation}
    we obtain
    \begin{eqnarray}
        \nonumber
        \lefteqn{
        \begin{pmatrix}
            Y_{11}(x) \\
            Y_{21}(x)
        \end{pmatrix}
        =
        \frac{1}{\sqrt{w(x)}}\frac{e^{-\frac{\pi i}{4}}}{\sqrt 2
        (1-x^2)^{1/4}}
        2^{-n\sigma_3}R(x)D_\infty^{\sigma_3}
        } \\[2ex]
        \nonumber
        &&
            \qquad\qquad\qquad \times\,
            \begin{pmatrix}
                \varphi_+(x)^{1/2} & i \varphi_+(x)^{-1/2} \\
                -i \varphi_+(x)^{-1/2} & \varphi_+(x)^{1/2}
            \end{pmatrix}
            e^{i(n\arccos x+\psi(x))\sigma_3}
            \begin{pmatrix}
                1 \\ 1
            \end{pmatrix},
    \end{eqnarray}
    for $x\in K$. Since $\varphi_+(x)=\exp(i\arccos x)$ we then
    find after an easy calculation
    \begin{eqnarray}
        \nonumber
        \lefteqn{
        \begin{pmatrix}
            Y_{11}(x) \\
            Y_{21}(x)
        \end{pmatrix}
        =\frac{1}{\sqrt{w(x)}}\frac{\sqrt 2}{
        (1-x^2)^{1/4}}
        2^{-n\sigma_3}R(x)D_\infty^{\sigma_3}
        } \\[2ex]
        &&
            \qquad\qquad\qquad\times\,
            \begin{pmatrix}
                \cos\left((n+\frac{1}{2})\arccos
                x+\psi(x)-\frac{\pi}{4}\right) \\[1ex]
                -i \cos\left((n-\frac{1}{2})\arccos
                x+\psi(x)-\frac{\pi}{4}\right)
            \end{pmatrix}.
    \end{eqnarray}
    Since $Y_{11}=\pi_n$, and since $R_{11}(x)=1+O(1/n)$ and
    $R_{12}(x)=O(1/n)$ as $n\to\infty$ with error terms that hold uniformly for $x\in
    K$, and that have a full asymptotic expansion in powers of $1/n$,
    the theorem is then proven.
\end{varproof}

Next, the behavior of the first column of $Y$ near the endpoint 1
is given by \cite[Lemma 3.3]{KuijlaarsVanlessen}. Namely,
\begin{equation}\label{lemma: asymptotics edge: MJW: equation}
        \begin{pmatrix}
            Y_{11}(x) \\
            Y_{21}(x)
        \end{pmatrix}
        =
        \sqrt{\frac{2\pi n}{w(x)}}\, 2^{-n\sigma_3}M_+(x)
        \begin{pmatrix}
            J_\alpha(n\arccos x) \\
            \frac{i}{2}\arccos x J'_\alpha(n\arccos x)
        \end{pmatrix},
    \end{equation}
for $x\in (1-\delta,1)$, with $M(z)$ given by
    \begin{equation}\label{definition: M: MJW}
        M(z)
        =
        R(z)N(z)W(z)^{\sigma_3}\frac{1}{\sqrt{2}}
        \begin{pmatrix}
            1 & -i \\
            -i & 1
        \end{pmatrix}
        f(z)^{\sigma_3/4}.
    \end{equation}

\begin{varproof}\textbf{of Theorem \ref{theorem: monic orthogonal polynomial endpoint: MJW}.}
    We need to determine a convenient expression for $M_+(x)$ and insert this into
    (\ref{lemma: asymptotics edge: MJW: equation}).
    From (\ref{RHPNsolution}), (\ref{formulas a})
    and (\ref{definition: M: MJW}) we have
    \begin{eqnarray}
        \nonumber
        \lefteqn{
            M_+(x) = \frac{e^{-\frac{\pi i}{4}}}{2 (1-x^2)^{1/4}} R(x)D_\infty^{\sigma_3}
            \begin{pmatrix}
                \varphi_+(x)^{1/2} & i\varphi_+(x)^{-1/2} \\
                -i\varphi_+(x)^{-1/2} & \varphi_+(x)^{1/2}
            \end{pmatrix}
        } \\[2ex]
        \label{proof: theorem endpoint: MJW: eq1}
        && \qquad\qquad\qquad \times\, \left(\frac{W_+(x)}{D_+(x)}\right)^{\sigma_3}
            \begin{pmatrix}
                1 & -i \\
                -i & 1
            \end{pmatrix} f_+(x)^{\sigma_3/4}.
    \end{eqnarray}
    It may be shown that $D_+(x)=\sqrt{w(x)}\exp(-i\psi(x))$, cf.\ \cite[Remark
    3.2]{KuijlaarsVanlessen}. We then have by (\ref{W^2(z)})
    \[
        \frac{W_+(x)}{D_+(x)}=e^{
        i(\psi(x)+\frac{\alpha\pi}{2})}.
    \]
    Since $\varphi_+(x)=\exp(i\arccos x)$ we then obtain from (\ref{definition: zeta12:
    MJW}) and (\ref{proof: theorem endpoint: MJW: eq1})
    \begin{eqnarray}
        \nonumber
        M_+(x)
        &=&
            \frac{e^{-\frac{\pi i}{4}}}{2 (1-x^2)^{1/4}} R(x)D_\infty^{\sigma_3}
            \begin{pmatrix}
                e^{i\zeta_1(x)} & i e^{-i\zeta_1(x)} \\
                -i e^{i\zeta_2(x)} & e^{-i\zeta_2(x)}
            \end{pmatrix}
            \begin{pmatrix}
                1 & -i \\
                -i & 1
            \end{pmatrix} f_+(x)^{\sigma_3/4} \\[2ex]
        \label{proof: theorem endpoint: MJW: eq2}
        &=&
            \frac{e^{-\frac{\pi i}{4}}}{(1-x^2)^{1/4}}R(x)D_\infty^{\sigma_3}
            \begin{pmatrix}
                \cos\zeta_1(x) & \sin\zeta_1(x) \\
                -i\cos\zeta_2(x) & -i\sin\zeta_2(x)
            \end{pmatrix} f_+(x)^{\sigma_3/4}.
    \end{eqnarray}
    Inserting the fact that
    \[
        f_+(x)^{1/4}=\frac{e^{\pi i/4}}{\sqrt 2} (\arccos
        x)^{1/2}
    \]
    into (\ref{proof: theorem endpoint: MJW: eq2}) we find
    \begin{eqnarray}
        \nonumber
        \lefteqn{
            M_+(x)
            \begin{pmatrix}
                J_\alpha(n\arccos x) \\
                \frac{i}{2}\arccos x J'_\alpha(n\arccos x)
            \end{pmatrix}
        } \\[2ex]
        \label{proof: theorem endpoint: MJW: eq3}
        &=&
            \frac{(\arccos x)^{1/2}}{\sqrt 2 (1-x^2)^{1/4}}
            R(x)D_\infty^{\sigma_3}
            \begin{pmatrix}
                \cos\zeta_1(x) & \sin\zeta_1(x) \\
                -i\cos\zeta_2(x) & -i\sin\zeta_2(x)
            \end{pmatrix}
            \begin{pmatrix}
                J_\alpha(n\arccos x) \\
                J'_\alpha(n\arccos x)
            \end{pmatrix}.
    \end{eqnarray}
    If we now insert (\ref{proof: theorem endpoint: MJW: eq3}) into (\ref{lemma: asymptotics edge: MJW: equation}),
    and use the facts
    that $Y_{11}=\pi_n$,
    $R_{11}(x)=1+O(1/n)$ and $R_{12}(x)=O(1/n)$ as $n\to\infty$,
    we obtain (\ref{theorem: monic orthogonal polynomial endpoint: MJW: equation}).
    Note that the error terms hold uniformly for
    $x\in(1-\delta,1)$ and have a full asymptotic expansion in
    powers of $1/n$. This proves the theorem.
\end{varproof}

\end{document}